\documentclass[11pt,leqno]{article}
\usepackage{amsmath,amsthm,amsfonts,amssymb}

\newcommand{\CC}{{\mathbb C}}

\newcommand{\EE }{{\mathbb E}}

\newcommand{\GG}{{\mathbb G}}

\newcommand{\cA}{{\mathcal A}}
\newcommand{\cB}{{\mathcal B}}
\newcommand{\cC}{{\mathcal C}}

\newcommand{\cE}{{\mathcal E}}
\newcommand{\cF}{{\mathcal F}}
\newcommand{\cG}{{\mathcal G}}

\newcommand{\cI}{{\mathcal I}}
\newcommand{\cJ}{{\mathcal J}}

\newcommand{\cL}{{\mathcal L}}

\newcommand{\cN}{{\mathcal N}}
\newcommand{\cO}{{\mathcal O}}

\newcommand{\cQ}{{\mathcal Q}}
\newcommand{\cR}{{\mathcal R}}
\newcommand{\cS}{{\mathcal S}}
\newcommand{\cT}{{\mathcal T}}
\newcommand{\cU}{{\mathcal U}}

\newcommand{\cW}{{\mathcal W}}

\newcommand{\cZ}{\mathcal{Z}}

\newcommand{\oA}{\overline{A}}
\newcommand{\oB}{\overline{B}}

\newcommand{\oS}{\overline{S}}

\newcommand{\chcL}{{\check{\mathcal L}}}

\newcommand{\chc}{\check c}

\newcommand{\cbR}{\check{\textbf R}}

\newcommand{\rmR}{\mathrm{ R }}

\newcommand{\uh}{{\underline h}}

\newcommand{\up}{{\underline p}}
\newcommand{\uq}{{\underline q}}

\newcommand{\uw}{{\underline w}}
\newcommand{\uz}{{\underline z}}
\newcommand{\ua}{{\underline a}}
\newcommand{\us}{{\underline s}}

\newcommand{\doa}{\dot{a}}
\newcommand{\dop}{\dot{p}}
\newcommand{\dopsi}{\dot{\psi}}

\newcommand{\dov}{\dot{v}}
\newcommand{\dou}{\dot{u}}
\newcommand{\dow}{\dot{w}}
\newcommand{\dos}{\dot{s}}
\newcommand{\doh}{\dot{h}}

\newcommand{\bfR}{\textbf{R}}

\newcommand{\ibfR}{\check{\textit{\textbf R}}}

\newcommand{\Id}{{\rm Id }}

\newcommand{\D}{{\partial }}
\newcommand{\eps}{\varepsilon}

\newcommand{\tf}{\tilde{f}}

\newcommand{\tH}{\tilde{H}}

\newcommand{\tD}{\tilde{D}}

\newcommand{\tx}{\tilde{x}}

\newcommand{\ta}{\tilde{a}}

\newcommand{\te}{\tilde{e}}

\newcommand{\tGamma}{\tilde{\Gamma}}
\newcommand{\tPsi}{\tilde{\Psi}}

\newcommand{\bF}{\mathbb{F}}
\newcommand{\bR}{\mathbb{R}}

\newcommand{\bE}{\mathbb{E}}

\newcommand{\bS}{\mathbb{S}}
\newcommand{\bK}{\mathbb{K}}

\newcommand{\bD}{\mathbb{D}}
\newcommand{\bC}{\mathbb{C}}

\newcommand{\hzeta}{\hat{\zeta}}
\newcommand{\heta}{\hat{\eta}}
\newcommand{\hgamma}{\hat{\gamma}}

\newcommand{\htau}{\hat{\tau}}
\newcommand{\hGamma}{\hat{\Gamma}}

\newcommand{\ocA}{\overline{\mathcal{A}}}
\newcommand{\obR}{\overline{\mathbb{R}}}

\newcommand{\uV}{\underline{V}}

\newtheorem{thm}{Theorem}[section]
\newtheorem{prop}[thm]{Proposition}
\newtheorem{cor}[thm]{Corollary}
\newtheorem{lem}[thm]{Lemma}
\newtheorem{defn}[thm]{Definition}
\newtheorem{ass}[thm]{Assumption}

\newtheorem{ex}[thm]{Example}
\newtheorem{rem}[thm]{Remark}

\newtheorem{notation}[thm]{Notation}
\numberwithin{equation}{section}

 \title{     \Large\bf
Nonclassical multidimensional viscous and inviscid shocks}

\author{Olivier Gu\`es\footnote{Universit\'e de Provence, partially supported by European network HYKE,  HPRN-CT-2002-00282},
Guy M\'etivier\footnote{Universit\'e de Bordeaux, partially
supported by European network HYKE,  HPRN-CT-2002-00282.}, Mark
Williams\footnote{University of North Carolina, partially
supported  by NSF grants DMS-0070684 and DMS-0401252},
 Kevin Zumbrun\footnote{Indiana University, partially supported  by NSF grants  DMS-0070765 and DMS-0300487. }}

\date{June 27, 2006}

\begin{document}

\maketitle
\begin{abstract}
%existence and stability of curved multi-d shocks for
%nonclassical of possibly nonconservative systems.
%In process, extend inviscid theory of Mokrane, Coulombel,
%and one-dimensional inviscid theory of Xiao-biao Lin.
%
Extending our earlier work on Lax-type shocks of systems of
conservation laws, we establish existence and stability of curved
multidimensional shock fronts in the vanishing viscosity limit for
general Lax- or undercompressive-type shock waves of nonconservative
hyperbolic systems with parabolic regularization. The hyperbolic
equations may be of variable multiplicity and the parabolic
regularization may be of ``real'', or partially parabolic, type.  We
prove an existence result for inviscid nonconservative shocks that
extends to multidimensional shocks a one-dimensional result of X.
Lin proved by quite different methods.   In addition, we construct
families of smooth viscous shocks converging to a given inviscid
shock as viscosity goes to zero, thereby justifying the small
viscosity limit for multidimensional nonconservative shocks.

In our previous work on shocks we made use of conservative form
especially in parts of the low frequency analysis. Thus, most of the
new analysis of this paper is concentrated in this area.    By
adopting the more general nonconservative viewpoint, we are able to
shed new light on both the viscous and inviscid theory.   For
example, we can now provide a clearer geometric motivation for the
low frequency analysis in the viscous case.  Also, we show that one
may, in the treatment of inviscid stability of nonclassical and/or
nonconservative shocks, remove an apparently restrictive technical
assumption made by Mokrane and Coulombel  in their work on,
%CHANGED (one more try... ! -K)
%respectively, shock-like nonconservative boundary problems and
respectively, shock-type nonconservative boundary problems and
conservative undercompressive shocks.
%ENDCHANGED
Another advantage of the nonconservative perspective is that Lax and
undercompressive shocks can be treated by exactly the same analysis.

%TODO: Did Xiaobiao actually carry out the inviscid limit problem,
%or just work out the inviscid theory (I think only the latter).
%If the former, put in something like this:
%Our results generalize a one-dimensional analysis of X. Lin,
%carried out by quite different methods.
%
%TODO: something along these lines? (I say no, better to keep
%it short-K)
%the main new features of the analysis appear in the low frequency
%estimates for the linearized problem, specifically in determining
%reduced problems and appropriate modified Evans function.
 \end{abstract}

\tableofcontents

\section{Introduction}

%LF standard Evans implies LF modified Evans.  LF standard is the one
%for which we almost have [PZ].  LF modified Evans needed for good
%viscous estimate as in GMWZ3.

%New perspective on GMWZ3.

%Remove Coulombel hypothesis for inviscid theory.

%Don't assume continuity of decaying eigenspaces or constant
%multiplicities in Evans/Lopatinski analysis.

%Main differences with conservative case occur (only?) in LF analysis
%(check).

%In GMWZ4 used divergence form to get profile equation for $w_2$,
%where $w=(w_1,w_2)$. Now?   Schechter treats the case of degenerate
%$B$ matrices.

%No point in suppressing continuity assumption in small viscosity
%limit. This is ok for Evans analysis though.

%\section{Introduction}

\textbf{\quad\quad}In this paper we develop a theory of
multidimensional inviscid and viscous shocks without assuming
conservative form for the underlying inviscid system.   One
motivation, of course, is to provide a
 theory of multidimensional nonconservative shocks.    Another
motivation is less obvious:  by dropping the assumption of
conservative form it turns out that we are able to treat both Lax
shocks and undercompressive shocks by \emph{exactly the same
analysis}.   Under the appropriate spectral stability conditions,
the construction of curved inviscid shocks, and the proof that such
shocks can be obtained as vanishing viscosity limits of viscous
shocks, can be accomplished by an analysis that does not distinguish
between Lax and undercompressive shocks.\footnote{ We do not treat
overcompressive shocks, for which behavior in the inviscid limit is
more complicated. See \cite{F1, F2, FL, ZS, Z1} for further
discussion. }We regard classical (i.e., conservative Lax) shocks as
special cases of nonclassical shocks and conservative shocks as
special cases of nonconservative ones. Recall that systems in
nonconservative form arise in models of deformation of
%CHANGED(added Gavrilyuk ref.- it is to relaxation, but never mind...-K)
%elastic-plastic solids \cite{TC}, two-phase flow \cite{SW}, spray
elastic-plastic solids \cite{TC}, two-phase flow \cite{SW, SGR, Se},
spray
%ENDCHANGED
dynamics \cite{RS}, and other applications.

In our previous work on shocks \cite{GMWZ1}-\cite{GMWZ4} we made use
of conservative form especially in parts of the low frequency
analysis. Thus, most of the new analysis of this paper is
concentrated in this area. The nonconservative viewpoint provides a
clearer geometric motivation for the low frequency analysis and
leads, for example, to a simpler construction of high order
approximate solutions.  In passing, we simplify also the
nonclassical inviscid theory, removing a technical assumption of
\cite{Mo, Cou}; see Remark \ref{mormk}. As we did for viscous
boundary problems involving a given fixed boundary in \cite{GMWZ5},
we allow the underlying inviscid system to have characteristics of
variable multiplicity; thus, our analysis is applicable, for
example, to (Lax or undercompressive) MHD shocks.

%CHANGED -0.
When $d=1$, let us recall the result of [BB] which applies to the
vanishing viscosity approach of 1D quasilinear strictly hyperbolic
systems (including non conservative ones), in a context which
extends the usual one of ``entropy weak solutions'' of the
conservative case.
%ENDCHANGED
%CHANGED(added-K):
This result is stronger than the ones derived here in that it is for
``unprepared data'' (i.e., it states that viscous solutions for {\it
any} nearby data converge to the inviscid solution) whereas ours are
for ``prepared data'' (viscous solutions for {\it some} nearby data
converge to the inviscid solution), is {\it global} in time (hence
accommodates interaction of shock fronts), and is framed in a weaker
norm (B.V. vs. $H^s$). On the other hand, the ``Glimm-type''
analysis
%CHANGED(added)
of \cite{BB}
%ENDCHANGED
based on approximate decoupling
of scalar modes uses strongly the assumptions of small
variation of the background solution (hence small shock amplitude),
Laplacian viscosity, strict hyperbolicity, and a single space dimension,
whereas one of the main motivations for our approach is to
remove such restrictions.
%ENDCHANGED

The approach we have used in \cite{GMWZ2,GMWZ3,GMWZ4} to construct
families of smooth, exact, viscous shocks converging to a given
curved inviscid shock in the vanishing viscosity limit has four main
steps.

The first step is to construct viscous profiles for planar shocks;
this amounts to solving an ODE with prescribed endstates at
$\pm\infty$.  These exact solutions to the viscous problem  describe
the fast, shock layer, transition between two constant states.

The second step is to linearize the full, parabolic (or partially
parabolic) problem about a profile solution and define  appropriate
spectral stability determinants or Evans functions for this
linearized problem. Suitable nonvanishing conditions for these
determinants give necessary (and sometimes sufficient) conditions
for linear and nonlinear stability.  We define two sorts of Evans
functions, standard and modified; the first exhibits the usual
translational degeneracy at zero frequency, while the second is
typically nonvanishing at zero frequency.  Part of this step is to
clarify the connection between the two Evans functions (see Theorem
\ref{j3}, for example).

The third step is to compute the standard Evans function for a given
profile and check whether the nonvanishing conditions are actually
satisfied. This can be done numerically (\cite{HZ}, e.g.) and
sometimes even analytically (\cite{PZ,FS}).

The final step is to construct approximate, curved viscous shock
solutions and then prove the linear and nonlinear stability of those
solutions, assuming the Evans condition is satisfied.  This step
involves understanding the link between the viscous Evans function,
the Lopatinski determinant that governs stability of the inviscid
hyperbolic problem, and the transversality properties of profiles
(see Theorem \ref{K}, for example). It also involves the
construction of Kreiss-type symmetrizers and the proof of ``maximal"
estimates for the linearized parabolic problem.

Most of the work of this paper is concerned with the second and
fourth steps. For some of the problems we consider, step one has not
been done yet. For example, we are not aware of any viscous profile
constructions for nonconservative, undercompressive shocks (except
in some trivial cases). For such shocks the present paper reduces
the full nonlinear stability problem to the construction of profiles
and the verification of the Evans hypothesis.

In other problems step one has been done, but not step three. For
example, viscous profiles have been constructed for nonconservative,
Lax shocks (\cite{S,Sa}, for example), and for some conservative,
undercompressive shocks (\cite{AMPZ,AMPZ2, IMP, Sh,SSh1,SSh2}), but
the Evans hypothesis has not yet been carefully verified.   In these
cases our paper reduces the full nonlinear stability problem to the
verification of the Evans hypothesis.  In some problems, such as the
Navier-Stokes regularization of Euler shocks studied in
\cite{GMWZ4}, all four steps are now complete.

Finally, let us clarify the relation of section \ref{genRH} of this
paper to the first step.  If one \textbf{starts with} a single
transversal profile for a given planar shock, the generalized
Rankine-Hugoniot condition identifies ``nearby" shocks for which
profiles must also exist.

\subsection{Inviscid $\cC$-shocks}\label{inviscidC}
We begin by defining a notion of inviscid shock that is  more
general than the standard concept of nonconservative shock (used,
e.g., in \cite{Lin} and recalled below). Consider an $N\times N$
system on $\bR^{d+1}$
\begin{align}\label{i1}
\sum^d_{j=0}A_j(u)\partial_ju=0.
\end{align}
The system is said to be \emph{conservative} when $A_j(u)=f_j'(u)$
for $\bR^N$-valued functions $f_j$.  This requirement is dropped in
the following assumptions.

The Assumptions \ref{i1a} and \ref{i10} stated below are in force
throughout the paper.  The other Assumptions are made in a given
Theorem, Proposition, etc., only if stated explicitly there.

\begin{ass}\label{i1a}
(H0)  The $A_j(u)$ are $C^\infty$ functions from a connected open
set $\cU^*\subset\bR^N$ to $\bR^{N\times N}$.    For all $u\in\cU^*$
the matrix $A_0(u)$ is invertible.

(H1)(Hyperbolicity near endstates)\;   Let
\begin{align}\label{i1b}
\oA_j:=A_0^{-1}A_j\text{ and }\oA(u,\xi):=\sum^d_{j=1}\xi_j\oA_j(u).
\end{align}
For $u_\pm$ in connected open sets $\cU_\pm\subset\cU^*$ and
$\xi\in\bR^d\setminus\{0\}$, the eigenvalues of $\oA(u_\pm,\xi)$ are
real.

\end{ass}

For $(s,h)\in\bR^d$ let
\begin{align}\label{i2}
\cA_d(u,s,h):=A_d(u)-sA_0(u)-\sum^{d-1}_{j=1}h_iA_i(u) \text{ and
}\ocA_d:=A_0^{-1}\cA_d.
\end{align}
Given  $q=(u_+,u_-,s,h)\in\cU_+\times\cU_-\times\bR^d$, let $R_-(q)$
(resp. $L_+(q)$) denote the number of negative (resp. positive)
eigenvalues of $\ocA_d(u_+,s,h)$ (resp. $\ocA_d(u_-,s,h)$).

\begin{ass}\label{i3}
Let $k$ be an integer such that  $0\leq k\leq N+1$.   We are given a
connected $N+d-k$ dimensional $C^\infty$ submanifold
\begin{align}\notag
\cC\subset\cU_+\times\cU_-\times\bR^d
\end{align}
such that for all $q=(u_+,u_-,s,h)\in\cC$, $\ocA_d(u_\pm,s,h)$ is
invertible.  $R_-$ and $L_+$ are now independent of $q\in\cC$ and we
suppose $R_-+L_+=N+1-k$. We'll refer to $\cC$ as a \emph{shock
manifold} and $k$ as the \emph{undercompressive index}.
\end{ass}

\begin{defn}\label{i4}
Let $t=x_0$ and $y=(x_1,\dots,x_{d-1})$.   For $T>0$ an
\emph{inviscid $\cC$-shock} on $[0,T]\times \bR^d$ is a triple of
functions
\begin{align}\label{i6}
(u_+(t,y,x_d),u_-(t,y,x_d),\psi(t,y))
\end{align}
taking values in $\bR^N\times \bR^N\times\bR$, with $\psi\in
C^1([0,T]\times\bR^{d-1})$ and $u_\pm$ of class $C^1$ in
$\pm(x-\psi(t,y))\geq 0$ respectively, and satisfying
\begin{align}\label{i5}
\begin{split}
&\sum^d_{j=0}A_j(u_\pm)\partial_ju_\pm=0\text{ in }\pm(x_d-\psi(t,y))\geq 0\\
&(u_+(t,y,\psi(t,y)),u_-(t,y,\psi(t,y)),d\psi(t,y))\in\cC\text{ for
}(t,y)\in[0,T]\times\bR^{d-1}.
\end{split}
\end{align}
When $k=0$ the $\cC$-shock is said to be of \emph{Lax type};
otherwise it is called \emph{undercompressive}.
\end{defn}

The problem \eqref{i5} can be regarded as a transmission problem
with a free interface for the unknowns \eqref{i6}.    To each point
$(u_+,u_-,s,h)\in\cC$ we can associate the \emph{planar $\cC$-shock}
\begin{align}\label{i7}
(u_+,u_-,st+hy).
\end{align}
With slight abuse we'll sometimes refer to the point $(u_+,u_-,s,h)$
itself as a planar $\cC$-shock.   Later we will often use the
notation $(p,s,h):=(p_+,p_-,s,h)$ for planar $\cC$-shocks.

Making the change of variable $\tx:=x-\psi(t,y)$ we see that the
problem \eqref{i5} is equivalent to the transmission problem with
flat interface
\begin{align}\label{i8}
\begin{split}
&(a)\;\sum^{d-1}_{j=0}A_j(u_\pm)\partial_ju_\pm+\cA_d(u_\pm,d\psi)\partial_du_\pm=0\text{
in }\pm x\geq 0\\
&(b)\;(u_+(t,y,0),u_-(t,y,0),d\psi(t,y))\in\cC\text{ for
}(t,y)\in[0,T]\times\bR^{d-1}.
\end{split}
\end{align}
We shall also refer to solutions of $(u_+,u_-,\psi)$ of \eqref{i8}
as inviscid $\cC$-shocks.

\begin{defn}\label{i8y}
Given a shock manifold $\cC$ as in Assumption \ref{i3} and
$q\in\cC$, a \emph{local defining function} for $\cC$ near $q$ is a
$C^\infty$ function $\chi:\cO\to\bR^{N+k}$, where
$\cO\subset\bR^{2N+d}$ is an open neighborhood of $q$, such that
$\nabla_q\chi(q)$ has full rank $N+k$ and
\begin{align}\label{i8x}
\cC\cap\cO=\{(p_+,p_-,s,h):\chi=0\}.
\end{align}

\end{defn}

\begin{rem}\label{i8z}
\textup{1. The point of introducing inviscid $\cC$-shocks is to
separate out the part of the construction of nonconservative
inviscid shocks that can be done without reference to  viscosity. We
are mainly interested in the special cases where $\cC$ represents
the set of endstates $u_\pm$, speeds $s$, and directions $h$ for
which there exists an associated viscous connection, and where the
defining function $\chi$ is derived from an associated Melnikov
separation function (see Section \ref{inviscidCB}).}

\textup{2.  In section \ref{Lopatinski} we define a spectral
stability condition for the problem \eqref{i8}, the \emph{uniform
Lopatinski condition}, which generalizes the classical uniform
stability condition of Majda \cite{Ma1}.    We show that the
validity of the condition at a point $q\in\cC$ depends only on the
inviscid operator in \eqref{i1} and the manifold $\cC$; in
particular, it is independent of the choice of local defining
function for $\cC$.}

\textup{3. Observe that it follows immediately from Assumption
\ref{i3} and Definition \ref{i4} that $\cC$-shocks are always
\emph{noncharacteristic}.}

\end{rem}

In section \ref{existence} we prove the following existence theorem
for inviscid $\cC$-shocks using results of \cite{MZ2, Mo,Cou,Ma2}
and this paper.   The theorem assumes the existence of a \emph{
$K$-family of smooth inviscid symmetrizers} for the linearization of
the interior problem \eqref{i8}(a) (see Remark \ref{i9z}.) In the
following statement  we set $\obR^d_\pm=\{(y,x_d):\pm x_d\geq 0\}$
and let
\begin{align}\label{i8h}
H^s(\obR^d_\pm)=\{u\in H^s_{loc}(\obR^d_\pm):u \text{ is constant
outside some compact subset of }\obR^d_\pm\}.
\end{align}
We define $H^s$ spaces on other unbounded domains similarly.
%The notion of \emph{compatible shock front initial data of order
%$s$} is defined in section \ref{existence}.

\begin{thm}\label{i8p}
Suppose the inviscid operator satisfies Assumption \ref{i1a} and,
for a given integer $k$ with $0\leq k\leq N-1$, let $\cC$ be an
$N+d-k$ dimensional $\cC$-shock manifold as in Assumption \ref{i3}.

1.  Consider a planar shock $q=(p_+,p_-,s,h)\in\cC$. Suppose that
the uniform Lopatinski condition is satisfied at $q$ and that a
$K$-family of smooth inviscid symmetrizers exists on a neighborhood
of $q$ in $\cU_+\times\cU_-\times\bR^d$.   Then for any finite
$T_0>0$ there exist nonplanar $\cC$-shocks \eqref{i6} on
$[0,T_0]\times\bR^d$ that are near $(p_+,p_-,st+hy)$ in $C^1$ norm.

2.    Let $s>\frac{d}{2}+1$ and suppose that for $u^0_\pm(y,x_d)\in
H^{s+1}(\obR^d_\pm)$, $\psi^0(y)\in H^{s+2}(\bR^{d-1})$, there
exists an $\bR$-valued function $\sigma(y)\in H^{s+1}(\bR^d)$ such
that
\begin{align}\label{i8i}
(u^0_+(y,0),u^0_-(y,0),\sigma(y),\nabla_y\psi^0(y))\in\cC \text{ for
}y\in\bR^{d-1},
\end{align}
and the uniform Lopatinski condition holds at all points
\eqref{i8i}. Suppose also that $u^0_\pm(y,x_d)$, $\sigma(y)$, and
$\psi^0(y)$ determine  shock front initial data compatible to order
$s-1$ (Definition \ref{sfd}). Then there exists a time $T_0>0$ and
functions
\begin{align}\notag
u_\pm(t,y,x_d)\in H^s([0,T_0]\times\obR^d_\pm),\;\; \psi(t,y)\in
H^{s+1}([0,T_0]\times\bR^{d-1})
\end{align}
satisfying \eqref{i8} on $[0,T_0]\times \obR^d_\pm$ with initial
data
\begin{align}\label{i8j}
u_\pm|_{t=0}=u^0_\pm,\;\psi|_{t=0}=\psi^0.
\end{align}
\end{thm}

\begin{rem}\label{i9z}

\textup{1. $K$-families of smooth inviscid symmetrizers are defined
in \cite{MZ2}, Defn. 4.9, and sufficient conditions on the symbol of
the linearized system are given there (\cite{MZ2}, Theorems 3.4 and
5.6) for the existence of such families.   These conditions are
satisfied, for example whenever the eigenvalues of $\oA(u_\pm,\xi)$
(as in \eqref{i1b}) are real and semisimple with constant
multiplicity for $u_\pm\in\cU_\pm\subset\cU^*$,
$\xi\in\bR^d\setminus \{0\}$. This is the case with Euler shocks
\cite{Ma1}.  In addition, inviscid $K$-families exist in some cases
where characteristics of variable multiplicity are allowed; for
example, in Friedrichs symmetrizable hyperbolic systems where all
real characteristic roots are either \emph{geometrically regular} or
\emph{totally nonglancing} (\cite{MZ2}, Theorem 5.6). Examples of
this sort are provided by both fast and slow inviscid MHD shocks
under generically satisfied conditions on the size of the magnetic
field (\cite{MZ2}, Appendix A.)}

\textup{The requirement that all real characteristic roots are
geometrically regular is shown in \cite{MZ2}, Theorem 3.4 to be
equivalent to Majda's block structure condition \cite{Ma1}. Note
also that the sufficient conditions described above for the
existence of inviscid $K$-families apply equally well to
nonconservative problems.}

\textup{The construction of a $K$-family of inviscid symmetrizers
$\Sigma_K$ is an intermediate step in the construction of a Kreiss
symmetrizer that is independent of the boundary (or transmission)
condition. When the uniform Lopatinski condition is satisfied, a
Kreiss symmetrizer is obtained from a $K$-family by taking $K$
sufficiently large (\cite{MZ2}, Prop. 4.10).}

\end{rem}

\subsection{Inviscid  $\cC_\cB$-shocks}\label{inviscidCB}

We will be especially interested in the case where the shock
manifold $\cC$ is derived from a viscosity.    Consider the $N\times
N$ viscous system on $\bR^{d+1}$
\begin{equation}\label{i9}
\sum_{j=0}^d A_j(u^\eps)\partial_j(u^\eps)  - \eps \sum_{j,k= 1}^d
\D_j \big( B_{j,k}(u^\eps) \D_k u^\eps \big) = 0,
\end{equation}
where $\epsilon>0$ and the $A_j$ satisfy Assumption \ref{i1a}. Set
\begin{align}\label{i10z}
\oB_{j,k}:=A_0^{-1}B_{j,k}\text{ and
}\oB(u,\xi):=\sum^d_{j,k=1}\oB_{j,k}(u)\xi_j\xi_k.
\end{align}

\begin{ass}\label{i10}
(H2)\;(Artificial Viscosity) The $B_{j,k}$ are $C^\infty$ functions
on $\cU^*$ valued in $\bR^{N\times N}$.   There is a positive
constant $c$ such that for all $u\in\cU^*$ and $\xi\in\bR^d$ the
eigenvalues $\mu$ of $\oB(u,\xi)$ satisfy $\Re\mu\geq c|\xi|^2$.

(H3)\;(Strict dissipativity near endstates) There is a positive
constant $c$ such that for all $u\in\cU_\pm$ and $\xi\in\bR^d$ the
eigenvalues $\mu$ of $i\oA(u,\xi)+\oB(u,\xi)$ satisfy $\Re\mu\geq
c|\xi|^2$.
\end{ass}

\begin{rem}\label{i11}
\textup{In order to treat Navier-Stokes or viscous MHD shocks, for
example,  we must of course allow partially parabolic (or
\emph{real}) viscosities. The case of real viscosities involves
substantial additional difficulties which have already been
addressed in \cite{GMWZ4} for conservative Lax shocks. These
difficulties are confined to the high frequency analysis and are
unaffected by the distinctions conservative/nonconservative or
Lax/undercompressive.  The exposition is lighter for artificial
viscosity; we discuss the changes needed to handle real viscosities
in an appendix.}
\end{rem}

We can look for exact travelling wave solutions to \eqref{i9} of the
form
\begin{align}\label{i12}
u^\eps(t,y,x_d)=w\left(\frac{x_d-st-hy}{\eps}\right).
\end{align}
Setting $\nu=\nu(s,h):=(-s,-h_1,\dots,-h_{d-1},1)$ we note that
\begin{align}\label{i13}
\cA_d(u,s,h)=\sum^d_{j=0}A_j(u)\nu_j\text{ and define
}\cB_{d,d}(u,s,h):=\sum^d_{j,k=1}B_{j,k}(u)\nu_j\nu_k.
\end{align}
Then $u^\eps$ is a solution of the parabolic problem \eqref{i9} if
and only if $w(z)$ is a solution of the profile equation
\begin{align}\label{i14}
\cA_d(w,s,h)\partial_z w-\partial_z(\cB_{d,d}(w,s,h)\partial_zw)=0
\text{ on }-\infty < z<\infty.
\end{align}
Given $q=(u_+,u_-,s,h)\in\bR^{2N+d}$, the solution $w(z)$ is said to
be a \emph{viscous profile associated to $q$} when in addition to
\eqref{i14} we have
\begin{align}\label{i15}
\lim_{z\to\pm\infty}w(z)=u_\pm,
\end{align}
and we write $w(z)=W(z,q)$.

We are mainly interested in problems where Assumption \ref{i3} is
replaced by the following stronger assumption.

\begin{ass}\label{i16}
For a fixed $k$ we are given a shock manifold $\cC$ as in Assumption
\ref{i3} together with a $C^\infty$ function $W(z,q)$ from
$\bR\times\cC$ to $\cU^*$ such that for all $q=(u_+,u_-,s,h)\in\cC$,
$W(z,q)$ satisfies \eqref{i14} and \eqref{i15}; that is, for all
$q\in\cC$, $W(z,q)$ is a viscous profile associated to $q$.
\end{ass}

\begin{defn}\label{i16a}
When the points of a  shock manifold $\cC$  are associated to
viscous profiles $W(z,q)$ corresponding to a particular viscosity
\begin{align}\label{i17}
\cB(u):=\sum^d_{j,k=1}\partial_j(B_{j,k}(u)\partial_k u)
\end{align}
as in Assumption \ref{i16}, we'll write $\cC=\cC_\cB$ and refer to
the associated shocks (Definition \ref{i4}) as \emph{inviscid
$\cC_\cB$-shocks}.   We call the condition
\begin{align}\label{i17z}
(u_+(t,y,0),u_-(t,y,0),d\psi(t,y))\in\cC_\cB\text{ for
}(t,y)\in[0,T]\times\bR^{d-1}
\end{align}
the \emph{generalized Rankine-Hugoniot (GRH) condition} determined
by $\cB$.

\end{defn}

\begin{rem}\label{i14z}

1.  \textup{The definition of inviscid $\cC_\cB$-shock given above
is the same as the notion of \emph{nonconservative shock} associated
to a viscosity $\cB$ used, for example, in \cite{Lin}.   In that
paper Lin proves the existence of (what we call) one-dimensional
inviscid $\cC_\cB$-shocks by methods that are quite different from
those we use to prove Theorem \ref{i8p}.}

2. \textup{It follows from (H2), (H3) that for
$q=(u_+,u_-,s,h)\in\cC$, the matrix
\begin{align}\notag
\cB_{d,d}^{-1}(u_+,s,h)\cA_d(u_+,s,h)
\end{align}
has no eigenvalues on the imaginary axis.  The number of eigenvalues
with negative real part is the same as for $\ocA_d(u_+,s,h)$, namely
$R_-$.  The analogous statement holds for $u_-$ with $L_+$ replacing
$R_-$ (\cite{Me1}, Lemma 5.1.3).   Consequently, for $q$ in a
compact subset of $\cC$, standard ODE theory (e.g., \cite{Me1},
Lemma 5.3.3) implies that the associated profiles $W(z,q)$ satisfy
\begin{align}\label{a1}
|\partial_zW(z,q)|\leq Ce^{-\delta|z|}\text{ for some }\delta>0.
\end{align}}
3. \textup{Planar Lax-type shock profiles
may be constructed by center manifold reduction using the method of
Schechter \cite{S}. Both artificial
%NOTES: Schecter doesn't construct undercompressives,
%and in fact it is much more difficult I think.  One must assume failure
%of genuine nonlinearity, to start with, else Schecter's construction
%shows that one gets Lax type shocks.  At least deserves a bit of discussion-
%how about ``we expect that one could prove this by a combination of the
%method of Schecter and the three-parameter bifurcation study of
%\cite{AMPZ} (Azevedo-Marchesin-Plohr-Zumbrun: bifurcation of.....).
and real viscosities are treated there.
Small-amplitude undercompressive profiles might be constructible by
a combination of the methods of Schecter \cite{S} and the
bifurcation analysis used in \cite{AMPZ} to construct such profiles
%CHANGED(added this parenthetical remark)
%in the conservative case.
in the conservative case.\footnote{
Schecter's analysis requires both strict hyperbolicity and
genuine nonlinearity of the associated hyperbolic system.
As pointed out in \cite{AMPZ},
at least one of these must fail in the undercompressive case
by examination of characteristic speeds, leading to a more
complicated, codimension-three bifurcation as compared
to the codimension-one bifurcation of the Lax case.
}
%ENDCHANGED
Constructions using fixed-point arguments are given in Sainsaulieu
\cite{Sa} for the Lax case. Starting with a single transversal
profile associated to a planar shock
$\uq=(\up_+,\up_-,\us,\uh)\in\bR^{2N+d}$, we show in Proposition
\ref{b30} how to construct a $\cC_\cB$-manifold near $\uq$. }

%3. \textup{Planar shock profiles in both the Lax and
%undercompressive cases may be constructed by center manifold
%reduction using the method of Schechter \cite{S}. Both artificial
%and real viscosities are treated there. Constructions using
%fixed-point arguments are given in Sainsaulieu \cite{Sa}. Starting
%with a single transversal profile associated to a planar shock
%$\uq=(\up_+,\up_-,\us,\uh)\in\bR^{2N+d}$, we show in Proposition
%\ref{b30} how to construct a $\cC_\cB$-manifold near $\uq$. }

4.  \textup{In \cite{Mo} Mokrane studies strictly hyperbolic
nonconservative boundary problems involving an additional front-like
unknown $\psi$, together with boundary conditions that are uniformly
stable in the sense of Majda.  This class of problems includes the
case of conservative Lax shocks.     Although Mokrane makes no
reference to undercompressive shocks and although his class does not
strictly include what we call inviscid $\cC$-shocks (since his
boundary conditions are not fully nonlinear in $\nabla\psi$), the
estimates and methods contained in his proofs allow for the
treatment of undercompressives in cases where block structure in the
sense of Kreiss/Majda can be attained. In particular his discussion
of the adjoint problem is fully adequate to yield existence of
inviscid undercompressive shocks.
}

5.  \textup{ An interesting example of (conservative)
undercompressive shocks arises in the study of isothermal phase
transitions in a van der Waals fluid.  Here the underlying
hyperbolic problem is the system of isothermal Euler equations with
the van der Waals pressure law
\begin{align}\label{vdw}
p(\rho)=P(v)=\frac{RT}{v-b}-\frac{a}{v^2};\quad\;\;v=\frac{1}{\rho},\;a>0,b>0.
\end{align}
In \cite{B} Benzoni-Gavage proves the existence of planar phase
transitions admissible in the sense of Slemrod's
viscosity-capillarity criterion \cite{S}.  These phase transitions
are undercompressive shock solutions to the Euler equations ($d=2$
or $3$) with undercompressive index $k=1$, and they satisfy the
uniform Lopatinski condition.    They provide an example of
$\cC$-shocks where the shock manifold $\cC$  comes from a
viscosity-capillarity regularization term rather than a viscosity
regularization term as in the case of $\cC_\cB$ shocks.}

\end{rem}

\subsection{Viscous $\cC_\cB$-shocks}
If we introduce an unknown front $x_d=\psi^\eps(t,y)$ and change
variables
\begin{align}\label{i20}
\tx_d=x_d-\psi^\eps(t,y),
\end{align}
the problem \eqref{i9} can be rewritten
\begin{align}\label{i21}
\sum^{d-1}_{j=0}A_j(u^\eps)\partial_j
u^\eps+\cA_d(u^\eps,d\psi^\eps)\partial_{\tx_d}u^\eps-\eps\sum^d_{j,k=1}D_j(B_{j,k}(u^\eps)D_k
u^\eps)=0\text{ on }\bR^{d+1},
\end{align}
where $D_j=\partial_j-(\partial_j\psi^\eps)\partial_{\tx_d}$, for
$j=1,\dots,d-1$ and $D_d=\partial_{\tx_d}$.     The introduction of
the \emph{viscous front} $\psi^\eps$ allows us as in \cite{GW,GMWZ3}
to reformulate \eqref{i20} as a parabolic transmission problem:
\begin{align}\label{i22}
\begin{split}
&\sum^{d-1}_{j=0}A_j(u^\eps)\partial_j
u^\eps+\cA_d(u^\eps,d\psi^\eps)\partial_{\tx_d}u^\eps-\eps\sum^d_{j,k=1}D_j(B_{j,k}(u^\eps)D_k
u^\eps)=0\text{ on }\obR^{d+1}_\pm,\\
&[u^\eps]=0, [\partial_{\tx_d}u^\eps]=0 \text{ on }\tx_d=0.
\end{split}
\end{align}
We can regard $\psi^\eps$ as an additional unknown, and then we are
forced to add an extra transmission condition as in \eqref{h1} to
obtain a well-posed problem \cite{GMWZ3}.

Henceforth we drop the tildes appearing in \eqref{i22}.

\begin{defn}\label{i23}
Let $(u_+(t,y,x_d),u_-(t,y,x_d),\psi(t,y))$ be an inviscid
$\cC_\cB$-shock satisfying the hyperbolic transmission problem
\eqref{i8}.  A \emph{viscous $\cC_\cB$-shock} is a family of exact
solutions $(u^\eps_\pm,\psi^\eps)$ of \eqref{i22} such that
\begin{align}\label{i24}
\begin{split}
&u^\eps_\pm\to u_\pm\text{ in }L^2_{loc}([0,T]\times\obR^d_\pm)\\
&u^\eps_\pm\to u_\pm\text{ in }L^\infty_{loc}([0,T]\times\bR^d_\pm)\\
&\psi^\eps\to\psi\text{ in }L^\infty_{loc}([0,T]\times\bR^{d-1})
\end{split}
\end{align}
as $\eps\to 0$.   When the inviscid shock is planar, the associated
viscous shock is called a \emph{planar viscous $\cC_\cB$-shock}.

\end{defn}

\begin{ex}\label{i25}
Let $q=(p_+,p_-,s,h)\in\cC_\cB$ be a planar inviscid
$\cC_\cB$-shock. Then the family
\begin{align}\label{i26}
\begin{split}
&u^\eps_\pm(t,y,x_d):=W\left(\frac{x_d}{\eps},q\right)|_{\pm x_d\geq 0}\\
&\psi^\eps(t,y)=st+hy
\end{split}
\end{align}
is a planar viscous $\cC_\cB$-shock satisfying \eqref{i24}.

\end{ex}

We are interested in the existence of nonplanar viscous $\cC_\cB$
shocks.  In section \ref{standard} we recall the definition of the
\emph{standard Evans function} $D_s(q,\hzeta,\rho)$.  For
$q\in\cC_\cB$ this is a spectral stability function for the
transmission problem obtained by linearizing \eqref{i22} with
respect to $u^\eps$ about $W(\frac{x_d}{\eps},q)$.

\begin{defn}\label{i29}
We say that the \emph{standard uniform Evans condition} is satisfied
at $q$ when there exist positive constants $c$, $\rho_0$ such that
\begin{align}\label{i30}
\begin{split}
&(a)\;|D_s(q,\hzeta,\rho)|\geq c\rho\text{ for }\hzeta\in\oS^d_+,
\;0<\rho\leq\rho_0\\
&(b)\;D_s(q,\hzeta,\rho)\neq 0\text{ for }\rho>0.
\end{split}
\end{align}
Here
$\oS^d_+=\{\zeta=(\tau,\gamma,\eta)\in\bR\times\bR\times\bR^{d-1}:|\zeta|=1,
\gamma\geq 0\}.$
\end{defn}

The next Theorem extends the main result of \cite{GMWZ3}, which
applied to conservative Lax shocks in cases where the hyperbolic
characteristics were of constant multiplicity and the viscosity was
fully parabolic, to nonclassical (nonconservative, undercompressive)
shocks.   The result works in certain cases where characteristics
have variable multiplicities.  The extension to real viscosity is
discussed in an appendix.

\begin{thm}\label{i27}
Consider the viscous transmission problem \eqref{i22} under the
structural Assumptions \ref{i1a} and \ref{i10}, and let the shock
manifold $\cC=\cC_\cB$ be as in Assumption \ref{i16}. Let
$(u_+(t,y,x_d),u_-(t,y,x_d),\psi(t,y))$ be an inviscid
$\cC_\cB$-shock satisfying the hyperbolic transmission problem
\eqref{i8} on $[-T_0,T_0]\times \obR^d_\pm$.  Let
\begin{align}\label{i28}
q(t,y):=(u_+(t,y,0),u_-(t,y,0),d\psi(t,y))\in\cC_\cB\text{ for all }
(t,y)\in [-T_0,T_0]\times\bR^{d-1}
\end{align}
be constant outside a compact set.    Suppose for all $(t,y)$ that
the standard uniform Evans condition holds at $q(t,y)$  and that a
K-family of smooth viscous symmetrizers exists on a neighborhood of
$q(t,y)$ in $\cU_+\times\cU_-\times\bR^d$.   Then provided
$(u_+,u_-,\psi)$ is sufficiently smooth, there exists a family of
viscous $\cC_\cB$-shocks $(u^\eps_+,u^\eps_-,\psi^\eps)$ satisfying
\eqref{i24}.
\end{thm}

\begin{rem}\label{i31}
\textup{$K$-families of smooth viscous symmetrizers are defined in
\cite{GMWZ5}, Definition 3.5.   Under the structural assumptions of
sections \ref{inviscidC} and \ref{inviscidCB} (or section
\ref{inviscidC} and Appendix A), such families always exist when the
hyperbolic characteristics are real and semisimple with constant
multiplicity.  This is the case for Euler shocks with artificial
viscosity or Navier-Stokes regularization \cite{GMWZ3,GMWZ4}.}

\textup{More generally, the main result of \cite{GMWZ5}, Theorem
3.7, shows that viscous K-families exist in either one of the
following two situations, both of which allow variable
multiplicities and nonconservative form:}

\textup{(a)\;all real characteristic roots satisfy a \emph{block
structure condition} (\cite{GMWZ5}, Defn. 4.9, condition BS}),

\textup{(b)\;the system is symmetric dissipative in the sense of
Kawashima (\cite{GMWZ5}, Defn. 2.5) and real characteristic roots
are either totally nonglancing \cite{GMWZ5}, Defn. 4.3) or satisfy
the above block structure condition.}

\textup{The block structure condition just referred to is now more
complicated than the one in the inviscid case.  In the inviscid case
the condition is the same as Majda's condition \cite{Ma1} and is
shown in \cite{MZ2} to be equivalent to geometric regularity of the
characteristic root (\cite{MZ2}, Theorem 3.4).  In the viscous case,
in addition to geometric regularity one must require that the
viscous perturbation ``respect" the decoupling between incoming and
outgoing modes (\cite{GMWZ5}, Defn. 4.9).}

\textup{An example of situation (b) is that of fast Lax shocks for
the equations of viscous MHD (\cite{GMWZ5}, section 8).  It is also
shown in \cite{GMWZ5} that smooth $K$-families do \emph{not} exist
for slow shocks for the viscous MHD equations. In fact,
\emph{viscous continuity} (see Defn. \ref{j1} below) is a necessary
condition for the existence of smooth $K$-families of viscous
symmetrizers, and this condition is shown to fail for slow MHD
shocks.}

\textup{Fast Lax shocks with a small magnetic field  are
perturbations of acoustic, gas dynamical shocks, so there is good
reason to expect the standard Evans condition to be satisfied in
this case, at least for ideal gas laws.}

\end{rem}

%\subsection{Standard and modified Evans functions}

\subsection{Overview of the main results}
\textbf{\quad}In most of the paper starting with section 2 we
restrict our exposition to the case of parabolic systems of the form
\begin{align}\label{i32}
\sum^d_{j=0}A_j(u)\partial_ju-\eps\triangle u=0 \text{ with
}A_0(u)=I.
\end{align}
The general fully parabolic case where the $B_{j,k}$ are just
required to satisfy Assumption \ref{i10}  does, in fact, involve
substantial additional difficulties in comparison with \eqref{i32},
but these have already been dealt with in \cite{MZ3,GMWZ5}.  In
particular, the reduction to generalized block structure
(\cite{GMWZ5}, Defn. 4.22) is much harder in the general case. Since
these additional difficulties are unrelated to the distinctions
Lax/undercompressive or conservative/nonconservative, we prefer to
lighten the exposition and focus on the case of Laplacian viscosity.
For the part of the theory we are presenting here, only a
straightforward notational adjustment is needed for passage to the
general fully parabolic case. In Appendix A we discuss the further
changes needed to treat real (partially parabolic) viscosities.

\subsubsection{Constructing $\cC_\cB$ and characterizing
$T_q\cC_\cB$}

\textbf{\quad}With $\cB(u)=\triangle u$ now and starting with a
single given transversal viscous profile $\uw(z)$ (Definition
\ref{c11}) associated to a planar shock $\uq=(\up_+,\up_-,0,0)$, we
show in Proposition \ref{b30} how to construct a shock manifold
$\cC_\cB$ near $\uq$, thereby providing a local verification of
Assumption \ref{i16}. Regarding the profile equation equivalently as
a transmission problem
\begin{align}\label{i33}
\begin{split}
&(a)\;(1+|h|^2)w''=\cA_d(w,s,h)w'\text{ on }\pm z\geq 0\\
&(b)\;[w]=0,\; [w_z]=0\text{ on } z=0
\end{split}
\end{align}
with unknowns $(w_\pm(z),s,h)$, we first obtain in Proposition
\ref{b17a} a convenient parametrization of all possible solutions of
\eqref{i33}(a) with endstates $p_\pm$ for $q=(p_+,p_-,s,h)$ near
$\uq$.  The parametrization is given in terms of the functions
$\phi_\pm(z,p_\pm,s,h,a_\pm)$ defined in \eqref{bn2}.   After adding
a third transmission condition to \eqref{i33},
\begin{align}\label{i33z}
s+w_+\cdot\uw_z-\uw_z\cdot\uw=0\text{ on }z=0,
\end{align}
chosen so that certain rank conditions \eqref{b31} are satisfied, we
construct $\cC_\cB$  by a direct application of the implicit
function theorem. This yields local defining functions $\chi(q)$ for
$\cC_\cB$.  In an appendix we prove that the manifold $\cC_\cB$
obtained by this process is independent of the choice of the third
transmission condition.

Just as membership in $\cC_\cB$ defines the GRH condition,
membership in $T_q\cC_\cB$ defines the linearized GRH condition at
$q$.   One has, of course, the obvious characterization of
$T_q\cC_\cB$ in terms of the implicitly defined function $\chi$ as
\begin{align}\label{i34}
T_q\cC_\cB=\{(\dop_+,\dop_-,\dos,\doh)\in\bR^{2N+d}:\chi'(q)(\dop_+,\dop_-,\dos,\doh)=0\}.
\end{align}
For the later analysis of the relationship between the standard and
modified Evans functions in the low frequency regime, it is
important to have a more explicit description of $T_q\cC_\cB$. Such
a description can be obtained by considering the \emph{full}
linearization with respect to $(w,s,h)$ of \eqref{i33}, \eqref{i33z}
at $(W(z,q),s,h)$:
\begin{align}\label{i36}
\begin{split}
&(a)\cL_0(z,q,\partial_z)\dow:=-(1+|h|^2)\dow''+\cA_d(W,s,h)\dow'+\partial_w\cA_d(W,s,h)\dow W'=\\
&\qquad\qquad\dos A_0(W)W'+\sum_{j=1}^{d-1}\doh_jA_j(W)W'+2\sum^{d-1}_{j=1}h_j\doh_jW'':=\cL_{0,1}(z,q)(\dos,\doh)\\
&(b)\Gamma_a(\dow,\dos,\doh):=\begin{pmatrix}[\dow]\\
[\dow_z]\\\dos+\dow_+\cdot\uw_z\end{pmatrix}=0\text{ on }z=0
\end{split}
\end{align}
Using the fact that solutions of \eqref{i36}(a) can be constructed
from the derivatives $\nabla_{p_\pm,s,h,a_\pm}\phi_\pm$, in
 Proposition \ref{c30} we show
\begin{align}\label{i37}
\begin{split}
&T_q\cC_\cB=\\
&\{(\dop_+,\dop_-,\dos,\doh):\text{ there exists a solution
}(\dow_+,\dow_-,\dos,\doh)\text{ of }\eqref{i36} \text{ with
}\lim_{z\to\pm\infty}\dow_\pm=\dop_ \pm\}.
\end{split}
\end{align}

The proof of \eqref{i37} also involves the construction (Prop.
\ref{c6k}) of  smooth functions $\cR_\pm(z,q,\dos,\doh)$,
exponentially decaying to zero as $z\to\pm\infty$ and linear in
$(\dos,\doh)$, such that if we define
\begin{align}\label{i38}
\dov:=\dow-\cR(z,q,\dos,\doh),
\end{align}
then $(\dow,\dos,\doh)$ satisfies \eqref{i36} if and only if
$(\dov,\dos,\doh)$ satisfies
\begin{align}\label{i39}
\begin{split}
&(a)\;\cL_0(z,q,\partial_z)\dov=0\text{ on }\pm z\geq 0\\
&(b)\Gamma_b(q)(\dov,\dos,\doh):=\;\begin{pmatrix}[\dov]+[\cR(z,q,\dos,\doh]\\
[\dov_z]+[\cR_z(z,q,\dos,\doh)]\\\dov_+\cdot\uw_z\end{pmatrix}=0\text{
on }z=0.
\end{split}
\end{align}
This reduces the study of the fully linearized profile problem to
the study of the partially linearized problem, but with a more
complicated boundary condition.  All solutions to \eqref{i39}(a) may
be constructed by rewriting the problem as a $2N\times 2N$ first
order system as in \eqref{c6d}, and then conjugating to block form
as in \eqref{c6p}.

Next we describe  a \emph{reduced} transmission operator,
$\Gamma_{0,red}(q)$, constructed in a simple way from the explicit
operator $\Gamma_b(q)$, with the property
\begin{align}\label{i39z}
\ker\Gamma_{0,red}(q)=T_q\cC_\cB.
\end{align}
Let us denote by $\cS_0(q)$ the space of solutions of \eqref{i39}(a)
that decay to $0$ as $z\to\pm\infty$ and by $\cS(q)$ the space of
bounded solutions of \eqref{i39}(a). By Remark \ref{c9b} we have
\begin{align}\label{i39y}
\dim\cS_0(q)=N+1-k,\;\;\;\dim\cS(q)=2N+(N+1-k).
\end{align}
Consider for a moment the partial linearization with respect to $w$
of \eqref{i33} at $W(z,q)$:
\begin{align}\label{i40}
\begin{split}
&(a)\;\cL_0(z,q,\partial_z)\dov=0\text{ on }\pm z\geq 0\\
&(b)\;\Gamma_c(\dov):=\begin{pmatrix} [\dov] \\
[\dov_z] \end{pmatrix}=0\text{ on }z=0.
\end{split}
\end{align}
Transversality of $W(z,q)$ and translation invariance imply that the
restriction $\Gamma_c|_{\cS_0(q)}$ has a nontrivial kernel spanned
by $W_z(z,q)$ (Prop. \ref{b39}).
%The extra transmission condition used in
%\eqref{i36} was chosen to be  a low frequency ($\rho\to 0$) limit of
%the extra transmission condition used later in \eqref{c47} for a
%similar purpose of desingularization in the definition of the
%modified Evans function $\tD(q,\hzeta,\rho)$.
Since $W_z(z,q)$ is near $\uw_z$ for $q$ near $\uq$, we deduce
immediately that
\begin{align}\label{i39w}
\bF_P(q):=\{\Gamma_b(q)(\dov,0,0):\dov\in\cS_0(q)\}\subset\bR^{2N+1}
\end{align}
has dimension $N+1-k$, since the third transmission condition
``removes the kernel" of $\Gamma_c|_{\cS_0(q)}$.   Choose an
arbitrary $N+k$ dimensional complementary subspace $\bF_{H,\cR}(q)$
such that
\begin{align}\label{i39x}
\bR^{2N+1}=\bF_{H,\cR}(q)\oplus\bF_P(q),
\end{align}
and let $\pi_{H,\cR}(q)$, $\pi_P(q)$ be the associated projections.
Letting $\dov(z,\dop_\pm)$ denote \emph{any} solution of
\eqref{i39}(a) satisfying $\lim_{z\to\pm\infty}\dov=\dop_\pm$, we
have a well-defined map
\begin{align}\label{i39v}
\Gamma_{0,red}(q):\bR^{2N+d}\to\bF_{H,\cR}(q)\subset\bR^{2N+1}
\end{align}
given by
\begin{align}\label{i39wa}
\Gamma_{0,red}(q)(\dop_+,\dop_-,\dos,\doh):=\pi_{H,\cR}(q)\left(\Gamma_b(q)(\dov(z,\dop_\pm),\dos,\doh)\right).
\end{align}
Transversality implies that $\Gamma_{0,red}(q)$ has full rank $N+k$
even when restricted to the subspace $\doh=0$ (Cor. \ref{c300}).
Finally, it is readily shown that $\ker\Gamma_{0,red}(q)$ is equal
to the right side of \eqref{i37}, and this gives \eqref{i39z} (Prop.
\ref{c30}). We remark that the solutions $\dov(z,\dop_\pm)$ can be
expressed explicitly in terms of block form coordinates $(u_H,u_P)$
\eqref{c10}.

\subsubsection{Stability determinants}
\textbf{\quad}At the beginning of section 4 we return to the general
context of $\cC$ shocks and define the Lopatinski and modified
Lopatinski determinants, $D_{Lop}(q,\hzeta)$ and
$D_{Lop,m}(q,\hzeta)$.   The corresponding uniform (resp., modified
uniform) Lopatinski condition at $q\in\cC$ is the condition that
\begin{align}\label{i40a}
|D_{Lop}(q,\hzeta)|\geq c \text{ (resp. } |D_{Lop,m}(q,\hzeta)|\geq
c)
\end{align}
for some $c>0$ independent of $\hzeta\in
S^d_+=\oS^d_+\cap\{\hgamma>0\}$.  The determinant $D_{Lop}$ is the
one that naturally governs the stability of inviscid $\cC$ shocks.
The need to define $D_{Lop,m}$ becomes apparent only later in the
low frequency analysis of the modified Evans function. Although the
definition of  $D_{Lop}$ itself is not independent of the choice of
defining function $\chi$ for $\cC$, we show in Proposition
\ref{c403} that the validity of the uniform Lopatinski condition is
independent of the choice of $\chi$, and thus depends just on the
inviscid operator and $\cC$.    Moreover, we show (Prop.
\ref{c410}):

\begin{prop}\label{i41}
Suppose the uniform Lopatinski condition holds at $q\in\cC$ and
$\chi$ is any defining function for $\cC$ near $q$. Then

(a)\;$\chi_{p,s}'(q)$ has full rank $N+k$;

(b)\;if $d\geq 2$, $\chi'_p(q)$ has full rank $N+k$.
\end{prop}

\begin{rem}\label{mormk}
\textup{This Proposition shows that it is not necessary, as in the
treatments \cite{Mo,Cou}, to introduce the full rank condition in
part (b) as an extra hypothesis when $d\geq 2$.
For $d=1$, the full rank condition is not necessary for other
reasons, since the analysis may be carried out by alternative
%CHANGED(Put back modified version of offending sentence... OK NOW?-K)
%methods that do not require it.
methods that do not require it.\footnote{ For example, one may
(after doubling to obtain a problem on a half-space) replace the
nonstandard front variable $\psi$ with a standard interior variable
$v:=\Psi_t$, where $\Psi:=\psi e^{-z}$ is an extension to $z>0$ of
the front variable $\psi$ (defined only at $z=0$), and add the
corresponding artificial interior equation $v_z+ v=0$, to obtain a
standard hyperbolic initial--boundary-value problem in $(u,v)$ (no
longer involving $\psi$), with linear boundary condition $b_1u +
b_2v=b_0$ (no longer involving $\psi_t$), treatable by the
techniques of \cite{CP}. }
%ENDCHANGED
At the expense of further effort, one may dispense with the rank
condition altogether by defining a pseudodifferential adjoint
problem as in \cite{GMWZ6}. (The rank condition is used only to
cleverly define an adjoint equation that is differential.)
}
\end{rem}

%\begin{rem}\textup{This Proposition shows that it is not necessary, as in the
%treatments \cite{Mo,Cou}, to introduce the full rank condition in
%part (b) as an extra hypothesis when $d\geq 2$.}
%\end{rem}

We show that the uniform Lopatinski condition at $q$ always implies
the modified Lopatinski condition at $q$, and in Proposition
\ref{c408} we give geometric conditions under which the converse
holds. As a corollary of these results for $\cC$-shocks, in the case
when $\cC=\cC_\cB$ and $W(z,q)$ is the viscous profile associated to
$q\in\cC_\cB$, we show that certain transversality assumptions on
$W(z,q)$ imply the equivalence of the two Lopatinski conditions
(Cor. \ref{cb3}).

 Section \ref{stability} continues with the definition of  the standard and modified Evans functions,
$D_s(q,\hzeta,\rho)$ and $D_m(q,\hzeta,\rho)$, for viscous $\cC_\cB$
shocks. We start with the rescaled transmission problem
\begin{align}\label{i42}
\begin{split}
&\sum^{d-1}_{j=0}A_j(u)\partial_j
u+\cA_d(u,d\psi)\partial_zu-\sum^d_{j=1}(\partial_j-\partial_j\psi\partial_z)^2u=0\text{
on }\pm z\geq 0\\
&[u]=0,\;[u_z]=0\text{ on }z=0
\end{split}
\end{align}
for which we have an exact solution given by a profile $W(z,q)$ and
front $\psi=st+hy$.  The determinant $D_s$ is defined by considering
the (Fourier-Laplace transform of the) partial linearization of
\eqref{i42} with respect to $u$ about $W(z,q)$,
\begin{align}\label{i43}
\begin{split}
&\cL(z,q,\zeta,\partial_z)u=f\text{ on }\pm z\geq 0\\
&[u]=0, [u_z]=0\text{ on }z=0,
\end{split}
\end{align}
while $D_m$ is defined by considering the full linearization of
\eqref{i42} with respect to $(u,\psi)$ about $(W(z,q),st+hy)$ and
adding a third transmission condition:
\begin{align}\label{i44}
\begin{split}
&\cL(z,q,\zeta,\partial_z)u-\psi\cL_1(z,q,\zeta)=f\\
&[u]=0,\;[u_z]=0,\; c_0(\zeta)\psi+\uw_z(0)\cdot u^+=0.
\end{split}
\end{align}
Here $\cL$ and $\cL_1$ are given explicitly in \eqref{c5a} and
$c_0(\zeta)=i\tau+\gamma+|\eta|^2$. Moreover, we suppose $q$ is near
a basepoint $\uq$ and   $\uw(z):=W(z,\uq)$.   It is important to
note the correspondence between the operators appearing in
\eqref{i44} and those ($\cL_0$ and $\cL_{0,1}$) appearing in the
fully linearized profile transmission problem \eqref{i36}.   Writing
$\cL_1(z,q,\zeta)=\rho\check{\cL}_1(z,q,\zeta,\rho)$  and
$c_0(\zeta)=\rho\check{c}_0(\hzeta,\rho)$, we have
\begin{align}\label{i45}
\begin{split}
&\cL_0(z,q,\partial_z)=\cL(z,q,0,\partial_z),\;
\cL_{0,1}(z,q)(i\htau+\hgamma,i\heta)=\check{\cL}_1(z,q,\hzeta,0),\\
&\text{ and }i\htau+\hgamma=\check{c}_0(\hzeta,0).
\end{split}
\end{align}

Corresponding to the two Evans functions we have the standard (Defn.
\ref{cc1}) and modified (Defn. \ref{c51}) uniform Evans conditions.
The standard uniform Evans condition is the one that is easier to
verify analytically \cite{PZ,FS} or numerically \cite{B,HZ}.  The
vanishing of $D_s(q,\hzeta,\rho)$ at $\rho=0$ reflects the
translational degeneracy of \eqref{i40} pointed out earlier. On the
other hand the modified uniform Evans condition, which was
introduced in \cite{GMWZ3}, implies maximal $L^2$ estimates for the
fully linearized problem \eqref{i44} and is essential for our
construction of viscous $\cC_\cB$ shocks.   Thus, a key result of
this paper, discussed below, is  the following theorem:

\begin{thm}\label{i45b}
Under Assumptions \ref{i1a}, \ref{i10}, and \ref{i16}, the standard
uniform Evans condition at $q\in\cC_\cB$ implies the modified
uniform Evans condition at $q$.
\end{thm}
This was proved for conservative problems with (hyperbolic)
characteristics of constant multiplicity in \cite{GMWZ3}, but a new
argument is needed for nonconservative, variable multiplicity
systems.
%In addition, this result is essential for our curved
%viscous shock theorem, Theorem \ref{i27}.

The third transmission condition in \eqref{i44} yields a well-posed
problem and, roughly speaking, removes the translational degeneracy
of \eqref{i43} at $\rho=0$.    Parallel to the earlier passage from
\eqref{i36} to \eqref{i39}, the fully linearized problem \eqref{i44}
is most easily studied by reducing it to the partially linearized
problem \eqref{i43} with modified transmission conditions.  In
Proposition \ref{c53a} we recall from \cite{GMWZ3} the construction
of  functions $\ibfR_\pm(z,q,\hzeta,\rho)$ with the property that
$(u,\psi)$ satisfies \eqref{i44} if and only if $(v,\phi)$ defined
by
\begin{align}\label{i46}
\phi:=\rho\psi,\qquad v_\pm=u_\pm-\phi\ibfR_\pm
\end{align}
satisfy
\begin{equation}
\label{i47}
\begin{aligned}
& \cL(z, q, \zeta, \partial_z) v_\pm   = f_\pm \, \quad \mathrm{ on
} \;\pm z \ge 0 \,,
\\
& [v (0) ] + \phi [\ibfR(0)] = 0 \quad
 [ \partial_z v (0) ] + \phi [\partial_z \ibfR (0)]  = 0 ,
\\
 &\uw_z(0) \cdot v^+ (0)  = 0 \,.
\end{aligned}.
\end{equation}
Moreover, the functions $\cR_\pm(z,q,\dos,\doh)$ appearing in
\eqref{i38} and the functions $\ibfR_\pm(z,q,\hzeta,\rho)$ can be
chosen to satisfy
\begin{align}\label{i48}
\cR_\pm(z,q,i\htau+\hgamma,i\heta)=\ibfR(z,q,\hzeta,0),
\end{align}
as might be anticipated from \eqref{i45}.

\begin{rem}\label{i49}
\textup{Considering the equality \eqref{i48}, we see that the
introduction of the functions $\cR_\pm(z,q,\dos,\doh)$ and their use
in characterizing $T_q\cC_\cB$ can be viewed as a ``missing step"
that helps  to clarify and provide a geometric motivation for some
of the low frequency analysis in \cite{GMWZ3}. This step  was not
needed there, and hence overlooked,  because of the assumption of
conservative form.}
\end{rem}

At the end of section \ref{stability} we discuss conditions for the
continuity of decaying eigenspaces of the operators that define the
linearized inviscid and viscous problems. Inviscid continuity
implies, for example, that the uniform Lopatinski condition must
hold for $q\in\cC$ near $\uq$ when it holds at $\uq$. Similarly
viscous continuity implies that the standard uniform Evans condition
holds for $q\in\cC$ near $\uq$ when it holds at $\uq$. Moreover,
when inviscid or viscous continuity holds, sometimes a converse can
be proved (as in the nonconservative Zumbrun-Serre theorem, Theorem
\ref{K}) or a proof can be substantially simplified (see Remark
\ref{sim}).    Viscous continuity (and thus inviscid continuity)
always holds when the linearized hyperbolic problem has
characteristics of constant multiplicity \cite{MZ1}. A sufficient
condition for inviscid (resp., viscous) continuity at $\uq$ is the
existence of a smooth $K$-family of inviscid (resp., viscous)
symmetrizers near $\uq$. In any given proposition we do not assume
inviscid or viscous continuity unless we explicitly say so.

The first main result of section \ref{low frequency} is the
nonconservative Zumbrun-Serre theorem, Theorem \ref{K}:

\begin{thm}\label{i50}
Consider a shock profile $\uw(z)=W(z,\uq)$, where $\uq=(\up,0,0)$,
and suppose the low frequency standard Evans condition
\eqref{i30}(a) holds at $\uq$.  Then $\uw(z)$ is transversal (Defn.
\ref{c11}) and the uniform Lopatinski condition holds at $\uq$.  In
fact, for $\hgamma>0$
\begin{align}\label{i51}
D_s(\uq,\hzeta,\rho)=\rho\alpha(q)D_{Lop}(\uq,\hzeta)+O_{\hgamma}(\rho^2),
\end{align}
where $O_{\hgamma}(\rho^2)\leq C_{\hgamma}\rho^2$ and $\alpha(\uq)$
is a constant whose nonvanishing is equivalent to $a$-transversality
of $\uw$ (Defn. \ref{c11}).
\end{thm}

This theorem does not assume the prior existence of a shock manifold
$\cC_\cB$.  However, the transversality conclusion allows us to
construct a unique $\cC_\cB$ manifold near $\uq$ as in Proposition
\ref{b30}. The proof of Theorem \ref{i50} is quite different from
the argument in the conservative case \cite{ZS} and relies heavily
on the functions $\phi_\pm(z,p_\pm,s,h,a_\pm)$ \eqref{bn2}.   One
advantage of the  new argument is that it generalizes almost
verbatim to the case of real viscosities.

\begin{cor}\label{i51a}
Even in situations, like that of slow MHD shocks, where viscous
continuity fails, the uniform Lopatinski condition is a necessary
condition for the standard uniform Evans condition to hold.

\end{cor}

The remainder of section \ref{low frequency} is mainly concerned
with obtaining block decompositions of the modified and standard
Evans functions for $\rho$ small, Proposition \ref{f6} and Theorem
\ref{g3}.  Roughly speaking, the blocks in the decomposition of the
determinants correspond to the $H$ and $P$ blocks in the conjugated
form \eqref{c6i} of the linearized parabolic problem.  The block
decompositions are useful for relating the standard and modified
Evans to each other and to the Lopatinski and modified Lopatinski
determinants.  Theorem \ref{g3} implies, for example:
\begin{thm}\label{i52}
Assume the profile $W(z,q)$ is strongly transversal (Defn.
\ref{c11}).   Then for some $\rho_0>0$
\begin{align}\label{i53}
D_s(q,\hzeta,\rho)=\rho\alpha(q,\zeta)D_m(q,\hzeta,\rho)+O(\rho^2),
\end{align}
where the error is uniform for $\hzeta\in\oS^d_+$, $0<\rho\leq
\rho_0$ and $\alpha(q,\zeta)$ is $C^\infty$ and bounded away from
zero for $\rho$ small.
\end{thm}

In combination with Theorem \ref{i50} and Proposition \ref{i41} this
implies Theorem \ref{i45b}.   A result like Theorem \ref{i52} was
proved in \cite{GMWZ3} in the conservative, constant multiplicity
case by giving a low frequency expansion of both $D_s$ and $D_m$ in
terms of $D_{Lop}$.  If viscous continuity fails,  this much simpler
type of argument does not work.

In Theorem \ref{j3} we summarize a number of the connections that
hold  between the different Evans and Lopatinski determinants.

Section \ref{approximate} is devoted to the construction of high
order approximate solutions (which converge as viscosity tends to
zero to a given curved inviscid $\cC_\cB$ shock) to the nonlinear
small viscosity transmission problem \eqref{h1}.  The
characterizations of transversality and of $T_q\cC_\cB$ in terms of
properties of the fully linearized profile problem given in
Propositions \ref{c10a} and \ref{c30} lead to a simpler and shorter
construction of higher order profiles than was given in
\cite{GW,GMWZ4}.   In section 7 we complete the construction of
curved inviscid $\cC$-shocks (Theorem \ref{i8p}), and of families of
curved viscous $\cC_\cB$-shocks converging to a given inviscid
$\cC_\cB$-shock as viscosity approaches zero (Theorem \ref{i27}). In
particular, Theorem \ref{i27} shows that the approximate solutions
constructed in section \ref{approximate} are close to true exact
solutions of the parabolic transmission problem \eqref{h1}.

Finally, in Appendix A we explain the changes needed to handle the
partially parabolic case of real viscosities, and in Appendix B we
prove the local uniqueness of $\cC_\cB$ manifolds.

\section{Generalized Rankine-Hugoniot condition}\label{genRH}

\textbf{\;\;}  We let  $\cB(u)=\triangle u$ now and set
$\uq=(\up_+,\up_-,0,0)$, where $\up_\pm\in\cU_\pm$ with
$A_d(\up_\pm)$ invertible.  Suppose $\uw(z)$ is a viscous profile
associated to $\uq$.  Our main task in this section is to determine
when and how we can construct a shock manifold $\cC_\cB$ near $\uq$.

\subsection{The connection problem as a transmission
problem}\label{connection}

Consider again the profile transmission problem \eqref{i33}:
\begin{align}\label{t1}
\begin{split}
&(a)\;w''=(1+|h|^2)^{-1}\cA_d(w,s,h)w'\text{ on }\pm z\geq 0\\
&(b)\;[w]=0,\; [w_z]=0\text{ on } z=0.
\end{split}
\end{align}

As in \eqref{t1} we'll often suppress $\pm$ subscripts on unknowns.
Set
\begin{align}\label{t2}
\begin{pmatrix}u\\v\end{pmatrix}=\begin{pmatrix}w\\w'\end{pmatrix};\;\;b(h):=(1+|h|^2)^{-1}
\end{align}
and rewrite \eqref{t1} as a $2N\times 2N$ first-order system:
\begin{align}\label{t3}
\begin{split}
&(a)\;\begin{pmatrix}u\\v\end{pmatrix}'=\begin{pmatrix}0&1\\0&G_d(u,s,h)\end{pmatrix}\begin{pmatrix}u\\v\end{pmatrix}=\begin{pmatrix}v\\G_d(u,s,h)v\end{pmatrix},
\text{ where }G_d(u,s,h):=b(h)\cA_d(u,s,h)\\
&(b)\;\Gamma\begin{pmatrix}u_+\\u_-\\v_+\\v_-\end{pmatrix}=\begin{pmatrix}
[u]\\ [v]\end{pmatrix}=0\text{ on }z=0,\;\Gamma:\bR^{4N}\to\bR^{2N}.
\end{split}
\end{align}
We seek conditions on $q=(p_+,p_-,s,h)$ that will allow us to find
solutions to \eqref{t3} such that
\begin{align}\label{b13}
\lim_{z\to\pm\infty}u(z)=p_\pm.
\end{align}

\begin{notation}\label{b13z}
Given a matrix $G:\bC^p\to\bC^p$ with no eigenvalues on the
imaginary axis, let $\bE_\pm(G)$ denote the invariant subspace of
$\bC^p$ generated by the generalized eigenvectors of $G$ associated
to eigenvalues $\mu$ such that $\pm\Re\mu>0$.   We denote by
$\Pi_\pm$ the corresponding spectral projectors associated to the
decomposition
\begin{align}\label{b13y}
\bC^p=\bE_+(G)\oplus\bE_-(G).
\end{align}
When $G$ is real (that is, when each entry of $G$ is real), the
$\Pi_\pm$ are real and so we can also regard the $\bE_\pm(G)$ as
real vector spaces.   This gives a decomposition
\begin{align}\label{b13x}
\bR^p=\bE_+(G)\oplus\bE_-(G).
\end{align}
Note that the spaces appearing in \eqref{b13y} are just the
complexifications of the corresponding spaces in \eqref{b13x}.

\end{notation}

In the first step we ignore the transmission conditions and
construct solutions $u_\pm$ to the interior problems in $\pm z\geq
0$ such that \eqref{b13} holds.  We do this by regarding
\eqref{t3}(a) in $\pm z\geq 0$ as a perturbation, quadratic in
$(u_\pm-p_\pm,v_\pm)$, of
\begin{align}\label{b14}
\begin{pmatrix}u\\v\end{pmatrix}'=\begin{pmatrix}v\\G_d(p_\pm,s,h)v\end{pmatrix}.
\end{align}
%Denote by  $\bE_\mp(G_d(p_\pm,s,h))$ the negative (resp., positive)
%eigenspace of $G_d(p_+,s,h)$ (resp., $G_d(p_-,s,h))$ and let
%\Pi_\mp(p_\pm,s,h)$ be the spectral projectors on
%$\bE_\mp(G_d(p_\pm,s,h))$.
Let $\Pi_-(p_+,s,h)$ be the projection on $\bE_-(G_d(p_+,s,h))$ with
respect to the decomposition
\begin{align}\label{b14g}
\bR^N=\bE_-(G_d(p_+,s,h))\oplus\bE_+(G_d(p_+,s,h)),
\end{align}
and define $\Pi_+(p_-,s,h)$ similarly.  For $p_\pm$ in  relatively
compact open neighborhoods $\omega_\pm\subset\cU_\pm$ of $\up_\pm$
and for $|s,h|\leq\epsilon_\omega$, where $\epsilon_\omega$ is a
small enough positive constant,  we fix isomorphisms linear in
$a_\pm\in\bE_\mp(G_d(\up_\pm,0,0))$ and $C^\infty$ in $(p_\pm,s,h)$:
\begin{align}\label{b15}
\alpha_\pm(p_\pm,s,h;a_\pm):\bE_\mp(G_d(\up_\pm,0,0))\to\bE_\mp(G_d(p_\pm,s,h)).
\end{align}

\begin{prop}\label{b17a}
(a)\;Let $\omega_\pm$, $\epsilon_\omega$, and $\alpha_\pm$ be as
just defined.    There are positive constants $R$ and $r$ such that
for all $p_\pm\in\omega_\pm$, $|s,h|\leq \epsilon_\omega$, and
$a_\pm\in\bE_\pm(G_d(\up_\pm,0,0))$ with $|a_\pm|\leq r$, the
equation \eqref{t1}(a) has a unique solution
$u_\pm=\Phi_\pm(z,p_\pm,s,h,a_\pm)$ satisfying
\begin{align}\label{b17b}
\begin{split}
&(i)\;\lim_{z\to\pm\infty}u_\pm=p_\pm\\
&(ii)\;\Pi_\mp(p_\pm,s,h)\partial_zu_\pm(0)=\alpha_\pm(p_\pm,s,h;a_\pm)\\
&(iii)\;\|\partial_zu_\pm\|_{L^1(\pm z\geq 0)}\leq R\text{ and
}\|\partial_z u_\pm\|_{L^\infty(\pm z\geq 0)}\leq R.
\end{split}
\end{align}
The function $\Phi_\pm$ is $C^\infty$ on $\{\pm z\geq 0\}\times
\Omega_\pm$, where
\begin{align}\label{b17c}
\Omega_\pm=\{(p_\pm,s,h,a_\pm):p_\pm\in\omega_\pm, |s,h|\leq
\epsilon_\omega, \text{and }a_\pm\in\bE_\pm(G_d(\up_\pm,0,0))\text{
with }|a_\pm|\leq r\}.
\end{align}
It satisfies
\begin{align}\label{b17d}
\;\Phi_\pm(z,p_\pm,s,h,a_\pm)=p_\pm+e^{zG_d(p_\pm,s,h)}G^{-1}_d(p_\pm,s,h)\alpha(p_\pm,s,h;a_\pm)+O(|a_\pm|^2)
\end{align}
uniformly with respect to $(z,p_\pm,s,h)$.  Moreover, there exist
positive constants $\delta$ and $C_\beta$ such that for all $\pm
z\geq 0$ and $(p_\pm,s,h,a_\pm)\in\Omega_\pm$:
\begin{align}\label{b17e}
|\partial_{z,p_\pm,s,h,a_\pm}^\beta
\left(\Phi_\pm(z,p_\pm,s,h,a_\pm)-p_\pm\right)|\leq C_\beta
e^{-\delta|z|},\;\;|\beta|\leq 2.
\end{align}
We will also denote by $\Phi_\pm(z,p_\pm,s,h,a_\pm)$ the maximal
extension of $\Phi_\pm$ to $\pm z\leq 0$ as a solution of
\eqref{t1}(a).

(b)\;Suppose $u_+$ is a solution of \eqref{t1}(a)$_+$,
\eqref{b17b}(i)$_+$ on $[z_1,+\infty)$ for some $p_+\in\omega_+$ and
$|s,h|\leq\epsilon_\omega$.  Then for all $z_0\geq z_1$ large
enough, we have
\begin{align}\label{b17f}
u_+(z)=\Phi_+(z-z_0,p_+,s,h,a_+)\text{ with
}\alpha(p_+,s,h;a_+)=\Pi_-(p,s,h)\partial_z u_+(z_0).
\end{align}
We may (and shall) always take $z_0\geq 0$.  The analogous statement
holds for solutions $u_-$ on $(-\infty,z_2]$.
\end{prop}

\begin{proof}

\textbf{1. }The proof of part (a) (resp., (b)) is essentially the
same as that of Prop. 5.3.5 (resp., Proposition 5.3.6) in
\cite{Me1}, to which we refer for more detail.  We provide a sketch
for the $+$ case that describes the slight differences.

\textbf{2. }Let $\Pi^+(p_+,s,h)$ denote the projection on the second
summand in \eqref{b14g}.  Note there is a $\theta>0$ such that
\begin{align}\label{b17g}
\begin{split}
&|e^{(z-y)G_d(p_+,s,h)}\Pi_-(p_+,s,h)|\leq Ce^{-\theta(z-y)}, \text{ for }z\geq y\\
&|e^{(z-y)G_d(p_+,s,h)}\Pi^+(p_+,s,h)|\leq Ce^{-\theta|z-y|}, \text{
for }z\leq y.
\end{split}
\end{align}
Define integral operators
\begin{align}\label{b17h}
\begin{split}
&I(v)(z)=-\int^{+\infty}_zv(y)dy\\
&\cI_0(F)(z)=\\
&\int^z_0e^{(z-y)G_d(p_+,s,h)}\Pi_-(p_+,s,h)F(y)dy-\int^{+\infty}_ze^{(z-y)G_d(p_+,s,h)}\Pi^+(p_+,s,h)F(y)dy.
\end{split}
\end{align}
With $F(u,v):=\left(G_d(u,s,h)-G_d(p_+,s,h)\right)v$ we contruct
$(u_+,\partial_zu_+)$ as a fixed point of the map
\begin{align}\label{b17i}
\cT_+(u_+,v_+):=\left(p_++I(v_+),e^{zG_d(p_+,s,h)}\alpha_+(p_+,s,h;a_+)+\cI_0(F(u,v))\right)
\end{align}
on
\begin{align}\label{b17j}
B_R=\{(u_+,v_+):\|u_+-p_+\|_{L^\infty(z\geq 0)}+\|v_+\|_{L^1\cap
L^\infty}\leq R\}.
\end{align}
For $R$ and $r$ small enough $\cT_+$ is a contraction on $B_R$.

Uniqueness follows from the fixed point theorem.    By construction
the solution satisfies $\|v_+\|_{L^1}=O(|a_+|)$; thus,
$\|u_+-p_+\|_{L^\infty}=O(|a_+|)$.  Since
$F(u_+,v_+)=O(|u_+-p_+||v_+|)$ it follows that
\begin{align}\label{b17k}
\|v_+-e^{zG_d(p_+,s,h)}\alpha_+(p_+,s,h;a_+)\|_{L^1\cap
L^\infty}=O(|a_+|^2),
\end{align}
which in turn yields
\begin{align}\label{b17l}
\|u_+-p_+-e^{zG_d(p_+,s,h)}G_d^{-1}(p_+,s,h)\alpha_+(p_+,s,h;a_+)\|_{L^\infty}=O(|a_+|^2).
\end{align}
Exponential decay as in \eqref{b17e} follows by a similar
application of the fixed point theorem in spaces $e^{-\delta
z}L^\infty$.

\textbf{3. } To prove part (b) we note that if $u_+$ is a solution
of \eqref{t1}(a)$_+$, \eqref{b17b}(a)$_+$ on $[z_1,+\infty)$, then
for all $z_0\geq z_1$ large enough there holds
\begin{align}\label{b17m}
\begin{split}
&\|\partial_zu_+\|_{L^1(z\geq
z_0)}+\|\partial_zu_+\|_{L^\infty(z\geq z_0)}\leq
R,\\
&\Pi_-(p_+,s,h)\partial_zu_+(z_0)=\alpha_+(p_+,s,h;a_+)\text{ for
some }|a_+|\leq r,
\end{split}
\end{align}
where $R$ and $r$ are the constants determined in part (a).   By
translation invariance $u_+(z+z_0)$ satisfies \eqref{t1}(a)$_+$,
\eqref{b17b}(a)$_+$ on $z\geq 0$, and using part (a) we obtain
\begin{align}\label{b17n}
u_+(z+z_0)=\Phi_+(z,p_+,s,h,a_+)\text{ for }z\geq 0.
\end{align}

\end{proof}

\begin{defn}\label{b19}
Let $R_\pm$ be the number of eigenvalues $\mu$ of
$\cA_d(\up_+,0,0)$ with $\pm\mu>0$ ($R$ refers to the right
endstate $\up_+$). Let $L_\pm$ be the number of eigenvalues $\mu$
of $\cA_d(\up_-,0,0)$ with $\pm\mu>0$.
\end{defn}

\begin{defn}\label{b20}
  We shall only consider shocks satisfying $R_-+L_+=N+1-k$, where $N$ is the
dimension of the system \eqref{i32} and $k\geq 0$ is defined to be
the \emph{undercompressive index}. Thus,
\begin{align}\label{b21}
\dim
\left(\bE_-(G_d(p_+,s,h))\times\bE_+(G_d(p_-,s,h)\right):=\cN_-=R_-+L_+=N+1-k.
\end{align}
%2.  For an undercompressive shock we have $R_-+L_+\leq N$.  Note
%that we can have a zero-amplitude (and noncharacteristic)
%undercompressive  shock, in which case one necessarily has
%\begin{align}\label{b22}
%R_-+L_+=N, \text{ and thus }\cN_-=N.
%\end{align}
%since $\up_+=\up_-$.
%It is impossible to have a zero-amplitude (and noncharacteristic)
%Lax shock.
When $k=0$ the shock is called a Lax shock; when $k>0$ it is called
undercompressive.

\end{defn}

Observe that if we take the given profile $\uw(z)$ and define
\begin{align}\label{t4}
\begin{pmatrix}u\\v\end{pmatrix}=\begin{pmatrix}\uw(z)\\\uw_z(z)\end{pmatrix},
\end{align}
then $\begin{pmatrix}u\\v\end{pmatrix}$ satisfies \eqref{t3}(a)and
\eqref{b13} with $(s,h)=(0,0)$ and $p_\pm=\up_\pm$. Thus, by part
(b) of Prop. \ref{b17a} there exist $\uz=-z_0\leq 0$ and $\ua_\pm$
such that
\begin{align}\label{b24}
\uw(z)=\begin{cases}\Phi_+(z+\uz,\up_+,0,0,\ua_+), \;z\geq
0\\\Phi_-(z-\uz,\up_-,0,0,\ua_-),\;z\leq 0\end{cases}.
\end{align}

\begin{notation}\label{bn1}
1.  For $\Phi_\pm$ as constructed in Prop. \ref{b17a} it is
convenient to define
\begin{align}\label{bn2}
\begin{split}
&\phi_\pm(z,p_\pm,s,h,a_\pm):=\Phi_\pm(z\pm\uz,p_\pm,s,h,a_\pm),\\
&\phi(z,p,s,h,a):=\begin{cases}\phi_+(z,p_+,s,h,a_+),\;z\geq
0\\\phi_-(z,p_-,s,h,a_-),\;z\leq 0\end{cases}.
\end{split}
\end{align}

2.  Let
\begin{align}\label{bn3}
\bE_-(G_d(p,s,h)):=\bE_-(G_d(p_+,s,h))\times\bE_+(G_d(p_-,s,h).
\end{align}

\end{notation}

The following corollary of Proposition \ref{b17a} will be used in
the proofs of Propositions \ref{b39} and \ref{c26}.

\begin{cor}\label{bp1}
For
$(\dop_\pm,\dos,\doh,\doa_\pm)\in\bR^{N}\times\bR^d\times\bE_\mp(G_d(\up_\pm,0,0))$
we have
\begin{align}\label{bp3}
\lim_{z\to\pm\infty}\nabla_{p_\pm,s,h,a_\pm}\phi_\pm(z,p_\pm,s,h,a_\pm)(\dop_\pm,\dos,\doh,\doa_\pm)=\dop_\pm.
\end{align}
Decreasing $r$ in the definition  of $\Omega_\pm$ \eqref{b17c} if
necessary, we have for $(p_\pm,s,h,a_\pm)\in\Omega_\pm$ that the map
\begin{align}\label{bp2}
(\dop_\pm,\dos,\doh,\doa_\pm)\to\begin{pmatrix}\nabla_{p_\pm,s,h,a_\pm}\phi_\pm(z,p_\pm,s,h,a_\pm)(\dop_\pm,\dos,\doh,\doa_\pm)\\\dos\\\doh\end{pmatrix}
\end{align}
is injective.

\end{cor}

\begin{proof}
Consider the $+$ case.  The limit \eqref{bp3} follows directly from
\eqref{b17e}. Thus, if the image of $(\dop_+,\dos,\doh,\doa_+)$ in
\eqref{bp2} is $0$ we find $(\dop_+,\dos,\doh)=0$.   Since
\eqref{b17d} implies
\begin{align}\label{bp5}
\nabla_{a_+}\phi_+(-\uz,p_+,s,h,a_+)\doa_+=G^{-1}_d(p_+,s,h)\alpha(p_+,s,h;\doa)+O(|a_+|)\doa_+,
\end{align}
we conclude for $r$ small enough that if the right side of
\eqref{bp5} is $0$, then $\doa_+=0$.
\end{proof}

\subsection{Extra transmission condition}\label{extra}

For $p:=(p_+,p_-)$ near $\up=(\up_+,\up_-)$, $(s,h)$ near $(0,0)$,
and $a:=(a_+,a_-)$ near $\ua=(\ua_+,\ua_-)$, define
\begin{align}\label{b25}
\Psi(p,s,h,a)=\Gamma\begin{pmatrix}\phi_\pm(0,p_\pm,s,h,\pm
a)\\\phi_{\pm,
z}(0,p_\pm,s,h,a_\pm)\end{pmatrix}=\begin{pmatrix}\phi_+(0,\cdot)-\phi_-(0,\cdot)\\\phi_{+,
z}(0,\cdot)-\phi_{-, z}(0,\cdot)\end{pmatrix}(p,s,h,a)\in\bR^{2N},
\end{align}
and note that since $\uw(z)$ is a connection, we have
\begin{align}\label{b26}
\Psi(\up,0,0,\ua)=0.
\end{align}

We shall add a nonhomogeneous boundary condition to the one in
\eqref{t3} in order to remove the translational indeterminacy
present in the shock case. Anticipating the later low frequency
analysis, we consider the augmented boundary condition
\begin{align}\label{b27}
\tilde{\Gamma}\begin{pmatrix}u_\pm\\v_\pm\\s\\h\end{pmatrix}=\begin{pmatrix}u_+-u_-\\v_+-v_-\\s+u_+\cdot\uw_z-\uw\cdot\uw_z\end{pmatrix}\text{
on }z=0.
\end{align}
Parallel to \eqref{b25} define
\begin{align}\label{b28}
\begin{split}
&\tilde{\Psi}(p,s,h,a)=\\
&\tilde{\Gamma}\begin{pmatrix}\phi_\pm(0,p_\pm,s,h,a_\pm)\\\phi_{\pm,
z}(0,p_\pm,s,h,a_\pm)\\s\\h\end{pmatrix}=\begin{pmatrix}\phi_+(0,\cdot)-\phi_-(0,\cdot)\\\phi_{+,
z}(0,\cdot)-\phi_{-,
z}(0,\cdot)\\s+\phi_+(0,\cdot)\cdot\uw_z(0)-\uw(0)\cdot\uw_z(0)\end{pmatrix}(p,s,h,a)\in\bR^{2N+1}.
\end{split}
\end{align}
Observe that
\begin{align}\label{b29}
\tilde{\Psi}(\up,0,0,\ua)=0.
\end{align}

\begin{rem}\label{melnikov}
\textup{ The function $\tilde \Psi$ plays a role here similar to
that of the extended Melnikov separation function defined in
\cite{ZS} for the undercompressive conservative case. Namely, it
provides a convenient full rank (implicit) representation of the
defining function $\chi$.}
\end{rem}

\subsection{Rank conditions and the manifold $\cC_\cB$}\label{rank}

\textbf{\quad}With $q=(p_+,p_-,s,h)$ the next Proposition
%is an analogue of \cite{Me1}, Prop. 5.4.3 for
%Lax shocks and
gives conditions under which there exist connections $W(z,q)$  near
$\uw(z)=W(z,\uq)$ satisfying \eqref{t1} and
\begin{align}\label{b299}
\lim_{z\to\pm\infty}W(z,q)=p_\pm.
\end{align}

\begin{prop}[Connections near a given one]\label{b30}
\textbf{}

(1) Let $\uw(z)$ be a connection corresponding to
$(p,s,h)=(\up,0,0)$ and suppose that
\begin{align}\label{b31}
\begin{split}
&(a)\mathrm{rank}\;\nabla_a\tilde{\Psi}(\up,0,0,\ua)=\cN_-=N+1-k,\\
&(b)\mathrm{rank}\;\nabla_{a,p}\tilde{\Psi}(\up,0,0,\ua)=2N+1.
\end{split}
\end{align}
Then in a neighborhood $\cO\subset\bR^{2N}\times\bR^d$ of
$(\up,0,0)$ there is a smooth manifold $\cC_\cB$ of dimension
$d+N-k$ and there are smooth mappings $W_\pm(z,p,s,h)$ on $\{\pm
z\geq 0\}\times\cC_\cB$ such that $(W_\pm,W_{\pm, z})$ satisfies the
profile equation \eqref{t3}, the endstate condition \eqref{b13}, and
the boundary condition
\begin{align}\label{b32}
\tilde{\Gamma}\begin{pmatrix}W_\pm\\W_{\pm,
z}\\s\\h\end{pmatrix}=0\text{ on }z=0.
\end{align}
With $q:=(p,s,h)=(p_+,p_-,s,h)\in\cC_\cB$, if we set
\begin{align}\label{b36}
W(z,q)=\begin{cases}W_+(z,q),\;z\geq 0\\W_-(z,q),\;z\leq
0\end{cases},
\end{align}
then $W(z,q)$ satisfies \eqref{t1}, \eqref{b299}.

The manifold $\cC_\cB$ can be defined by a \emph{generalized
Rankine-Hugoniot condition}
\begin{align}\label{b32e}
\chi(p,s,h)=0
\end{align}
for a smooth function $\chi:\cO\to \bR^{N+k}$ such that
$\mathrm{rank}\;\chi_p(\up,0,0)=N+k$.

(2)  The same conclusions hold if the second rank condition
\eqref{b31}(b) is replaced by
\begin{align}\label{b356}
\mathrm{rank}\;\nabla_{a,p,s}\tilde{\Psi}(\up,0,0,\ua)=2N+1.
\end{align}
The only difference is that now we have
$\mathrm{rank}\;\chi_{p,s}(\up,0,0)=N+k$.

\end{prop}

\begin{proof}
%The proof is like that of Prop. 5.4.3 (except for the appearance
%of $\Pi_-(0)$ in that proof, which I think is a cut and paste from
%the proof of Prop. 5.4.1 that shouldn't be there).
There is a reordering
$(p_\alpha,p_\beta)\in\bR^{N+k}\times\bR^{N-k}$ of the original
$(p_+,p_-)$ coordinates such that
$\nabla_{a,p_\alpha}\tilde{\Psi}(\up,0,0,\ua)$ is an isomorphism.
Applying the implicit function theorem to $\tilde{\Psi}(p,s,h,a)=0$
yields functions $p_\alpha(p_\beta,s,h)$, $a_\pm(p_\beta,s,h)$ and a
manifold $\cC_\cB$ parametrized in the new coordinates by
$(p_\alpha(p_\beta,s,h),p_\beta,s,h)$.   Reordering again yields
smooth functions $p_\pm(p_\beta,s,h)$ such that $\cC_\cB$ is given
in the original coordinates by
\begin{align}\label{b355}
\cC_\cB=\{(p_+(p_\beta,s,h),p_-(p_\beta,s,h),s,h):(p_\beta,s,h)\text{
near }(\up_\beta,0,0)\}.
\end{align}
We define
\begin{align}\label{b35}
W_\pm(z,p_+(p_\beta,s,h),p_-(p_\beta,s,h),s,h):=\phi_\pm(z,p_\pm(p_\beta,s,h),s,h,a_\pm(p_\beta,s,h)).
\end{align}
A generalized RH condition is given by the defining equation for
$\cC_\cB$ in $(p,s,h)$ space:
\begin{align}\label{b34e}
\chi(p,s,h):=p_\alpha-p_\alpha(p_\beta,s,h)=0.
\end{align}

The proof of part (2) is essentially the same.

\end{proof}

\begin{rem}\label{b34}
\textbf{1. }  \textup{Using the condition \eqref{b34e} we formulate
later a nonconservative analogue of the curved inviscid shock
problem and construct curved multidimensional nonconservative
shocks.}

\textbf{2. } \textup{When $k=0$, if in place of \eqref{b31}(b) we
assume
\begin{align}\label{b36z}
\mathrm{rank}\;\nabla_{a,p_+}\tilde{\Psi}(\up,0,0,\ua)=2N+1,
\end{align}
then we can prescribe $(p_-,s,h)$ and solve for $p_+$ as in the
conservative case.}

%\textbf{3. }  \textup{The factors $b(h)$ in \eqref{t3} affect the
%form of the function $p_\alpha(p_\beta,s,h)$, giving the expected
%influence of the viscosity on the form of the generalized RH
%condition.}
\end{rem}

We proceed to restate the hypotheses \eqref{b31} equivalently in
terms of rank conditions on $\Psi$.  We expect to be able to do
this, since $\tilde{\Psi}$ is just designed to remove the
translational indeterminacy left by $\Psi$.   A first step is to
relate such hypotheses to properties of solutions of the
linearized problem.

Consider the linearization of \eqref{t3}(a) with respect to $(u,v)$
at $(\uw,\uw')$, $(s,h)=(0,0)$.  Write this homogeneous $2N\times
2N$ linearized system  as
\begin{align}\label{b37}
\cL\begin{pmatrix}\dot{u}\\\dot{v}\end{pmatrix}:=\begin{pmatrix}\dot{u}\\\dot{v}\end{pmatrix}'-\cG(z)\begin{pmatrix}\dot{u}\\\dot{v}\end{pmatrix}=0\text{
on }\pm z\geq 0,
\end{align}
where
\begin{align}\label{b37b}
\cG(z)=\begin{pmatrix}0&I\\O(e^{-\delta|z|})&A_d(\uw)\end{pmatrix}.
\end{align}
The following conjugation lemma shows that solutions
$U:=(\dou,\dov)$ of \eqref{b37} can be conjugated to solutions $V$
of the limiting constant-coefficient systems
\begin{align}\label{b37c}
\partial_z V=\cG_\ell(\up_\pm)V, \text{ where
}\cG_\ell(\up_\pm):=\lim_{z\to\pm\infty}\cG(z)=\begin{pmatrix}0&I\\0&A_d(\up_\pm)\end{pmatrix}.
\end{align}

\begin{lem}[\cite{MZ3}, Lemma 2.6]\label{b37d}
For $\delta>0$ as in \eqref{a1}, there exist $2N\times 2N$ matrices
$Y_\pm(z)$ on $\pm z\geq 0$ and positive constants $C$,
$\delta'<\delta$ such that

(i)\;$Y_\pm$ and $Y^{-1}_\pm$ are $C^\infty$ and bounded with
bounded derivatives,

(ii)\;$|Y_\pm(z)-I|+|\partial_z Y_\pm(z)|\leq C e^{-\delta'|z|}$,

(iii)\;$Y_\pm$ satisfy
\begin{align}\label{b37e}
\partial_z Y_\pm=\cG(z)Y_\pm-Y_\pm\cG_\ell(\up_\pm)\text{ on }\pm z\geq 0.
\end{align}

\end{lem}
Observe that $U$ satisfies \eqref{b37} on $\pm z\geq 0$ if and only
if $V$ defined by $U=YV$ satisfies \eqref{b37c} on $\pm z\geq 0$.

Let $\cS_\pm$ denote the space of bounded solutions of \eqref{b37}
in $\pm z\geq 0$.  We'll refer to
\begin{align}\label{bz2}
\cS=\cS_+\times \cS_-
\end{align}
as the space of bounded solutions of \eqref{b37}.    Similarly, let
$\cS^0_\pm$ denote the space of solutions of \eqref{b37} that decay
to $0$ as $z\to\pm\infty$ and set $\cS_0=\cS^0_+\times\cS^0_-$.

\begin{prop}\label{bz1}
(a)\;Let $R_-$, $L_+$  be as in Definition \ref{b19}.  The
dimensions of $\cS_\pm$ are $N+R_-$ and $N+L_+$ respectively. Thus,
$\dim\cS=2N+(N+1-k)$.

(b)\;The dimensions of $\cS^0_\pm$ are $R_-$ and $L_+$ respectively.
Thus, $\dim \cS_0=N+1-k$.

\end{prop}

\begin{proof}
Using the conjugators $Y_\pm$ we can obtain the Proposition
immediately by proving the analogous statements for solutions
$V=(v_1,v_2)$ of the limiting problem \eqref{b37c}.   Solutions in
$z\geq 0$ are given by
\begin{align}\label{bz4}
V_+=\begin{pmatrix}v_{1+}\\v_{2+}\end{pmatrix}=\begin{pmatrix}v_{1+}(0)+(e^{zA_d(\up_+)}-I)A^{-1}_d(\up_+)v_{2+}(0)\\e^{zA_d(\up_+)}v_{2+}(0)\end{pmatrix}.
\end{align}
Clearly, $V_+$ is bounded in $z\geq 0$ if and only if
$v_{2+}(0)\in\bE_-(A_d(\up_+))$ (with $v_{1+}(0)\in\bR^N$
arbitrary), and $V_{+}$ decays to $0$ as $z\to +\infty$ if and only
if
\begin{align}\label{bz5}
v_{2+}(0)\in\bE_-(A_d(\up_+))\text{ and
}v_{1+}(0)-A^{-1}_d(\up_+)v_{2+}(0)=0.
\end{align}
The ``$-$" case is similar.

\end{proof}

In the next Proposition  we consider rank conditions on $\Psi$ as in
\eqref{b25}. With $\Gamma$ as in \eqref{t3} observe that
\begin{align}\label{b37a}
(\uw',\uw'')
\end{align}
satisfies the problem
\begin{align}\label{b38}
\cL\begin{pmatrix}\dot{u}\\\dot{v}\end{pmatrix}=0\text{ on }\pm
z\geq 0,
\;\Gamma\begin{pmatrix}\dot{u}_\pm\\\dot{v}_\pm\end{pmatrix}=0\text{
on }z=0,\; \lim_{z\to \pm\infty}\dot{u}=0.
\end{align}

\begin{prop}\label{b39}
(a)The condition
\begin{align}\label{b40}
\mathrm{rank}\nabla_a\Psi(\up,0,0,\ua)=\cN_--1=N-k
\end{align}
holds if and only if the problem \eqref{b38} has a one dimensional
kernel spanned by \eqref{b37a}.

(b)  The condition
\begin{align}\label{b41}
\mathrm{rank}\nabla_{a,p}\Psi(\up,0,0,\ua)=2N
\end{align}
holds if and only if for all $g\in\bR^{2N}$ the problem
\begin{align}\label{b42}
\cL\begin{pmatrix}\dot{u}\\\dot{v}\end{pmatrix}=0,
\;\Gamma\begin{pmatrix}\dot{u}_\pm\\\dot{v}_\pm\end{pmatrix}=g
\end{align}
has a bounded solution.

\end{prop}

\begin{proof}
Proposition \ref{bz1} shows that  the space of $\cS$ of bounded
solutions of $\cL\dot{U}=0$ has dimension $2N+(N+1-k)$, while
$\cS_0$ has dimension $N+1-k$. Thus, Corollary \ref{bp1} implies
that the map
\begin{align}\label{b43}
(\dot{p},\dot{a})\to
\begin{pmatrix}\nabla_a\phi(z,\up,0,0,\ua)\dot{a}+
\nabla_p\phi(z,\up,0,0,\ua)\dot{p}\\\nabla_a\phi_z(z,\up,0,0,\ua)\dot{a}+
\nabla_p\phi_z(z,\up,0,0,\ua)\dot{p}\end{pmatrix}
\end{align}
is an isomorphism of $\bR^{2N}\times \bE_-(G_d(\up,0,0))$ onto
$\cS$, and that
\begin{align}\label{b44}
\dot{a}\to\begin{pmatrix}\nabla_a\phi(z,\up,0,0,\ua)\dot{a}\\\nabla_a\phi_z(z,\up,0,0,\ua)\dot{a}\end{pmatrix}
\end{align}
is an isomorphism onto $\cS_0$.  It follows that the dimension of
the space of solutions of \eqref{b38} is the dimension of
\begin{align}\label{b45}
\{\dot{a}\in\bE_-(G_d(\up,0,0)):\nabla_a\Psi(\up,0,0,\ua)\dot{a}=0\}.
\end{align}
This dimension is one if and only if
$\mathrm{rank}\nabla_a\Psi(\up,0,0,\ua)=N-k$.

Similarly, the map $\dot{U}\to \Gamma\dot{U}_\pm(0)$ from $\cS$ to
$\bR^{2N}$ is onto if and only if
\begin{align}\label{b46}
\nabla_{a,p}\Psi(\up,0,0,\ua)
\end{align}
has rank $2N$.

\end{proof}

\begin{prop}\label{b47}
\begin{align}\notag
\begin{split}
&(a)\mathrm{rank}\;\nabla_a\tilde{\Psi}(\up,0,0,\ua)=N+1-k\Leftrightarrow\mathrm{rank}\;\nabla_a\Psi(\up,0,0,\ua)=N-k\\
&(b)\mathrm{rank}\;\nabla_{a,p}\tilde{\Psi}(\up,0,0,\ua)=2N+1\Leftrightarrow\mathrm{rank}\;\nabla_{a,p}\Psi(\up,0,0,\ua)=2N\\
&(c)\mathrm{rank}\;\nabla_{a,p,s}\tilde{\Psi}(\up,0,0,\ua)=2N+1\Leftrightarrow\mathrm{rank}\;\nabla_{a,p,s}\Psi(\up,0,0,\ua)=2N.
\end{split}
\end{align}
\end{prop}

\begin{proof}
\textbf{1}.  Using the isomorphism \eqref{b44} we find an element
$\underline{\dot{a}}\in\bE_-(G_d(\up,0,0))$ such that
\begin{align}\label{b48}
\nabla_a\phi(z,\up,0,0,\ua)\underline{\dot{a}}=\uw_z(z).
\end{align}
We always have
\begin{align}\label{b49}
\underline{\dot{a}}\in\mathrm{ker}\nabla_a\Psi(\up,0,0,\ua).
\end{align}

\textbf{2}.  Statement (a) is equivalent (suppressing evaluation
at $(\up,0,0,\ua)$) to
\begin{align}\label{b50}
\mathrm{dim\; ker}\nabla_a\Psi=1\Leftrightarrow\mathrm{dim\;
ker}\nabla_a\tilde{\Psi}=0.
\end{align}
Also observe from \eqref{b28} that
\begin{align}\label{b51}
\dot{a}\in\mathrm{ker
}\nabla_a\tilde{\Psi}\Leftrightarrow\dot{a}\in\mathrm{ker
}\nabla_a\Psi\text{ and }
(\nabla_a\phi_+\;\dot{a})\cdot\uw_z(0)=0.
\end{align}
Now, if $\mathrm{dim\; ker}\nabla_a\Psi=1$, then $\mathrm{ker
}\nabla_a\Psi$ is spanned by $\underline{\dot{a}}$.  But by
\eqref{b48}
\begin{align}\label{b52}
\nabla_a\phi_+(0,\up,0,0,\ua)\underline{\dot{a}}=\uw_z(0),
\end{align}
so $\mathrm{dim\; ker}\nabla_a\tilde{\Psi}=0$.

Suppose $\mathrm{dim\; ker}\nabla_a\tilde{\Psi}=0$.   Then by
\eqref{b51} the linear functional defined on $\mathrm{ker
}\nabla_a\Psi$ by
\begin{align}\label{b53}
T\dot{a}=(\nabla_a\phi_+\;\dot{a})\cdot\uw_z(0)
\end{align}
satisfies $\mathrm{dim\;ker}\;T=0$. The rank of $T$ is one, so
$\mathrm{dim\; ker}\nabla_a\Psi=1$.   This proves (a).

\textbf{3. }Statement (b) is equivalent to
\begin{align}\label{b54}
\mathrm{dim\; ker}\nabla_{a,p}\Psi=N+1-k\Leftrightarrow\mathrm{dim\;
ker}\nabla_{a,p}\tilde{\Psi}=N-k.
\end{align}
Observe that $\ker\nabla_{a,p}\tilde{\Psi}\subset\ker
\nabla_{a,p}\Psi$ and
\begin{align}\label{b55}
(\dot{a},\dot{p})\in\mathrm{ker
}\nabla_{a,p}\tilde{\Psi}\Leftrightarrow(\dot{a},\dot{p})\in\mathrm{ker
}\nabla_{a,p}\Psi\text{ and }
(\nabla_a\phi_+\;\dot{a}+\nabla_p\phi_+\;\dot{p})\cdot\uw_z(0)=0.
\end{align}
We have
\begin{align}\label{b56}
(\underline{\dot{a}},0)\in\mathrm{ker}\;\nabla_{a,p}\Psi\setminus\mathrm{ker}\;\nabla_{a,p}\tilde{\Psi}.
\end{align}

Suppose $\mathrm{dim\; ker}\nabla_{a,p}\Psi=N+1-k$ and let a basis
for that kernel be $\{(\underline{\dot{a}},0),e_1,\dots,e_{N-k}\}$.
Using \eqref{b52} and \eqref{b55} we obtain a basis for
$\ker\nabla_{a,p}\tilde{\Psi}$ of the form
$\{e_i+s_i(\underline{\dot{a}},0):i=1,\dots,N-k\}$ for appropriate
$s_i\in\bR$.

Suppose $\mathrm{dim\; ker}\nabla_{a,p}\tilde{\Psi}=N-k$.  By
\eqref{b55} the functional defined on $\mathrm{
ker}\nabla_{a,p}\Psi$ by
\begin{align}\label{b56t}
S(\doa,\dop)=(\nabla_a\phi_+\;\dot{a}+\nabla_p\phi_+\;\dot{p})\cdot\uw_z(0)
\end{align}
has rank one and kernel of dimension $N-k$.  So $\mathrm{
ker}\nabla_{a,p}\Psi$ has dimension $N+1-k$.

The proof of (c) is similar to that of (b).
\end{proof}

\begin{defn}\label{c11}
1.  When both of the rank conditions
\begin{align}\label{ct2}
\begin{split}
&\mathrm{rank}\nabla_a\tilde{\Psi}(\up,0,0,\ua)=N+1-k\;\;(a-transversality)\\
&\mathrm{rank}\nabla_{a,p,s}\tilde{\Psi}(\up,0,0,\ua)=2N+1\;\;((a,p,s)-transversality)
\end{split}
\end{align}
are satisfied, we say the profile $\uw(z)$ is \emph{transversal}.

2. The profile $\uw$ is said to be \emph{strongly transversal} if
both  a-transversality and
\begin{align}\label{ct3}
\mathrm{rank}\nabla_{a,p}\tilde{\Psi}(\up,0,0,\ua)=2N+1\;\;((a,p)-transversality)
\end{align}
are satisfied.

3.  For $q=(p,s,h)\in\cC$ the profile $W(z,q)$ is transversal
(resp., strongly transversal) if conditions \eqref{ct2} (resp.,
(a) and (b) of Prop. \ref{b47}) hold with $(\up,0,0,\ua)$ replaced
by $(p,s,h,a)$, where $a=a(p_\beta,s,h)$ is as in \eqref{b35}.
\end{defn}

\begin{rem}\label{b56y}

\textbf{1. }\textup{Transversality (resp., strong transversality) of
$\uw$ implies transversality (resp., strong transversality) of
$W(z,q)$ for $q$ near $\uq$.}

\textbf{2. Geometric transversality. }  \textup{Consider the
first-order system \eqref{t3} as a system on $\bR$ for the moment.
Let $W_s$, $W_{cs}$ be the stable and center-stable manifolds of the
rest point $(\up_+,0)$, and let $W_u$, $W_{cu}$ be the unstable and
center-unstable manifolds of the rest point $(\up_-,0)$.  Then
a-transversality corresponds to the statement $W_s\pitchfork W_u$ at
$\begin{pmatrix}\uw(0),\uw_z(0)\end{pmatrix}$, while the combination
of a-transversality and (a,p)-transversality (which we call
\emph{strong transversality}) corresponds to $W_{cs}\pitchfork
W_{cu}$ at $\begin{pmatrix}\uw(0),\uw_z(0)\end{pmatrix}$}.

\textup{Next consider the system obtained by augmenting \eqref{t3}
with the equation $s'=0$, and let $W^a_{cs}$  (resp., $W^a_{cu}$) be
the center-stable (resp., center-unstable)  manifold of the rest
point $(\up_+,0,0)$ (resp., $(\up_-,0,0)$) for the augmented system.
Then the combination of a-transversality and (a,p,s)-transversality
(which we call \emph{transversality}) corresponds to
$W^a_{cs}\pitchfork W^a_{cu}$ at
$\begin{pmatrix}\uw(0),\uw_z(0),0\end{pmatrix}$.}

%\textbf{2. }It is easy to check that the analogue of Proposition
%\ref{b47} still holds if the second rank hypothesis \eqref{b31}(b)
%is replaced by the weaker hypothesis
%\begin{align}\label{be56}
%\mathrm{rank} \nabla_{a,p,s}\tilde{\Psi}(\up,0,0,\ua)=2N+1.
%\end{align}
%The same applies to Proposition \ref{b47} as well.   This remark
%is important, for example, in the proof of Theorem \ref{K}.

\textbf{3. }\textup{One could apply the implicit function theorem
directly to the equation
\[
\Psi(p,s,h,a)=0
\]
using the rank conditions on $\nabla_a\Psi$ and $\nabla_{a,p}\Psi$
in Prop. \ref{b47}.   But then instead of obtaining a function
$p_\alpha(p_\beta,s,h)$ as in \eqref{b34e}, we'd obtain a function
$p_\alpha(p_\beta,s,h,a_i)$, where $a_i$ is one of the $a$
components.   This is not the form a generalized RH condition should
have, since $a_i$ does not correspond to any of the unknowns in the
inviscid hyperbolic problem; so we use the extra boundary condition
of $\tilde{\Psi}$ instead.}

%\textbf{. Small amplitude undercompressive shocks.} One might try
%to construct small amplitude undercompressive shock profiles using
%the argument for small amplitude boundary layers in Prop. 5.4.1 of
%Guy's book.   The needed rank conditions are satisfied, but now
%the implicit function theorem yields only zero amplitude profiles
%($p_+=p_-$) near a given noncharacteristic constant state
%($\up_+=\up_-$).

%\textbf{4. }  \textup{See Schechter \cite{S} for a construction of
%small amplitude nonconservative shocks by center manifold
%reduction.}

\end{rem}

\section{Linearized GRH derived from transmission
conditions}\label{linGRH}

\textbf{\;\;}In this section we derive more explicit
characterizations of $T_q\cC_\cB$ and the linearized GRH conditions
that are useful in the later low frequency analysis of the standard
and modified Evans functions.

\subsection{Linearized hyperbolic problem}\label{linhyp}

\textbf{\;\;} Suppose we are given a shock manifold $\cC$ as in
Assumption \ref{i3} and a planar shock
$\uq=(\up_+,\up_-,0,0)\in\cC$. To construct a curved nonconservative
shock as a perturbation of $\uq$, we solve the transmission problem
\begin{align}\label{c1}
\begin{split}
&\sum^{d-1}_{j=0}A_j(u)\partial_ju+\cA_d(u,d\psi)\partial_du=0\\
&\chi(u,d\psi)=0\text{ on }x_d=0
\end{split}
\end{align}
where $u$ is near $\up$, $d\psi$ is near $0$, and $\chi(p,s,h)$
\eqref{b34e} is the defining function for $\cC$ near $\uq$.

To solve this nonlinear problem we need to solve linearized problems
at $q=(p,s,h)$ near $\uq$.  The (fully) linearized problem at
$q=(p_+,p_-,s,h)$ is
\begin{align}\label{c2}
\begin{split}
&(a)\sum^{d-1}_{j=0}A_j(p)\partial_j\dou+\cA_d(p,s,h)\partial_d\dou=f\text{ on }\pm x_d\geq 0\\
&(b)\chi'(p,s,h)(\dou,d\dopsi)=g\text{ on }x_d=0.
\end{split}
\end{align}
The $p$ in \eqref{c2}(a) should be understood as $p_\pm$, while that
in \eqref{c2}(b) should be understood as $(p_+,p_-)$. Similar
notation is sometimes used below.  When $g=0$, the boundary
condition in \eqref{c2}  is the requirement
\begin{align}
(\dou,d\dopsi)\in T_q\cC.
\end{align}

In the low frequency analysis of viscous stability determinants in
the case when $\cC=\cC_\cB$, we'll sometimes require a more explicit
form of the linearized boundary condition, namely,
\begin{align}\label{c3}
\Gamma_{0,red}(q)(\dou,d\dopsi)=0,
\end{align}
where $\Gamma_{0,red}(q)$ is defined in \eqref{c22}.  We construct
$\Gamma_{0,red}(q)$ by considering the full linearization with
respect to $w$, $s$, and $h$ of the profile transmission problem for
$w(z)$, where  $w=w_\pm$ on $\pm z\geq 0$.
\begin{align}\label{c4}
\begin{split}
&\cA_d(w,s,h)w'=(1+|h|^2)w''\text{ in }\pm z\geq 0\\
&\Gamma_1(w,w_z,s,h):=\begin{pmatrix}[w]\\
[w_z]\\s+w_+\cdot\uw_z-\uw_z\cdot\uw\end{pmatrix}=0\text{ on }z=0.
\end{split}
\end{align}

\subsection{Linearized parabolic problem and transversality
}\label{linpar}

 \textbf{\quad}Let us write $W(z,q)$ for the profile on $\bR$ given by
\eqref{b36}; thus $\uw(z)=W(z,\uq)$.   The full linearization of
\eqref{c4} at $(W(z,q),s,h)$ is
\begin{align}\label{c5}
\begin{split}
&(a)\cL_0(z,q,\partial_z)\dow:=-(1+|h|^2)\dow''+\cA_d(W,s,h)\dow'+\partial_w\cA_d(W,s,h)\dow W'=\\
&\qquad\qquad\dos A_0(W)W'+\sum_{j=1}^{d-1}\doh_jA_j(W)W'+2\sum^{d-1}_{j=1}h_j\doh_jW'':=\cL_{0,1}(z,q)(\dos,\doh)\\
&(b)\Gamma_2(\dow,\dow_z,\dos,\doh):=\begin{pmatrix}[\dow]\\
[\dow_z]\\\dos+\dow_+\cdot\uw_z\end{pmatrix}=0\text{ on }z=0
\end{split}
\end{align}

Before studying \eqref{c5} and in order to avoid later repetitions,
it is desirable at this point to consider the full linearization of
the interior parabolic problem in \eqref{i22}
\begin{align}\label{c5aa}
\sum^{d-1}_{j=0}A_j(u)\partial_j
u+\cA_d(u,d\psi)\partial_zu-\sum^d_{j=1}(\partial_j-\partial_j\psi\partial_z)^2u=0\text{
on }\pm z\geq 0
\end{align}
about the exact solution given by a profile $W(z,q)$ and front
$\psi=st+hy$.  The Laplace-Fourier transform in $(t,y)$ of that
linearization is
\begin{align}\label{c5b}
\cL(z,q,\zeta,\partial_z)\dou-\dopsi\cL_{1}(z,q,\zeta)=f\text{ in
}\pm z\geq 0.
\end{align}
Here $\zeta=(\tau,\gamma,\eta)$ where $\gamma >0$ and, with
$W'=W_z(z,q)$, the operators $\cL$ and $\cL_1$ are given explicitly
by
\begin{align}\label{c5a}
\begin{split}
&\cL(z,q,\zeta,\partial_z)\dou=-(1+|h|^2)\dou_{zz}+\left(\cA_d(W,s,h)+2\sum^{d-1}_{j=1}h_ji\eta_j\right)\dou_z+\partial_w\cA_d(W,s,h)\dou W'+\\
&\qquad\qquad\qquad A_0(W)(i\tau+\gamma)\dou+\sum^{d-1}_{j=1}A_j(W)i\eta_j\dou+|\eta|^2\dou,\\
&\cL_1(z,q,\zeta)=A_0(W)W'(i\tau+\gamma)+\sum^{d-1}_{j=1}A_j(W)W'i\eta_j+2\sum^{d-1}_{j=1}h_ji\eta_jW''+|\eta|^2W'.
\end{split}
\end{align}
It is important in the sequel to know the relationship between the
operators $\cL_0$ and $\cL_{0,1}$ appearing in \eqref{c5} and the
operators $\cL$ and $\cL_1$ in \eqref{c5a}.    Introducing polar
coordinates
\begin{align}\label{c6a}
\zeta=\rho\hzeta, \;\rho=|\zeta|,\; \hzeta=(\htau,\hgamma,\heta)\in
S^d_+=S^d\cap\{\hgamma>0\},\;\oS^d_+=S^d\cap\{\hgamma\geq 0\}
\end{align}
and writing $\cL_1(z,q,\zeta)=\rho\chcL_1(z,q,\hzeta,\rho)$ we
clearly have
\begin{align}\label{c6b}
\cL_0(z,q,\partial_z)=\cL(z,q,0,\partial_z)\text{ and
}\cL_{0,1}(z,q)(i\htau+\hgamma,i\heta)=\chcL_1(z,q,\hzeta,0).
\end{align}

\subsubsection{Conjugation to limiting and block diagonal
systems}\label{conj}

In order to understand the behavior of solutions to \eqref{c5} and
\eqref{c5b}, it is helpful  first to rewrite the partially
linearized problem
\begin{align}\label{c6c}
\cL(z,q,\zeta,\partial_z)\dou=f
\end{align}
as a first order system and then conjugate it to simpler forms. Here
we'll also establish notation that will be used in the rest of the
paper.   Recall $q=(p_+,p_-,s,h)$, where $p_\pm\in\omega_\pm$ and
$|s,h|\leq\epsilon_\omega$ for $\omega_\pm$, $\epsilon_\omega$ as in
Prop. \ref{b17a}.    Set
\begin{align}\label{c6dd}
\cQ:=\omega_+\times\omega_-\times\{(s,h):|s,h|\leq
\epsilon_\omega\}.
\end{align}

With $U:=(\dou,\dou_z)^t$   we can rewrite \eqref{c6c} as a
$2N\times 2N$ first order system
\begin{align}\label{c6d}
\partial_zU=G(z,q,\zeta)U+F\text{ on }\pm z\geq 0,
\end{align}
where
\begin{align}\label{c6e}
\begin{split}
&F=(0,-b(h)f)^t,\;\;b(h)=(1+|h|^2)^{-1},\;\;G=\begin{pmatrix}0&I\\G_{21}&G_{22}\end{pmatrix}\text{ with }\\
&G_{21}(z,q,\zeta)\dou=b(h)\left(\partial_w\cA_d(W,s,h)\dou
W'+A_0(W)(i\tau+\gamma)\dou+\sum^{d-1}_{j=1}A_j(W)i\eta_j\dou+|\eta|^2\dou\right)\\
&G_{22}(z,q,\zeta)=b(h)\left(\cA_d(W,s,h)+2\sum^{d-1}_{j=1}h_ji\eta_j\right).
\end{split}
\end{align}
Recalling that $W(z,q)\to p_\pm$ as $z\to\pm\infty$, we define
limiting systems $G_\pm(q,\zeta)$ by replacing $W$ by $p_\pm$ and
$W'$ by $0$ in \eqref{c6e}:
\begin{align}\label{c6ee}
G_\pm(q,\zeta)=\begin{pmatrix}0&I\\G^{21}_\pm(q,\zeta)&G^{22}_\pm(q,\zeta)\end{pmatrix}.
\end{align}

The following conjugation lemma is a version of Lemma \ref{b37d}
with parameters.

\begin{lem}[\cite{MZ3},\;Lemma 2.6]
\label{lemconj} For all $q_0\in \cQ$ and $\zeta_0 \in
\obR^{1+d}_+=\bR^{d+1}\cap\{\gamma\geq 0\}$ there is a neighborhood
$\Omega$ of $(q_0, \zeta_0)$ in $\cQ \times \obR^{1+d}_+$ and there
are matrices $Y_\pm$ defined and $C^\infty$ on $\{\pm z \ge 0\}
\times \Omega$  such that:

\quad i) $Y_\pm $ and  $(Y_\pm)^{-1} $ are uniformly bounded and for
$\delta>0$ as in \eqref{a1} there are positive constants $C_\alpha$,
$\delta'<\delta$ such that for $(q, \zeta) \in \Omega$:
 \begin{equation*}
\big\vert  \D^\alpha _{z, q, \zeta}
 \big( Y_\pm(z, q, \zeta) - \Id   \big)
 \big\vert
  \leq C_\alpha  e^{ - \delta_1 \vert  z\vert }
\,\text{ on } \pm z \ge 0\,\text{ for all }\alpha;
\end{equation*}

\quad ii) $Y_+$ and $Y_-$ satisfy
\begin{equation*}
\D_z Y_\pm (z, q, \zeta) = G  (z, q, \zeta)Y_\pm (z, q, \zeta) -
Y_\pm(z, q, \zeta) G_\pm( q, \zeta) \text{ on }\pm z\geq 0.
\end{equation*}

\end{lem}

Observe that $U$ satisfies \eqref{c6d} if and only if $V$ defined by
$U_\pm=Y_\pm V_\pm$ satisfies
\begin{align}\label{c6f}
\partial_z V=G_\pm(q,\zeta)V+Y^{-1}F\text{ on }\pm z\geq 0.
\end{align}
As in \cite{GMWZ3} for $\rho=|\zeta|$ small we perform an additional
conjugation of \eqref{c6f} to block diagonal (or ``HP") form using
\begin{align}\label{c6g}
V_\pm=\Lambda_\pm(q,\zeta)\begin{pmatrix} u_{H\pm},\\
u_{P\pm}\end{pmatrix}
\end{align}
where $\Lambda_\pm$ is $C^\infty$ and
\begin{align}\label{c6h}
\Lambda_\pm(q,0)=\begin{pmatrix}I&(G^{22}_\pm)^{-1}\\0&I\end{pmatrix}.
\end{align}
This transforms \eqref{c6f} to
\begin{align}\label{c6i}
\partial_z\begin{pmatrix}u_{H\pm}\\u_{P\pm}\end{pmatrix}=\begin{pmatrix}H_\pm(q,\zeta)&0\\0&P_\pm(q,\zeta)\end{pmatrix}\begin{pmatrix}u_{H\pm}\\u_{P\pm}\end{pmatrix}+\Lambda^{-1}Y^{-1}F\text{ on }\pm z\geq 0,
\end{align}
where
\begin{align}\label{c6j}
\begin{split}
&P_\pm(q,\zeta)=G^{22}_\pm(q,\zeta)+O(\rho)\\
&H_\pm(q,\zeta)=-(G^{22}_\pm)^{-1}G^{21}_{\pm}+O(\rho^2)=\\
&\qquad\qquad-\cA_d(p_\pm,s,h)^{-1}\left(A_0(p_\pm)(i\tau+\gamma)+\sum^{d-1}_{j=1}A_j(p_\pm)i\eta_j\right)+O(\rho^2).
\end{split}
\end{align}

%&\tGamma_H(q,\zeta)u_H+\tGamma_P(q,\zeta)u_P=0\text{ on }z=0.

\subsubsection{Other characterizations of
transversality.}\label{other}

Next we give characterizations of transversality that are useful in
formulating the reduced transmission conditions.    The first step
is to rewrite the fully linearized profile problem \eqref{c5} in a
form with modified transmission conditions and where $\cL_{0,1}$ no
longer appears.

\begin{prop}\label{c6k}

There exist $C^\infty$  functions $\cR_\pm(z,q,\dos,\doh)$ valued in
$\bR^N$ and  satisfying
\begin{align}\label{c6}
\begin{split}
&(a)\;\cL_0(z,q,\partial_z)\cR(z,q,\dos,\doh)=\cL_{0,1}(z,q)(\dos,\doh)\text{ on }\pm z\geq 0\\
&(b)\;\cR_\pm(z,q,\dos,\doh)\text{ is linear in }(\dos,\doh)\\
&(c)\;\uw_z(0)\cdot\cR_\pm(0,q,\dos,\doh)=-\dos\\
&(d)\;\nabla_{\dos,\doh}\cR_\pm=O(e^{-\delta |z|}) \text{ for some
}\delta>0.
\end{split}
\end{align}

\end{prop}

\begin{proof}
\textbf{1. }We first construct $\cR_1(z,q,\dos,\doh)$ such that
\begin{align}\label{c6l}
\cL_0(z,q,\partial_z)\cR_1(z,q,\dos,\doh)=\cL_{0,1}(z,q)(\dos,\doh)\text{
on }\pm z\geq 0,
\end{align}
by setting
\begin{align}\label{c6m}
\cR_1(z,q,\dos,\doh):=\cR_{10}(z,q)\dos+\sum^{d-1}_{j=1}\cR_{1j}(z,q)\doh_j,
\end{align}
where the $\cR_{1k}$ are $C^\infty$, exponentially decaying as
$z\to\pm\infty$, and satisfy
\begin{align}\label{c6n}
\begin{split}
&\cL_0(z,q,\partial_z)\cR^\pm_{10}=A_0(W)W'\text{ on }\pm z\geq 0\\
&\cL_0(z,q,\partial_z)\cR^\pm_{1j}=A_j(W)W'+2h_jW''.
\end{split}
\end{align}
Functions $\cR_{1k}$ as above may be constructed, for example, by
rewriting the equations \eqref{c6n} as first order systems and then
conjugating by $T_\pm(z,q,0)$, where
\begin{align}\label{c6o}
T_\pm(z,q,\zeta):=Y_\pm(z,q,\zeta)\Lambda_\pm(q,\zeta)
\end{align}
(recall \eqref{c6b}, Lemma \ref{lemconj}, and \eqref{c6g}).     This
yields easily solvable systems of the form
\begin{align}\label{c6p}
\partial_z\begin{pmatrix}u_{H\pm}\\u_{P\pm}\end{pmatrix}=\begin{pmatrix}0&0\\0&P_\pm(q,0)\end{pmatrix}\begin{pmatrix}u_{H\pm}\\u_{P\pm}\end{pmatrix}+F_{k\pm}
\end{align}
with $F_{k\pm}$ exponentially decaying, where
\begin{align}\label{c6q}
\begin{pmatrix}\cR_{1k}\\\partial_z\cR_{1k}\end{pmatrix}(z,q)=T_\pm(z,q,0)\begin{pmatrix}u_{H\pm}\\u_{P\pm}\end{pmatrix}(z,q)\text{ on }\pm z\geq 0.
\end{align}

\textbf{2. }Since
\begin{align}\label{c6q2}
\cL_0(z,q,\partial_z)W_z(z,q)=0,
\end{align}
we can arrange the transmission condition \eqref{c6}(c) by setting
\begin{align}\label{c6r}
\cR(z,q,\dos,\doh):=\cR_1(z,q,\dos,\doh)-\left(\frac{\uw_z(0)\cdot\cR_1(0,q,\dos,\doh)+\dos}{\uw_z(0)\cdot
W_z(0,q)}\right) W_z(z,q).
\end{align}
\end{proof}

\begin{rem}\label{c7}
%1.  The construction of $\cR$ does not use conservative structure.
%However, the proof of the ZS result in GMWZ3 did use conservative
%structure.  One advantage of the strategy in \cite{FF} is that one
%can prove a ZS result without conservative structure.
%TODO: discuss need for proofs without cons. structure in intro.
  \textup{Although $\cR$ is initially defined only for
$(\dos,\doh)\in\bR^d$, it is important for later applications to
note that it extends immediately to $(\dos,\doh)\in\bC^d$ (for
example, see \eqref{c22} and \eqref{c37z}).}

\end{rem}

An immediate consequence of Proposition \ref{c6k} is:
\begin{cor}\label{c6s}
Define
\begin{align}\label{c8}
\dov:=\dow-\cR(z,q,\dos,\doh).
\end{align}
Then $(\dow,\dos,\doh)$ satisfies \eqref{c5} if and only if
$(\dov,\dos,\doh)$ satisfies
\begin{align}\label{c9}
\begin{split}
&(a)\;\cL_0(z,q,\partial_z)\dov=0\text{ on }\pm z\geq 0\\
&(b)\;\Gamma_3(q)(\dov,\dov_z,\dos,\doh):=\begin{pmatrix}[\dov]+[\cR(z,q,\dos,\doh]\\
[\dov_z]+[\cR_z(z,q,\dos,\doh)]\\\dov_+\cdot\uw_z\end{pmatrix}=0\text{
on }z=0.
\end{split}
\end{align}
\end{cor}

Let us write $T(z,q,\zeta)$  in \eqref{c6o} as (suppressing $\pm$)
\begin{align}\label{c9a}
T(z,q,\zeta)=\begin{pmatrix}T_{11}&T_{12}\\T_{21}&T_{22}\end{pmatrix}(z,q,\zeta).
\end{align}

\begin{rem}\label{c9b}

\textup{By conjugation to  \eqref{c6p} with $F_\pm=0$ we see that
solutions of \eqref{c9}(a) can be written as
%(recall $W(z,q)\to p_\pm$ as $z\to\pm\infty$)
\begin{align}\label{c10}
\dov(z)=T_{11}(z,q,0)u_{H}+T_{12}(z,q,0)e^{zP(q,0)}u_{P}\text{ on
}\pm z\geq 0,
\end{align}
where $u_{H\pm}\in\bR^N$, $u_{P\pm}\in\bR^N$, and
$T_{11\pm}(z,q,0)\to I$ as $z\to\pm\infty$.   A solution $\dov$ is
bounded if and only if $u_{P\pm}\in\bE_\mp(P_\pm(q,0))$.  Thus, the
space of bounded solutions $\cS(q)$ of \eqref{c9}(a) has dimension
$2N+(N+1-k)$ and the space $\cS_0(q)$ of solutions tending to zero
as $z \to\pm\infty$ has dimension $R_-+L_+=N+1-k$.}
\end{rem}

  A small modification of the proof of Prop. \ref{b39} (replace
$\Psi$ by $\tilde{\Psi}$) now yields:

\begin{prop}\label{c10a}
The property of a-transversality holds for $W(z,q)$ (Definition
\ref{c11}) if and only if there are no nontrivial solutions of
\eqref{c9} which tend to zero as $z\to\pm\infty$ when
$(\dos,\doh)=0$. Similarly, (a,p)-transversality holds if and only
if for all $g\in\bR^{2N+1}$ the problem
\begin{align}\label{c13}
\begin{split}
&(a)\;\cL_0(z,q,\partial_z)\dov=0\text{ on }\pm z\geq 0\\
&(b)\;\Gamma_3(q)(\dov,\dov_z,0,0):=g\text{ on }z=0
\end{split}
\end{align}
has a bounded solution.

\end{prop}

We can use \eqref{c10} to write the boundary condition \eqref{c9}(b)
as
\begin{align}\label{c14}
\Gamma_{0,H}(q)u_H+\Gamma_{0,P}(q)u_P+\Gamma_\cR(\dos,\doh)=\Gamma_3(q)(\dov,\dov_z,\dos,\doh)=0,
\end{align}
where
\begin{align}\label{c15}
\begin{split}
&\Gamma_{0,H}(q)u_H=\begin{pmatrix}\begin{bmatrix}T_{11}(0,q,0)u_H\\T_{21}(0,q,0)u_H\end{bmatrix}\\(T_{11+}(0,q,0)u_{H+})\cdot\uw_z(0)\end{pmatrix},
\;\Gamma_{0,P}(q)u_P=\begin{pmatrix}\begin{bmatrix}T_{12}(0,q,0)u_P\\T_{22}(0,q,0)u_P\end{bmatrix}\\(T_{12+}(0,q,0)u_{P+})\cdot\uw_z(0)\end{pmatrix},\\
&\text{ and }\Gamma_\cR(\dos,\doh):=\begin{pmatrix}[\cR(0,q,\dos,\doh]\\
[\cR_z(0,q,\dos,\doh)]\\0\end{pmatrix}.
\end{split}
\end{align}

With  \eqref{c10} and \eqref{c14} we obtain directly the following
rephrasing of Prop. \ref{c10a}:
\begin{prop}\label{c16t}
With $P_{0\pm}(q):=P(q,0)=b(h)\cA_d(p_\pm,s,h)$,  a-transversality
holds for $W(z,q)$ if and only if
\begin{align}\label{c17}
\ker\Gamma_{0,P}(q)\cap
\left(\bE_-(P_{0+}(q))\times\bE_+(P_{0-}(q))\right)=\{0\},
\end{align}
while (a,p)-transversality means that the rank of the matrix
\begin{align}\label{c18}
(\Gamma_{0,H}(q),\Gamma_{0,P}(q),\Gamma_{\cR}):\bR^{2N}\times\left(\bE_-(P_{0+}(q))\times\bE_+(P_{0-}(q))\right)\times
\bR^d\to\bR^{2N+1}
\end{align}
is $2N+1$ when $(\dos,\doh)=0$.

Equivalent characterizations of a-transversality and
(a,p)-transversality are obtained by replacing each space appearing
in \eqref{c17} or \eqref{c18} by its complexification.

\end{prop}

In making the last assertion of the Proposition we used the
observations made in Notation \ref{b13z} together with the extension
of $\cR$ to $(\dos,\doh)\in\bC^d$.

\subsection{Reduced transmission conditions and
$T_q\cC_\cB$}\label{TC}

\textbf{\quad}The  reduced boundary operator constructed in this
section will play a key role in the later stability analyses.

If a-transversality holds there is a decomposition
\begin{align}\label{c19}
\bC^{2N+1}=\bF_{H,\cR}(q)\oplus\bF_{P}(q),
\end{align}
where
\begin{align}\label{c20}
\bF_{P}(q)=\Gamma_{0,P}(q)\left(\bE_-(P_{0+}(q))\times\bE_+(P_{0-}(q))\right)
\end{align}
has dimension $N+1-k$ and $\bF_{H,\cR}(q)$ is an arbitrary
complementary subspace (necessarily of dimension $N+k$).   Denote by
$\pi_{H,\cR}(q)$ and $\pi_P(q)$ the projections associated to this
splitting.

\begin{rem}[More careful choice of $\bF_{H,\cR}(q)$]\label{c20y}
\textup{Henceforth, we'll work with a choice of
$\bF_{H,\cR}(q)\subset\bC^{2N+1}$ made as follows.  Using the
remarks made in Notation \ref{b13z} and the fact that
$\Gamma_{0,P}(q)$ is a real matrix, we first choose an $N+k$
dimensional subspace $\bF_{H,\cR}(q)\subset\bR^{2N+1}$ such that
\begin{align}\label{c20z}
\bR^{2N+1}=\bF_{H,\cR}(q)\oplus\bF_{P}(q).
\end{align}
We then take $\bF_{H,\cR}(q)$ in \eqref{c19} to be the
complexification of $\bF_{H,\cR}(q)\subset\bR^{2N+1}$.}

\end{rem}

For $\dov$ as in \eqref{c10} we can eliminate $u_P$ from the
boundary conditions \eqref{c14}.    That is to say,
$(u_H,u_P,\dos,\doh)$ satisfies \eqref{c14} if and only if
\begin{align}\label{c21}
\Gamma_{0,red}(q)(u_H,\dos,\doh)=0,\text{ and
}\;\;u_P=R_{P}(q)(u_H,\dos,\doh)
\end{align}
where
\begin{align}\label{c22}
\Gamma_{0,red}(q)(u_H,\dos,\doh):=\pi_{H,\cR}(q)\left(\Gamma_{0,H}(q)u_H+\Gamma_\cR(\dos,\doh)\right)
\end{align}
and
\begin{align}\label{c23}
R_P(q)(u_H,\dos,\doh)=-\Gamma_{0,P}^{-1}(q)\pi_P(q)\left(\Gamma_{0,H}(q)u_H+\Gamma_\cR(\dos,\doh)\right).
\end{align}

\begin{rem}\label{c24}
\textup{1.  Using Remark \ref{c20y} we see that $\Gamma_{0,red}(q)$
can be regarded as a map defined between complex spaces
\begin{align}\label{c24y}
\Gamma_{0,red}(q):\bC^{2N}\times\bC^d\to\bF_{H,\cR}(q)\subset\bC^{2N+1}
\end{align}
or as a map between real spaces
\begin{align}\label{c24z}
\Gamma_{0,red}(q):\bR^{2N}\times\bR^d\to\bF_{H,\cR}(q)\subset\bR^{2N+1}.
\end{align}
The choice should be clear from the context.  For example, in
statements like \eqref{c25a} below, we take $\Gamma_{0,red}(q)$ as
in \eqref{c24z}.}

2.  \textup{In view of Proposition \ref{c16t}, if a-transversality
holds for $W(z,q)$, then (a,p)-transversality means that the map
\begin{align}\label{c25}
\bC^{2N}\ni u_H\to\Gamma_{0,red}(q)(u_H,0,0)\in\bC^{2N+1}
\end{align}
has rank $N+k$.  The same holds for the map \eqref{c24z}.}

\end{rem}

The next step is to show that transversality of $W(z,q)$ implies
\begin{align}\label{c25a}
\ker\Gamma_{0,red}(q)=T_q\cC_\cB.
\end{align}
For this we need

\begin{prop}\label{c26}
For $q\in\cC_\cB$ consider the set $\bS(q)$ of solutions
$(\dow_+,\dow_-,\dos,\doh)$ of the fully linearized interior
equation on $\pm z\geq 0$, \eqref{c5}(a), for which $\dow_\pm$ are
bounded as $z\to\pm\infty$. We have
\begin{align}\label{c27}
\dim \;\bS(q)=2N+d+(N+1-k)
\end{align}
and
\begin{align}\label{c28}
\bS(q)=\left\{\begin{pmatrix}\phi_+'(z,q,a;\dop_+,\dos,\doh,\doa_+)\\\phi_-'(z,q,a;\dop_-,\dos,\doh,\doa_-)\\\dos\\\doh\end{pmatrix}:(\dop,\dos,\doh,\doa)\in\bR^{2N}\times\bR^d\times\bE_-(G_d(\up,0,0))\right\},
\end{align}
where, for example,
\begin{align}\label{c28a}
\begin{split}
&(a)\phi_+'(z,q,a;\dop_+,\dos,\doh,\doa_+):=\nabla_{p_+,s,h,a_+}\phi_+(z,p_+,s,h,a_+)(\dop_+,\dos,\doh,\doa_+)\text{
and }\\
&(b)\lim_{z\to+\infty}\phi_+'(z,q,a;\dop_+,\dos,\doh,\doa_+)=\dop_+,
\end{split}
\end{align}
and $a=a(p_\beta,s,h)$ is as in \eqref{b35}.

\end{prop}

\begin{proof}

Corollary \ref{bp1} shows that the map
\begin{align}\label{c29}
(\dop,\dos,\doh,\doa)\to\begin{pmatrix}\phi_+'(z,q,a;\dop_+,\dos,\doh,\doa_+)\\\phi_-'(z,q,a;\dop_-,\dos,\doh,\doa_-)\\\dos\\\doh\end{pmatrix}
\end{align}
is injective.   The functions $\phi_\pm(z,p_\pm,s,h,a_\pm)$ satisfy
the interior (nonlinear) profile equation in \eqref{c4}, and
differentiation with $\nabla_{p_\pm,s,h,a_\pm}$ shows that the
column vector on the right in \eqref{c29} gives a bounded solution
of \eqref{c5}(a). On the other hand any bounded solution $\dow$ of
\eqref{c5}(a) can be written $\dow=\dov+\cR(z,q,\dos,\doh)$, where
$\dov$ is in the $2N+(N+1-k)$ dimensional space of bounded solutions
to the problem with $(\dos,\doh)=0$.   This implies \eqref{c27}
(first) and then \eqref{c28}.

\end{proof}

As a corollary we obtain the following analogues of Propositions
\ref{c10a} and \ref{c16t}.

\begin{cor}\label{c300}
(1) The condition
\begin{align}\label{c35ab}
\mathrm{rank}\nabla_{a,p,s}\tilde{\Psi}(p,s,h,a)=2N+1
\end{align}
holds if and only if for all $g\in\bR^{2N+1}$ the problem
\begin{align}\label{c35ac}
\begin{split}
&\cL_0(z,q,\partial_z)\dow=\cL_{0,1}(z,q)(\dos,0)\text{ on }\pm z\geq 0\\
&\Gamma_2(\dow,\dow_z,\dos,0)=g\text{ on }z=0   \;\; (\Gamma_2\text{
as in }\eqref{c5})
\end{split}
\end{align}
has a solution $(\dow,\dos)$ with $\dow$ bounded.

(2)  The condition \eqref{c35ab} means that the rank of the matrix
\begin{align}\label{c35ad}
(\Gamma_{0,H}(q),\Gamma_{0,P}(q),\Gamma_{\cR}):\bC^{2N}\times\left(\bE_-(P_{0+}(q))\times\bE_+(P_{0-}(q))\right)\times
\bC^d\to\bC^{2N+1}
\end{align}
is $2N+1$ when $(\dos,\doh)=(\dos,0)$.

(3)  If a-transversality holds for $W(z,q)$, then
(a,p,s)-transversality means that the map
\begin{align}\label{c35adz}
\bC^N\times\bC\ni(u_H,s)\to\Gamma_{0,red}(q)(u_H,s,0)\in\bC^{2N+1}
\end{align}
has rank $N+k$.
\end{cor}

\begin{proof}
Let $\phi'=\phi'(z,p,s,h,a;\dop,\dos,\doh,\doa)$ be given by the
first two components of the column vector in \eqref{c28}.  Part (1)
then follows directly from Prop. \ref{c26} and the observation that
\begin{align}\label{c301}
\Gamma_2(\phi',\phi'_z,\dos,\doh)=\nabla_{(p,s,h,a)}\tilde{\Psi}(p,s,h,a)(\dop,\dos,\doh,\doa).
\end{align}
Let $\Gamma_3(q)$ be as in \eqref{c9} and note that for
\begin{align}
\dov=\dow-\cR(z,q,\dos,\doh)
\end{align}
as in \eqref{c8}, we have
\begin{align}\label{c303}
\Gamma_2(\dow,\dow_z,\dos,\doh)=g\Leftrightarrow\Gamma_3(q)(\dov,\dov_z,\dos,\doh)=g.
\end{align}
After rewriting $\Gamma_3(q)$ as in \eqref{c14}, we deduce part (2)
from part (1).   Part (3) then follows directly from part (2) and
the definition of $\Gamma_{0,red}(q)$.

\end{proof}

\begin{prop}\label{c30}
\textbf{} (1)  Assume $W(z,q)$ is strongly transversal.   Then we
have
\begin{align}\label{c31}
\begin{split}
&(a)\ker\Gamma_{0,red}(q)=T_q\cC_\cB\\
&(b)T_q\cC_\cB=\\
&\qquad\{(\dop_+,\dop_-,\dos,\doh):\text{ there exists a solution
}(\dow_+,\dow_-,\dos,\doh)\text{ of \eqref{c5} with
}\lim_{z\to\pm\infty}\dow_\pm=\dop_\pm\}.
\end{split}
\end{align}

(2)  The same conclusions hold if we assume just that $W(z,q)$ is
transversal.

\end{prop}

\begin{proof}
\textbf{1. } We show that both sets in \eqref{c31}(a) are equal to
the set on the right in \eqref{c31}(b).   Assume that $W(z,q)$ is
strongly transversal; the proof in the other case is essentially the
same.

\textbf{2. }Suppose $(u_H,\dos,\doh)\in\ker\Gamma_{0,red}(q)$.  By
\eqref{c21}-\eqref{c23} there exists
\begin{align}\label{c32r}
u_P\in\bE_-(P_{0+}(q))\times\bE_+(P_{0-}(q))
\end{align}
such that $\dov_\pm$ as in \eqref{c10} satisfies \eqref{c9} with
\begin{align}\label{c32t}
\lim_{z\to\pm\infty}\dov_\pm=u_{H\pm}.
\end{align}
But then
\begin{align}\label{c32s}
\dow=\dov+\cR(z,q,\dos,\doh)
\end{align}
satisfies \eqref{c5} with
\begin{align}\label{c32}
\lim_{z\to\pm\infty}\dow_\pm=u_{H\pm}.
\end{align}

Conversely, if  $(u_H,\dos,\doh)$ is such that there exists $\dow$
satisfying \eqref{c5} and \eqref{c32}, then $\dov$ as in
\eqref{c32s} satisfies \eqref{c9} and \eqref{c32t}.  So there
exists $u_P$ as in \eqref{c32r} such that $\dov$ is given by
\eqref{c10} and \eqref{c14} holds.    Apply $\pi_{H,\cR}(q)$ to
\eqref{c14} to obtain $(u_H,\dos,\doh)\in\ker\Gamma_{0,red}(q)$.

\textbf{3. }Recall from \eqref{b34e} that $\cC_\cB$ is given as the
graph of
\begin{align}\label{c34}
p_\alpha=p_\alpha(p_\beta,s,h),
\end{align}
so its tangent space at $(p,s,h)$ is given by the graph of
\begin{align}\label{c35}
\dop_\alpha=\nabla_{(p_\beta,s,h)}p_\alpha(p_\beta,s,h)(\dop_\beta,\dos,\doh).
\end{align}

For $\tilde{\Psi}(p,s,h,a)$ as in \eqref{b28},  the equation
$\tilde{\Psi}=0$ also yielded the relation
\begin{align}\label{c36}
a=a(p_\beta,s,h),
\end{align}
whose linearization at $(p_\beta,s,h)$ is
\begin{align}\label{c38}
\doa=\nabla_{(p_\beta,s,h)}a(p_\beta,s,h)(\dop_\beta,\dos,\doh).
\end{align}
We can also differentiate the equation $\tilde{\Psi}(p,s,h,a)=0$
in all variables to obtain the linearized form
\begin{align}\label{c33}
\nabla_{(p,s,h,a)}\tilde{\Psi}(p,s,h,a)(\dop,\dos,\doh,\doa)=0,
\end{align}
Observe that since the rank conditions hold
\begin{align}\label{c34z}
(\dop,\dos,\doh,\doa)\text{ satisfies
\eqref{c33}}\Leftrightarrow\text{ both \eqref{c35} and \eqref{c38}
hold.}
\end{align}

Suppose now that $(\dop,\dos,\doh)\in T_q\cC_\cB$.  Then \eqref{c35}
holds.  If we take $\doa$ as in \eqref{c38}, then \eqref{c33} holds
and implies that for these choices of $(\dop,\dos,\doh,\doa)$ the
column vector in \eqref{c28} satisfies \eqref{c5}.   From
\eqref{c28a}(b) we see that $(\dop,\dos,\doh)$ is an element of the
set on the right in \eqref{c31}(b).

Conversely, suppose $(\dop,\dos,\doh)$ is such that there exists a
solution $(\dow_+,\dow_-,\dos,\doh)$ as on the right in
\eqref{c31}(b).   By Prop. \ref{c26} this solution must have the
form of the column vector in \eqref{c28} for some $\doa$.  Since
this solution satisfies the boundary conditions in \eqref{c5},
this means that \eqref{c33} holds, which implies \eqref{c35}.

\end{proof}

%\begin{rem}\label{c35aa}

%1.  It is easy to check that Prop. \ref{c30} still holds if we
%replace the assumption of (a,p)-transversality by the weaker
%condition
%\begin{align}\label{c35ab}
%\mathrm{rank}\nabla_{a,p,s}\tilde{\Psi} \text{ has full rank}.
%\end{align}

%2.  Arguments similar to those in the above proof yield an analogue
%of Prop. \ref{c10a} for the condition \eqref{c35ab}:

%The condition \eqref{c35ab} holds if and only if for all
%$g\in\bR^{2N+1}$ the problem
%\begin{align}\label{c35ac}
%\begin{split}
%&\cL_0(z,q,\partial_z)\dov=\cL_{0,1}(z,q)(\dos,0)\text{ on }\pm z\geq 0\\
%&\Gamma_2(\dov,\dov_z,\dos,0)=g\text{ on }z=0   \;\;
%(\Gamma_2\text{ as in }\eqref{c5})
%\end{split}
%\end{align}
%has a solution $(\dov,\dos)$ with $\dov$ bounded.

%3.  Parallel to Proposition \ref{c16t}, the condition \eqref{c35ab}
%means that the rank of the matrix
%\begin{align}\label{c35ad}
%(\Gamma_{0,H}(q),\Gamma_{0,P}(q),\Gamma_{\cR}):\bC^{2N}\times\left(\bE_-(P_{0+}(q))\times\bE_+(P_{0-}(q))\right)\times
%\bC^d\to\bC^{2N+1}
%\end{align}
%is $2N+1$ when $(\dos,\doh)=(\dos,0)$.
%\end{rem}

\begin{cor}\label{c35z}
Assume $\uw(z)$ is either transversal or strongly transversal. For
$q\in\cC_\cB$ near $\uq$ the linearized hyperbolic transmission
problem at $q$ (\eqref{c2} with $f=0$, $g=0$) can now be written
\begin{align}\label{c36z}
\begin{split}
&\sum^{d-1}_{j=0}A_j(p)\partial_j\dou+\cA_d(p,s,h)\partial_d\dou=0\\
&\Gamma_{0,red}(q)(\dou,d\dopsi)=0\text{ on }x_d=0.
\end{split}
\end{align}
\end{cor}

\section{Stability determinants}\label{stability}

\subsection{Lopatinski determinants}\label{Lopatinski}

\textbf{\quad}In this section we define and study the relationship
between two stability determinants, the Lopatinski determinant,
$D_{Lop}(q,\hzeta)$, and the modified Lopatinski determinant,
$D_{Lop,m}(q,\hzeta)$, for the linearized hyperbolic problem
\eqref{c2}.  Later we'll see how these arise in low frequency
expansions of, respectively, the standard and modified Evans
functions.   But for now we study them in the context of inviscid
$\cC$ shocks, where $\cC$ is a shock manifold as in Assumption
\ref{i3}

%In this section we define and study the relationship between two
%stability determinants, $D_{Lop}(q,\hzeta)$ and
%$D_{Lop,\chi}(q,\hzeta)$, that govern the solution of the inviscid
%hyperbolic shock problem \eqref{c1}. $D_{Lop,\chi}(q,\hzeta)$ is
%more directly connected to the construction of Kreiss symmetrizers
%for the inviscid shock problem, while $D_{Lop}(q,\hzeta)$  is more
%readily related to the low frequency limit of the modified Evans
%function defined later (see \eqref{g53}, for example).

\subsubsection{$\cC$ shocks}
 Dropping dots and some hats, we write the Laplace-Fourier
transform of \eqref{c2}, with $f=0$, $g=0$, as
\begin{align}\label{c37z}
\begin{split}
&(a)\;A_0(p)(i\htau+\hgamma)u+\sum^{d-1}_{j=1}A_j(p)i\heta_j
u+\cA_d(p,s,h)\partial_du=0\text{ in }\pm x_d\geq 0\\
%&\Gamma_{0,red}(q)(u,(i\htau+\hgamma)\psi,i\heta\psi)=0.
&(b)\;\hGamma_{\chi}(q,\hzeta)(u_+,u_-,\psi)=0,
\end{split}
\end{align}
where $\hGamma_\chi(q,\zeta):\bC^{2N+1}\to\bC^{N+k}$ is defined by
\begin{align}\label{cc40e}
\hGamma_\chi(q,\hzeta)(u_+,u_-,\psi)=\chi'_{p_+}(q)u_++\chi'_{p_-}(q)u_-+\chi'_s(q)(i\htau+\hgamma)\psi+\chi'_h(q)i\heta\psi.
\end{align}
\begin{defn}\label{cc38z}

1.  For $\psi\in\bC$, $u_\pm\in\bC^{N}$,  and $\hzeta\in
S^d_+=S^d\cap\{\hgamma>0\}$ we set
\begin{align}\label{c39}
\begin{split}
%&\hGamma_{0,red}(q,\hzeta)(u_+,u_-,\psi):=\Gamma_{0,red}(q)(u,(i\htau+\hgamma)\psi,i\heta\psi)\\
&H_{0\pm}(q,\hzeta)=-\cA_d(p_\pm,s,h)^{-1}\left(A_0(p_\pm)(i\htau+\hgamma)+\sum^{d-1}_{j=1}A_j(p_\pm)i\heta_j
\right)\\
&\bE_-(H_0(q,\hzeta)):=\bE_-(H_{0+}(q,\hzeta))\times
\bE_+(H_{0-}(q,\hzeta)),
\end{split}
\end{align}
where, as usual, $\bE_-(H_{0+}(q,\hzeta))$ denotes the generalized
eigenspace of $H_{0+}(q,\zeta)$ associated to eigenvalues $\mu$ with
$\Re\mu<0$.

2.  The \emph{Lopatinski determinant} for the problem  \eqref{c37z}
is the $ (N+k)\times (N+k)$ determinant
\begin{align}\label{c40c}
D_{Lop}(q,\hzeta):=\det\left(\chi'_p(q)\bE_-(H_0(q,\hzeta)),\chi'_s(q)(i\htau+\hgamma)+\chi'_h(q)i\heta\right).
\end{align}
The first $N+k-1$ columns of the matrix in \eqref{c40c} are computed
using an orthonormal basis of $\bE_-(H_0(q,\hzeta))$.

We'll say that the \emph{uniform Lopatinski condition} holds at
$q\in\cC$  when
\begin{align}\label{c40b}
\exists c>0\text{ such that }|D_{Lop}(q,\hzeta)|\geq c \text{ for
all }\hzeta\in S^d_+.
\end{align}
\end{defn}

\begin{defn}\label{c400}

1.  When $E$ and $F$ are subspaces of  $\bC^D$ with $\dim E +\dim
F=D$, $\det(E,F)$ denotes the determinant formed by taking
orthonormal bases in $E$ and $F$.  Up to a sign this determinant is
independent of the choice of bases.  Determinants of nonsquare
matrices are defined to be $0$.

2.  We define the \emph{modified Lopatinski determinant} for the
hyperbolic problem \eqref{c37z}
\begin{align}\label{c40}
D_{Lop,m}(q,\hzeta)=\det\left(\bE_-(H_0(q,\hzeta))\times\bC,
\ker\hGamma_\chi(q,\hzeta)\right).
\end{align}
The \emph{modified uniform Lopatinski condition} holds at $q\in\cC$
when there is a $c>0$ such that $|D_{Lop,m}(q,\hzeta)|\geq c$ for
all $\hzeta\in S^d_+$.

\end{defn}

\begin{rem}\label{c40z}
1.  \textup{Hyperbolicity (H1) implies that for all $\hzeta\in
S^d_+$
\begin{align}\label{c40y}
\dim\bE_-(H_{0+}(q,\hzeta))=N-R_-,\;\;\dim\bE_+(H_{0-}(q,\hzeta))=N-L_+
\end{align} for $R_-$, $L_+$ as in Defn. \ref{b19}.   Thus, the
determinant in \eqref{c40} is a $(2N+1)\times(2N+1)$ determinant
when $\dim\ker\hGamma_{\chi}(q,\hzeta)=N+1-k$.}

2.   \textup{The spaces $\bE_\mp(H_{0\pm}(q,\hzeta))$ define
$C^\infty$ vector bundles over $S^d_+$.}

3. \textup{The uniform Lopatinski condition allows one to construct
Kreiss symmetrizers and prove maximal estimates for the linearized
inviscid problem \eqref{c2}; see section \ref{existence}.}

4.  \textup{Both determinants have been defined only for $\hzeta\in
S^d_+$.
 The function $D_{Lop,m}(q,\cdot)$ has a continuous extension to any
subset of $\oS^d_+$ where $\bE_-(H_0(q,\hzeta))$ is continuous and
$\hGamma_{\chi}(q,\hzeta)$ maintains full rank, while
$D_{Lop}(q,\cdot)$ has a continuous extension to any subset of
$\oS^d_+$ where $\bE_-(H_0(q,\hzeta))$ is continuous.}

\end{rem}

Suppose that $\chi_1$ and $\chi_2$ are local defining functions
(Defn. \ref{i8y}) for $\cC$ near $q$.  Since
\begin{align}\label{c402}
\ker\chi_1'(q)=T_q\cC=\ker\chi'_2(q)\Rightarrow\ker\hGamma_{\chi_1}(q,\hzeta)=\ker\hGamma_{\chi_2}(q,\hzeta),
\end{align}
we see that $D_{Lop,m}(q,\hzeta)$, and hence the validity of the
modified uniform Lopatinski condition, is independent of the choice
of local defining function.   The determinant  $D_{Lop}(q,\hzeta)$
is clearly not independent of the choice of local defining function,
but we will show that the uniform Lopatinski condition is.  In the
proof we use the following general fact about determinants.

\begin{prop}\label{ca1}
Let $M$ be a $p\times q$ matrix, $p\neq q$,  $E$ a  $q\times p$
matrix, and denote by $\ker M$ a matrix whose columns form an
orthonormal basis of the kernel of $M$.  Then we have
\begin{align}\label{ca2}
\det ME=c(M)\det\begin{pmatrix}E&\ker M\end{pmatrix},
\end{align}
where $c(M)\neq 0$ if $M$ has full rank and is otherwise defined
to be $0$.

 If $p>q$ or if $p\leq q$ and $M$ is not full rank, the
matrix $\begin{pmatrix}E&\ker M\end{pmatrix}$ on the right is not
square and its determinant is then defined to be zero.

\end{prop}

\begin{proof}
\textbf{1. }When $p>q$ or when $p\leq q$ and $M$ is not full rank,
both determinants are zero.

\textbf{2. }Assume then that $p\leq q$ and $M$ has full rank $p$,
so the determinant on the right in \eqref{ca2} is $q\times q$.  In
this case $MM^*$ is invertible and
\begin{align}\label{ca3}
\begin{pmatrix}M\\(\ker
M)^*\end{pmatrix}^{-1}=\begin{pmatrix}M^*(MM^*)^{-1}&\ker
M\end{pmatrix}.
\end{align}
Note that $(\ker M)^*$ is a $(q-p)\times q$ matrix.

\textbf{3. }We have
\begin{align}\label{ca4}
\begin{pmatrix}M\\(\ker
M)^*\end{pmatrix}\begin{pmatrix}E&\ker
M\end{pmatrix}=\begin{pmatrix}ME&0\\(\ker
M)^*E&I_{(q-p)\times(q-p)}\end{pmatrix},
\end{align}
so \eqref{ca2} holds with
\begin{align}\label{ca5}
c(M):=\det\begin{pmatrix}M\\(\ker M)^*\end{pmatrix},
\end{align}
which is  nonzero when $M$ has full rank by part 2.
\end{proof}

\begin{prop}\label{c403}
(a)  The uniform Lopatinski condition is independent of the choice
of local defining function used to compute $D_{Lop}(q,\hzeta)$.

(b)The uniform Lopatinski condition at $q$ implies the modified
uniform Lopatinski condition at $q$.

\end{prop}

\begin{proof}
\textbf{1. }Let $\chi_1$ and $\chi_2$ be local defining functions
for $\cC$ near $q$, and let $D_{Lop,\chi_i}$ denote the
corresponding determinants \eqref{c40c}.  Suppose there exists $c>0$
such that
\begin{align}\label{c406}
|D_{Lop,\chi_1}(q,\hzeta)|\geq c \text{ for }\hzeta\in S^d_+.
\end{align}
Then $\hGamma_{\chi_1}(q,\hzeta)$ (and thus also
$\hGamma_{\chi_2}(q,\hzeta)$ by \eqref{c402}) has full rank $N+k$
for $\hzeta\in S^d_+$. Applying Prop. \ref{ca1} with
\begin{align}
M:=\hGamma_{\chi_i}(q,\hzeta)\text{ and
}E:=\bE_-(H_0(q,\hzeta))\times\bC,
\end{align}
we obtain
\begin{align}\label{c404}
D_{Lop,\chi_i}(q,\hzeta)=c_i(q,\hzeta)D_{Lop,m}(q,\hzeta), \;i=1,2
\end{align}
where $c_i(q,\hzeta)\neq 0$ and is continuous near $(q,\hzeta)$.
(More precisely, $E$ is a $(2N+1)\times (N+k)$ matrix whose columns
form a basis of the given space.)

\textbf{2. }To complete the proof of part (a) it suffices to show
that $\hGamma_{\chi_i}(q,\hzeta)$ cannot drop rank at a point
$\hzeta_0$ with $\hgamma_0=0$, for then $c_i(q,\hzeta)$ must be
continuous and nonvanishing on $\oS^d_+$.  When
$\bE_-(H_0(q,\hzeta))$ extends continuously to $\oS^d_+$, so does
$D_{Lop,\chi_1}(q,\cdot)$, and \eqref{c406} implies that
$\hGamma_{\chi_1}(q,\cdot)$ has full rank on $\oS^d_+$.

In the general case the Lopatinski determinants are not defined for
$\hgamma=0$, but we can substitute compactness of Grassmannians for
compactness of $\oS^d_+$.      If $\hGamma_{\chi_1}(q,\hzeta)$ drops
rank at $\hzeta_0$, then
\begin{align}\label{cb5}
\det(\hGamma_{\chi_1}(q,\hzeta_0)E_0)=0
\end{align}
for any $E_0$ in the set of limit points of
$\bE_-(H_0(q,\hzeta_j))\times\bC$ as $\hzeta_j\to\hzeta_0$ with
$\hgamma_j>0$.   But then we must have
\begin{align}\notag
D_{Lop,\chi_1}(q,\hzeta_j)\to 0
\end{align}
for some sequence $\hzeta_j$ with $\hgamma_j>0$, a contradiction.

\textbf{3. }Part (b) follows directly from \eqref{c404} and the fact
that $c_1(q,\hzeta)$ is continuous and nonvanishing on $\oS^d_+$.
\end{proof}

\begin{rem}\label{c407}

\textup{(a)\;A similar argument fails for the reverse implication in
part (b) above, because although
$\hGamma_\chi(q,\hzeta_j)\to\hGamma_\chi(q,\hzeta_0)$, we don't have
\begin{align}\notag
\ker\hGamma_\chi(q,\hzeta_j)\to\ker\hGamma_\chi(q,\hzeta_0)
\end{align}
when $\hGamma_\chi(q,\cdot)$ drops rank at $\hzeta_0$.}

\textup{(b)\;Observe that if $\chi$ is a local defining for $\cC$
near $q$ and if $\chi'_{p,s}(q)$ (or, alternatively, $\chi'_p(q)$)
has full rank $N+k$, then the same is true for any other defining
function near $q$. This is relevant to the next Proposition.}

\end{rem}

\begin{prop}\label{c408}
Let $\chi$ be a local defining function for $\cC$ near $q$.

(a) Suppose $\chi'_p(q)$ has full rank $N+k$.  Then
\begin{align}\label{c4088}
D_{Lop}(q,\hzeta)=c(q,\hzeta)D_{Lop,m}(q,\hzeta), \;i=1,2
\end{align}
with $c(q,\hzeta)$ continuous and nonvanishing on $\oS^d_+$.
Consequently, the modified uniform Lopatinski condition at $q$
implies the uniform Lopatinski condition at $q$.

(b)Suppose $\chi'_{p,s}(q)$ has full rank $N+k$ and $d=1$.  Then
\eqref{c4088} holds with $c(q,\hzeta)$ continuous and nonvanishing
on $\oS^d_+$.    Again, the modified uniform Lopatinski condition at
$q$ implies the uniform Lopatinski condition at $q$.

\end{prop}

\begin{proof}
\textbf{1. }If $\chi'_p(q)$ has full rank, then
$\hGamma_\chi(q,\hzeta)$ maintains full rank even when restricted to
the subspace $\psi=0$; hence it maintains full rank for $\hzeta\in
\oS^d_+$.  Thus, \eqref{c4088} holds  with $c(q,\hzeta)$ continuous
and nonvanishing on $\oS^d_+$.

\textbf{2. }Since $\chi'_{p,s}(q)$ has full rank, when $d=1$ it
follows that $\hGamma_\chi(q,\htau,\hgamma)$ must have full rank for
$\hzeta\in\oS^d_+$.  The result then follows as in part (a).

\end{proof}

\begin{prop}\label{c410}
Suppose the uniform Lopatinski condition holds at $q\in\cC$.

(a)  Then $\chi'_{p,s}(q)$ has full rank $N+k$ for any local
defining function $\chi$.

(b)   If $d\geq 2$, then $\chi'_p(q)$ has full rank $N+k$.

\end{prop}

\begin{proof}
\textbf{1.  (a)} This follows immediately from \eqref{c40b} by
taking $\heta=0$.

\textbf{2. (b)} If  $\chi'_p(q)$ does not have full rank,   we claim
\begin{align}\label{c412}
\dim\ker\hGamma_\chi(q,\hzeta)\geq N+2-k\text{ for some
}\hzeta\in\oS^d_+.
\end{align}
%and so the uniform Lopatinski condition must fail.
By part (a) we can choose coordinates
$(p_\alpha,p_\beta)\in\bR^{N-1+k}\times\bR^{N+1-k}$ such that
$\chi'_{p_\alpha,s}(q)$ is nonsingular; write
$u=(u_\alpha,u_\beta)\in\bC^{2N}$.  If  $\chi'_p(q)$ does not have
full rank, it must have rank $N-1+k$, so the map
\begin{align}\label{c40l}
u\to\chi'_p(q)u=\chi'_{p_\alpha}(q)u_\alpha+\chi'_{p_\beta}(q)u_\beta
\end{align}
has an $N+1-k$ dimensional kernel.  This yields an $N+1-k$
dimensional subspace $S(q)\subset\ker\hGamma_\chi(q,\hzeta)$
independent of $\hzeta$, namely,
\begin{align}\label{c40m}
\begin{split}
&S(q)=\{(u,\psi)=(u,0):u\in\ker\chi'_p(q)\}=\\
&\qquad\{(u_\alpha,u_\beta,0):(u_\alpha,0)=-(\chi'_{p_\alpha,s}(q))^{-1}\chi'_{p_\beta}(q)u_\beta,
\;u_\beta\in\bC^{N+1-k}\}.
\end{split}
\end{align}

\textbf{3. } To finish we show that for certain $\hzeta$, there are
elements of $\ker\hGamma_\chi(q,\hzeta)$ that have nonzero $\psi$
components.

Consider the $d$ vectors in $\bR^N$:
\begin{align}\label{c40n}
\chi'_s(q),\chi'_{h_i}(q), i=1,\dots,d-1.
\end{align}
If these vectors are linearly dependent over $\bR$, there clearly
exist $\htau$, $\heta$ such that for $\hzeta=(\htau,0,\heta)$,
\begin{align}\label{c40o}
(0,0,1)\in\ker\hGamma_\chi(q,\hzeta).
\end{align}
On the other hand if the vectors \eqref{c40n} are linearly
independent over $\bR$, and thus over $\bC$, the range of the map
\begin{align}\label{c40p}
(\hzeta,\psi)\to\chi'_s(q)(i\htau+\hgamma)\psi+\chi'_h(q)i\heta\psi
\end{align}
must have nontrivial intersection in $\bC^{N+k}$ with the $N-1+k$
dimensional range of the map
\begin{align}\label{c40q}
(u_\alpha,u_\beta)\to\chi'_{p_\alpha}(q)u_\alpha+\chi'_{p_\beta}(q)u_\beta.
\end{align}
This is because $d\geq 2$.  In view of \eqref{cc40e} this means that
in any case, for some $\hzeta$, there are elements of
$\ker\hGamma_\chi(q,\hzeta)$ that have nonzero $\psi$ components.
Since $S(q)\subset\ker\hGamma_\chi(q,\hzeta)$ as well, this implies
\eqref{c412}.

\textbf{4. }If the set of bad directions $\hzeta$ where \eqref{c412}
holds includes a point $\hzeta_0\in S^d_+$,  then
$D_{Lop,\chi}(q,\hzeta_0)=0$ since $\hGamma_\chi(q,\hzeta_0)$ does
not have full rank. So suppose all bad directions $\hzeta_0$ satisfy
$\hgamma_0=0$. Use compactness of Grassmannians as before to obtain
an $N+k$ dimensional limit point $E_0$ of
$\bE_-(H_0(q,\hzeta_j))\times\bC$ as $\hzeta_j\to\hzeta_0$ with
$\hgamma_j>0$.    We then obtain a sequence
\begin{align}\notag
D_{Lop,\chi}(q,\hzeta_j)\to 0
\end{align}
just as in the proof of Proposition \ref{c403}.
\end{proof}

\subsubsection{$\cC_\cB$ shocks}

\textbf{\quad} Suppose now that $\cC=\cC_\cB$ and that $W(z,q)$ is a
viscous profile corresponding to the planar shock $q\in\cC_\cB$.  We
have seen that for any local defining function $\chi$,  if $W(z,q)$
is transversal, then $\chi_{p,s}'(q)$ has full rank $N+k$, and if
$W(z,q)$ is strongly transversal, then $\chi'_p(q)$ has full rank.
Thus, we obtain the following immediate corollaries of Propositions
\ref{c403}, \ref{c408}, and \ref{c410}.

%Next we obtain an easy corollary which does not address uniform
%behavior as $\hgamma\to 0$.

%\begin{cor}\label{cb1} Fix $(q,\hzeta)$ with $\hgamma>0$ and assume
%$W(z,q)$ is transversal.  Then $\hGamma_\chi(q,\hzeta)$ as in
%\eqref{c40e} has full rank $N$ and
%\begin{align}\label{cb2}
%D_{Lop,\chi}(q,\hzeta)=c(q,\hzeta)D_{Lop}(q,\hzeta),
%\end{align}
%where $c(q,\hzeta)\neq 0$ and is continuous near $(q,\hzeta)$.

%(b)  If $\mathrm{rank}\;\hGamma_\chi(q,\hzeta)<N$ both
%determinants are zero.
%\end{cor}

%\begin{proof}
%\textbf{1. }a-transversality allows $\hGamma_{0,red}(q,\hzeta)$ to
%be defined, and since Proposition \ref{c30} remains true when
%$\chi'_{p,s}(q)$ has rank $N+k$, the equality in \eqref{c40f} still
%holds.

%\textbf{2. }  Since $\chi'_{p,s}(q)$ has rank $N$, so does
%$\hGamma_\chi(q,\hzeta)$ when $\hgamma>0$, so the assertion follows
%from Prop. \ref{ca1} if we take
%\begin{align}
%M:=\hGamma_\chi(q,\hzeta)\text{ and
%}E:=\bE_-(H_0(q,\hzeta))\times\bC.
%\end{align}
%More precisely, $E$ is a $(2N+1)\times (N+k)$ matrix whose columns
%form a basis of the given space.

%\end{proof}

\begin{cor}\label{cb3}
Assume  $W(z,q)$ is transversal.

(a)  When $d=1$,  the uniform Lopatinski condition holds at $q$ if
and only if the modified uniform Lopatinski condition holds at $q$.

(b) For $d\geq 1$ the uniform Lopatinski condition at $q$ implies
the modified uniform Lopatinski condition at $q$.

 (c) Suppose $W(z,q)$ is strongly transversal. Then the modified uniform
Lopatinski condition at $q$ implies the uniform Lopatinski condition
at $q$.

\end{cor}

\begin{cor}\label{c40i}

When $d\geq 2$ and $W(z,q)$ is transversal,  if the uniform
Lopatinski condition holds at $q$, then $W(z,q)$ is strongly
transversal.

\end{cor}

These results  are useful in the analysis  of  the viscous stability
determinants or Evans functions defined in the following sections;
see Proposition \ref{j3}.

%The situation is a little different when $d=1$:

%\begin{prop}\label{c40s}
%Let $d=1$ and assume $W(z,q)$ is transversal. Then the uniform
%Lopatinski$_\chi$ condition \eqref{c40b} is equivalent to the
%uniform Lopatinski condition, Definition \ref{c38z}, and
%\begin{align}\label{c40t}
%D_{Lop}(q,\htau,\hgamma)=c(q)D_{Lop,\chi}(q,\htau,\hgamma)
%\end{align}
%for some smooth nonvanishing function $c(q)$.

%\end{prop}

%\begin{proof}

%Choose coordinates $(p_\alpha,p_\beta)\in\bR^{N-1}\times\bR^{N+1}$
%such that $\chi'_{p_\alpha,s}(q)$ is nonsingular, and write
%$u=(u_\alpha,u_\beta)\in\bC^{2N}$. For
%$\hat{\lambda}=i\htau+\hgamma\in\oS^1_+$ define
%\begin{align}\label{c40u}
%T(q,\hlambda)(u_\alpha,\psi)=\chi'_{p_\alpha}(q)u_\alpha+\chi'_s(q)\hlambda\psi.
%\end{align}
%The equality \eqref{c40f} holds and, writing elements of
%$\ker\hGamma_\chi(q,\htau,\hgamma)$ as $(u_\alpha,\psi,u_\beta)$,
%we can form a basis for that kernel of the form:
%\begin{align}\label{c40v}
%\{\left(T(q,\hlambda)^{-1}(\chi'_{p_\beta}(q)e_\beta),e_\beta\right):e_\beta\in\cB_{N+1}\},
%\end{align}
%where $\cB_{N+1}$ is the standard basis of $\bC^{N+1}$.   After
%substituting this parametrization of
%$\ker\hGamma_{0,red}(q,\hzeta)$ into the determinant \eqref{c40},
%one obtains  a nonvanishing multiple of
%$D_{Lop,\chi}(q,\htau,\hgamma)$ after a few row and column
%operations.

%\end{proof}

\begin{rem}\label{c41}
\textup{(a) For $\psi\in\bC$, $u_\pm\in \bC^N$, and
$\hzeta\in\oS^d_+$, we define
\begin{align}\label{c420}
\hGamma_{0,red}(q,\hzeta)(u_+,u_-,\psi):=\Gamma_{0,red}(q)(u,(i\htau+\hgamma)\psi,i\heta\psi).
\end{align}}
\textup{If $\chi$ is a local defining function for $\cC_\cB$ near
$q$, then Proposition \ref{c30} implies
\begin{align}\label{c421}
\ker\hGamma_{0,red}(q,\hzeta)=\ker\hGamma_\chi(q,\hzeta)\text{ for
all }\hzeta\in\oS^d_+.
\end{align}
Thus, we can rewrite $D_{Lop,m}$ as
\begin{align}\label{c422}
D_{Lop,m}(q,\hzeta)=\det\left(\bE_-(H_0(q,\hzeta))\times\bC,
\ker\hGamma_{0,red}(q,\hzeta)\right).
\end{align}
This form of $D_{Lop,m}$ appears naturally in the low frequency
expansion of the modified Evans function (see, e.g., Cor.
\ref{f23}).}

\textup{(b) Since
\begin{align}\label{c42}
\hGamma_{0,red}(q,\hzeta)(u,\psi)=\pi_{H,\cR}(q)\left(\Gamma_{0,H}(q)u+\begin{pmatrix}[\cR(0,q,i\htau+\hgamma,i\heta)]\psi\\
[\cR_z(0,q,i\htau+\hgamma,i\heta)]\psi\\0\end{pmatrix}\right),
\end{align}
by considering $(u_+,u_-,\psi)=(0,0,1)$, we see that the  uniform
Lopatinski condition implies
\begin{align}\label{c43}
\begin{pmatrix} [\cR(0,q,i\htau+\hgamma,i\heta)]\\ [\cR_z(0,q,i\htau+\hgamma,i\heta)]\end{pmatrix}\neq
0\text{ for all }\hzeta\in\overline{S}^d_+.
\end{align}
This is immediate for $\hzeta\in S^d_+$.  It is true for
$\hzeta\in\oS^d_+$ as well, since $\hGamma_{0,red}(q,\hzeta)$ cannot
drop rank at points with $\hgamma=0$ by \eqref{c421} and the
argument used in the proof of Prop. \ref{c403}.}
\end{rem}

\subsection{The standard Evans function}\label{standard}
\textbf{}

In this section and the next we consider $q\in\cC_\cB$ and define
the Evans functions that turn out to govern such nonlinear stability
questions as the small viscosity limit of curved viscous shocks
(i.e., convergence of viscous to inviscid shocks as viscosity goes
to zero) and the long time stability of planar viscous shocks.

 For the partially linearized operator $\cL$ as in
\eqref{c5a} consider
\begin{align}\label{ce0}
\begin{split}
&\cL(z,q,\zeta,\partial_z)u=f \text{ on }\pm z\geq 0\\
&[u]=0, [u_z]=0\text{ on }z=0.
\end{split}
\end{align}
With $U:=(u,u_z)$ we can as in \eqref{c6d} rewrite \eqref{ce0} when
$f=0$ as a $2N\times 2N$ first-order problem:
\begin{align}\label{ce1}
\begin{split}
&(a)\;\partial_zU=G(z,q,\zeta)U\text{ on }\pm z\geq 0,\\\;\;
&(b)\;\Gamma_s(U_+,U_-)=0\text{ on }z=0,
\end{split}
\end{align}
where $\Gamma_s:\bC^{2N}\times\bC^{2N}\to\bC^{2N}$ is given by
\begin{align}\label{ce2}
\Gamma_s(U_+,U_-):=\begin{pmatrix} [u(0)]\\
[\partial_zu(0)]\end{pmatrix}.
\end{align}

For $\hzeta\in\oS^d_+$, $ \rho>0$, let $E_\pm(q,\hzeta,\rho)$ denote
the set of initial data, $U_\pm(0)$, of bounded solutions of
\eqref{ce1}(a) on $\pm z\geq 0$. Below we show that the spaces
$E_\pm(q,\hzeta,\rho)$ are $C^\infty$ in $(q,\hzeta,\rho)$ and
satisfy
\begin{align}\label{cc3}
\dim E_\pm(q,\hzeta,\rho)=N\text{ for }\hzeta\in\oS^d_+,\; \rho>0.
\end{align}

\begin{defn}\label{cc1}
1.  The \emph{standard Evans function} is defined for $\hzeta\in
\oS^d_+$ and $\rho>0$ by the $4N\times 4N$ determinant
\begin{align}\label{cc2}
D_s(q,\hzeta,\rho)=\det\left(E_-(q,\hzeta,\rho)\times
E_+(q,\hzeta,\rho),\ker\Gamma_s\right).
\end{align}

2.  The profile $W(z,q)$ satisfies the  \emph{standard low frequency
Evans condition at $q$} if and only if for some positive constants
$c$ and $\rho_0$:
\begin{align}\label{ce7}
|D_s(q,\hzeta,\rho)|\geq c\rho\text{ for }\hzeta\in\oS^d_+\text{ and
}0<\rho\leq\rho_0.
\end{align}

3.  The profile $W(z,q)$ satisfies the \emph{standard uniform Evans
condition at $q$} if and only if in addition to \eqref{ce7} we have
$D_s(q,\hzeta,\rho)\neq 0$ for $\hzeta\in\oS^d_+\text{ and }\rho>0$.

\end{defn}
Defining $V_\pm$ by $U_\pm=Y_\pm V_\pm$ as in \eqref{c6f}, we
transform \eqref{ce1} to
\begin{align}\label{cea2}
\partial_zV=G_\pm(q,\zeta)V\text{ on }\pm z\geq 0,\;\tilde{\Gamma}(q,\hzeta,\rho)(V_+,V_-)=0\text{ on
}z=0,
\end{align}
where
$\tilde{\Gamma}(q,\hzeta,\rho)(V_+,V_-):=\Gamma_s(Y_+V_+,Y_-V_-)$.
The argument of \cite{MZ3}, Lemma 2.5, shows that for
$\hzeta\in\oS^d_+$, $\rho>0$, each of $G_\pm(q,\zeta)$ has $N$
eigenvalues counted with their multiplicities in $\Re\mu>0$ and $N$
eigenvalues in $\Re\mu<0$. Together with the properties of
$Y_\pm(0,q,\zeta)$, this implies \eqref{cc3} and the smooth
dependence of $E_\pm(q,\hzeta,\rho)$ on $(q,\hzeta,\rho)$.

Next we give an alternative form of $D_s(q,\hzeta,\rho)$ for low
frequencies. For $\rho$ small conjugation of \eqref{cea2} to HP form
\eqref{c6i} using
\begin{align}\label{cec2}
V_\pm=\Lambda_\pm(q,\zeta)\begin{pmatrix} u_{H\pm}\\
u_{P\pm}\end{pmatrix}
\end{align}
transforms it to
\begin{align}\label{ce4}
\begin{split}
&\partial_z\begin{pmatrix}u_{H\pm}\\u_{P\pm}\end{pmatrix}=\begin{pmatrix}H_\pm(q,\zeta)&0\\0&P_\pm(q,\zeta)\end{pmatrix}\begin{pmatrix}u_{H\pm}\\u_{P\pm}\end{pmatrix}\\
&\tGamma_H(q,\zeta)u_H+\tGamma_P(q,\zeta)u_P=0\text{ on }z=0,
\end{split}
\end{align}
where
\begin{align}\label{cc7}
\tGamma_H(q,\zeta)u_H=\begin{bmatrix}T_{11}(0,q,\zeta)
u_H\\T_{21}(0,q,\zeta)u_H\end{bmatrix},\;\;\;\tGamma_P(q,\zeta)u_P=\begin{bmatrix}T_{12}(0,q,\zeta)
u_P\\T_{22}(0,q,\zeta)u_P\end{bmatrix}
\end{align}
(compare \eqref{c15}).

\begin{rem}\label{ce8}
1.  \textup{We note that for $q=(p_+,p_-,s,h)$,
$\uq=(\up_+,\up_-,0,0)$, $P_{0\pm}(q)$ as in \eqref{c16t}, and
$H_{0\pm}(q,\hzeta)$ as in \eqref{c39}, we have
\begin{align}\label{ce9}
\begin{split}
&(a)\;H_\pm(q,\zeta)=\rho\check{H}_\pm(q,\hzeta,\rho),\;
\check{H}_\pm(q,\hzeta,\rho)=H_{0\pm}(q,\hzeta)+O(\rho)\\
&(b)\;P_\pm(q,\zeta)=P_{0\pm}(q)+O(\rho).
\end{split}
\end{align}}

2.  \textup{For $\rho>0$ small  let us set
\begin{align}\label{ce10}
\begin{split}
&(a)\bE_-(H(q,\zeta)):=\bE_-(H_+(q,\zeta))\times\bE_+(H_-(q,\zeta))\subset\bC^{2N}\\
&(b)\bE_-(P(q,\zeta)):=\bE_-(P_+(q,\zeta))\times\bE_+(P_-(q,\zeta))\subset\bC^{2N}.
\end{split}
\end{align}}

3.   \textup{For $\check{H}_\pm(q,\hzeta,\rho)$ as in \eqref{ce9},
$\hzeta\in\overline{S}^d_+$, and $\rho>0$ small  we clearly have
\begin{align}\label{ce11}
\bE_-(H(q,\zeta))=\bE_-(\check{H}(q,\hzeta,\rho)).
\end{align}}

\end{rem}

\begin{prop}\label{cc8}
1.  The spaces appearing in \eqref{ce10}(a) have dimensions $N-1+k$,
$N-R_-$, and $N-L_+$, respectively.

2.   The spaces appearing in \eqref{ce10}(b) have dimensions
$N+1-k$, $R_-$, and $L_+$, respectively.

3.  Analogously, one can define spaces $\bE_+(H(q,\zeta))$ and
$\bE_+(P(q,\zeta))$ of dimensions $N+1-k$ and $N-1+k$ respectively.

4.  The spaces in \eqref{ce10} are $C^\infty$ in $(q,\hzeta,\rho)$
for $\hzeta\in\oS^d_+$ and $\rho>0$ small.
\end{prop}

\begin{proof}
The dimensions of $\bE_\mp(P_\pm(q,\zeta))$ follow directly from
Defn. \ref{b19} and the fact that
\begin{align}\label{cc9}
P_\pm(q,\zeta)=b(h)\cA_d(p_\pm,s,h)+O(\rho)  \;\;\eqref{c6j}.
\end{align}
The dimensions of $\bE_\mp(H_\pm(q,\zeta))$ then follow immediately
from
\begin{align}\label{cc10}
\dim \bE_-(G_+(q,\zeta))=N,\;\dim \bE_+(G_-(q,\zeta))=N.
\end{align}
The $C^\infty$ dependence of the spaces on $(q,\hzeta,\rho)$ follows
from the absence of pure imaginary eigenvalues of $H_\pm(q,\zeta)$
and $P_\pm(q,\zeta)$ for $\hzeta\in\oS^d_+$, $\rho>0$.

\end{proof}

\begin{rem}\label{ce5}
    \textup{Using the properties of the conjugators $T_\pm$ and
\eqref{ce4}, we see that up to a $C^\infty$ factor bounded away from
$0$, the \emph{standard Evans function} is given for $\hzeta\in
\oS^d_+$ and $\rho>0$ small by the $4N\times 4N$ determinant:
\begin{align}\label{cee6}
\tD_s(q,\hzeta,\rho):=\det\left(\bE_-(H(q,\zeta))\times\bE_-(P(q,\zeta)),\ker\tGamma_{H,P}(q,\zeta)\right),
\end{align}
where $\tGamma_{H,P}:\bC^{2N}\times\bC^{2N}\to\bC^{2N}$ is given
by
\begin{align}\label{ce6}
\tGamma_{H,P}(q,\zeta)(u_H,u_P):=\tGamma_H(q,\zeta)u_H+\tGamma_P(q,\zeta)u_P.
\end{align}}

\end{rem}

The operators
\begin{align}\label{cee11}
T_\pm(z,q,\zeta):=Y_\pm(z,q,\zeta)\Lambda_\pm(q,\zeta)
\end{align}
map solutions of \eqref{ce4} to solutions of \eqref{ce1}. For
$\rho>0$ choose smooth bases
\begin{align}\label{ce12}
\{u_{H,j}^\pm(q,\hzeta,\rho)\},\;\;\{u_{P,k}^\pm(q,\zeta)\}
\end{align}
of $\bE^-(H(q,\zeta))$ and $\bE^-(P(q,\zeta))$, respectively, and
set
\begin{align}\label{ce13}
\begin{split}
&S_j^\pm(z,q,\hzeta,\rho)=\begin{pmatrix}s_j^\pm\\\partial_zs^\pm_j\end{pmatrix}=T_\pm(z,q,\zeta)\begin{pmatrix}e^{zH_\pm(q,\zeta)}u_{H,j}^\pm(q,\hzeta,\rho)\\0\end{pmatrix}\\
&F_k^\pm(z,q,\zeta)=\begin{pmatrix}f_j^\pm\\\partial_zf^\pm_j\end{pmatrix}=T_\pm(z,q,\zeta)\begin{pmatrix}0\\e^{zP_\pm(q,\zeta)}u_{P,j}^\pm(q,\zeta)\end{pmatrix}.
\end{split}
\end{align}
The set of functions $\{S_j^+,j=1,\dots,N-R_-\}\cup\{F^+_k,
k=1,\dots,R_-\}$ is a basis for the space of solutions to
\eqref{ce1}(a) which decay to zero, when $\rho>0$, as
$z\to+\infty$.   We refer to the $S_j^+$, $F^+_k$ as \emph{slow
modes} and \emph{fast modes} respectively. Similar statements
apply to the $S^-_j$, $F^-_k$, where now $j=1,\dots,N-L_+$ and
$k=1,\dots, L_+$.

\begin{rem}\label{cee4}
1.  \textup{Define the $2N\times 2N$ determinant (suppressing
evaluation at $(0,q,\hzeta,\rho)$)
\begin{align}\label{ce15}
\bD_s(q,\hzeta,\rho):=\det\begin{pmatrix}S^+_j&S^-_k&F^+_l&F^-_m\end{pmatrix}.
\end{align}
Performing a few row and column operations
%and using \eqref{ce13}
shows that $\bD_s$ is equal, up to a sign, to the $4N\times 4N$
determinant
\begin{align}\label{ce16}
\det\begin{pmatrix}S^+_j&F^+_l&0&0&\te_n&\tf_n\\0&0&S^-_k&F^-_m&\te_n&\tf_n\end{pmatrix},
\end{align}
where, with $\{e_n,n=1,\dots,N\}$ a basis for $\bC^N$, we've set
$\te_n=\begin{pmatrix}e_n\\0\end{pmatrix},\;\tf_n=\begin{pmatrix}0\\e_n\end{pmatrix}$.
But the determinant in \eqref{ce16} is just $D_s(q,\hzeta,\rho)$
\eqref{cc2}, so
%\begin{align}\label{ce16z}
%\alpha(q,\zeta)=\det T_+(0,q,\zeta)\cdot\det T_-(0,q,\zeta),
%\end{align}
%and $D_s$ is given by \eqref{ce6}.   Since $\alpha(q,\zeta)$ is
%smooth and bounded away from zero for $\rho$ small,
the low frequency standard Evans condition can equivalently be
expressed in terms of $\bD_s$.}

%$=\alpha(q,\zeta)D_s(q,\hzeta,\rho)$

2.  \textup{For $\hgamma>0$ the functions
$u^\pm_{H,j}(q,\hzeta,\rho)$ can be chosen to extend smoothly to
$[0,\rho_0)$, since $H_{0,\pm}(q,\hzeta)$ \eqref{ce9} has no pure
imaginary eigenvalues when $\hgamma>0$. Hence the same is true for
the Evans functions $D_s$ and $\bD_s$.  Note also that the
$u^\pm_{H,j}(q,\hzeta,0)$ span $\bE_\mp(H_{0\pm}(q,\hzeta)$.}

3.   \textup{a-transversality is equivalent to the condition that
the kernel of
\begin{align}\label{ce18}
\tilde{\Gamma}_P(q,0):\bE_-(P(q,0))\to\bC^{2N}
\end{align}
be of dimension one.  We know that $(c_+(q),c_-(q))$ is in the
kernel, where
\begin{align}\label{ce19}
\begin{pmatrix}0\\c_\pm(q)\end{pmatrix}:=T_\pm^{-1}(0,q,0)\begin{pmatrix}W_z(0,q)\\W_{zz}(0,q)\end{pmatrix},
\end{align}
so $D_s(q,\hzeta,\rho)$ must vanish to at least first order at
$\rho=0$ when $\hgamma>0$.  When the low frequency Evans condition
holds, $D_s$ vanishes to precisely first order at $\rho=0$, so weak
transversality holds.}

\end{rem}

\subsection{The modified Evans function}\label{modified} \textbf{}

The translational degeneracy in the partially linearized problem
represented by the nontrivial kernel of \eqref{ce18} is a serious
obstacle to proving robust $L^2$ estimates. Thus, we are led as in
\cite{GMWZ3} to consider the fully linearized parabolic problem with
an extra transmission condition
\begin{align}\label{c46}
\begin{split}
&\cL(z,q,\zeta,\partial_z)u-\psi\cL_1(z,q,\zeta)=f\\
&[u]=0,\;[u_z]=0, \;\uw_z(0)\cdot u^++c_0(\zeta)\psi=0.
\end{split}
\end{align}
The operators the operators $\cL$ and $\cL_1$ are given explicitly
in \eqref{c5a} and $c_0(\zeta)=i\tau+\gamma+|\eta|^2$. The choice of
the third transmission condition here is related to the choice in
\eqref{c4} and \eqref{c5}; again, it serves to remove the
translational degeneracy.  In this section we define a modified
Evans function for \eqref{c46} which turns out to be bounded away
from zero for $\rho>0$ small when the standard low frequency Evans
condition \eqref{ce7} holds.

Setting
\begin{align}\label{c46a}
\cL_1(z,q,\zeta)=\rho\check{\cL}_1(z,q,\hzeta,\rho),\;\;c_0(\zeta)=\rho\chc_0(\hzeta,\rho),\;\;\phi=\rho\psi
\end{align}
we can rewrite \eqref{c46}:
\begin{align}\label{c47}
\begin{split}
&\cL(z,q,\zeta,\partial_z)u-\phi\check{\cL}_1(z,q,\hzeta,\rho)=f\text{
on }\pm z\geq 0,\\
&[u]=0,\;[u_z]=0, \;\uw_z(0)\cdot
u^++\check{c}_0(\hzeta,\rho)\phi=0.
\end{split}
\end{align}

Let $ \widetilde \EE(q, \hzeta, \rho)$ denote the space of triples
$(U^-_0, U^+_0,\phi) \in \CC^{2N} \times \CC^{2N} \times \CC $
with $U^\pm_0 = (u^\pm_0, v^\pm_0)$, such that the solutions $u^
\pm $ of
$$
\cL u^\pm - \phi \check{\cL}_1 = 0 \text{ on } \pm z \ge 0\, ,
\qquad u^\pm (0) = u^\pm_0 \,, \quad \D_z  u^\pm (0) = v^\pm_0
$$
are  bounded at infinity. Let $\ker \widetilde
\Gamma(q,\hzeta,\rho) $ denote the set of $(U^-_0, U^+_0,\phi) \in
\CC^{2N} \times \CC^{2N} \times \CC $ such that
\begin{equation}
\label{c47y} U^-_0 = U^+_0  \, , \qquad \uw_z(0) \cdot u^+(0)
+\chc_0(\hzeta, \rho) \phi = 0\,.
\end{equation}

\begin{defn}\label{c47z}
The \emph{modified Evans function}  is the $(4N+1)\times (4N+1)$
determinant
\begin{equation}\label{c47a}
\widetilde D(q, \hzeta, \rho) = \det \Big(\widetilde \EE(q,
\hzeta, \rho), \ker \widetilde \Gamma (q,\hzeta,\rho) \Big) \,.
\end{equation}
\end{defn}

Parallel to what we did in Proposition \ref{c6k} and Corollary
\ref{c6s} for the linearized profile equations, in the  next
Proposition we recast \eqref{c47} in an equivalent form where the
operator $\check{\cL}_1$ no longer appears.    For this we need a
good  extension of $W_z(z,q)$ to nonzero frequencies.

We recall that
\begin{align}\label{p1}
\cW_\pm(z,q,0):=(W_z(z,q),W_{zz}(z,q))\text{ on }\pm z\geq 0
\end{align}
satisfies
\begin{align}\label{p2}
\partial_z\cW=G(z,q,0)\cW\text{ on }\pm z\geq
0,\;\Gamma_s(\cW_+,\cW_-)=0\text{ on }z=0.
\end{align}
We can smoothly extend $\cW_\pm$ to $|\zeta|$ small as solutions of
\begin{align}\label{p3}
\partial_z\cW=G(z,q,\zeta)\cW\text{ on }\pm z\geq 0
\end{align}
by setting, for $T_\pm$ as in \eqref{cee11},
\begin{align}\label{p4}
\cW_\pm(z,q,\zeta)=T_\pm(z,q,\zeta)\begin{pmatrix}0\\e^{zP_\pm(q,\zeta)}\pi_\pm(q,\zeta)c_\pm(q)\end{pmatrix}.
\end{align}
Here $\pi_\pm(q,\zeta)$ are, respectively, the  projections of
$\bC^N$ onto $\bE_-(P_+(q,\zeta))$, $\bE_+(P_-(q,\zeta))$ (along
$\bE_+(P_+(q,\zeta))$, $\bE_-(P_-(q,\zeta))$ respectively), and
$c_\pm(q)$ are elements of $\bE_{\mp}(P_\pm(q,0))$ defined by
\eqref{p4} at $\zeta=0$.

\begin{prop}\label{c53a}
For $\hzeta\in\oS^d_+$ and $0\leq\rho\leq\rho_0$ small, there exist
$C^\infty$ functions $\bfR_\pm(z,q,\zeta)$ on $\pm z\geq 0$,
exponentially decaying to zero as $z\to\pm\infty$,  and satisfying:
\begin{equation}
\label{c53} \left\{ \begin{aligned} & \cL(z,q,  \zeta, \partial_z)
\bfR = \cL_1(z, q, \zeta)   \quad \mathrm{on}  \  \pm z \ge 0 ,
\\
&  \uw_z(0) \cdot \bfR_ \pm(0, q, \zeta  ) =  -
c_0(\zeta):=-(i\tau+\gamma+|\eta|^2)\,, \quad \bfR^\pm (z, q, 0) = 0
\,.
\end{aligned} \right.
\end{equation}
Define $\check{\bfR}_\pm(z,q,\hzeta,\rho)$ by
\begin{align}\label{c53z}
\bfR(z,q,\zeta)=\rho\check{\bfR}_\pm(z,q,\hzeta,\rho).
\end{align}
The functions $\bfR_\pm(z,q,\zeta)$ can be constructed so that
\begin{align}\label{c54}
\cR_\pm(z,q,i\htau+\hgamma,i\heta):=\check{\bfR}_\pm(z,q,\hzeta,0)
\end{align}
is linear in $(i\htau+\hgamma,i\heta)$.  It then makes sense to
consider $\cR_\pm(z,q,\dos,\doh)$ for $(\dos,\doh)\in\bR^d$, and
these functions have all the properties of the functions constructed
in Proposition \ref{c6k}. In particular, we have
\begin{align}\label{c55}
\uw_z(0)\cdot\cR_\pm(0,q,i\htau+\hgamma,i\heta)=\uw_z(0)\cdot\check{\bfR}_\pm(0,q,\hzeta,0)=-(i\htau+\hgamma).
\end{align}
\end{prop}

\begin{proof}
The proof is quite similar to that of Proposition \ref{c6k}.    In
particular, the first $N$-dimensional component of $\cW_\pm
(z,q,\zeta)$ \eqref{p4} plays the role of $W_z(z,q)$ in the earlier
proof. We construct
\begin{align}\label{c56}
\textit{\bfR}(z,q,\zeta)=S_0(z,q,\zeta)(i\tau+\gamma)+\sum^{d-1}_{j=1}S_j(z,q,\zeta)i\eta_j,
\end{align}
where the $S_j^\pm$ are $C^\infty$ functions chosen to satisfy:
\begin{align}\label{c56z}
\begin{split}
&\cL(z,q,\zeta,\partial_z)S^\pm_0=A_0(W)W'\text{ on }\pm z\geq 0\\
&\cL(z,q,\zeta,\partial_z)S^\pm_j=A_j(W)W'+2h_jW''-i\eta_jW',\;j=1,\dots,d-1.\\
&\uw_z(0)\cdot S^\pm_0(0,q,\zeta)=-1,\;\;\uw_z(0)\cdot
S^\pm_j(0,q,\zeta)=i\eta_j.
\end{split}
\end{align}
The interior forcing terms are exponentially decaying, so
exponentially decaying functions $S_j$ with these properties are
readily constructed by using the conjugators $T_\pm$ to reduce to
$HP$ form \eqref{ce4} (see \cite{GMWZ3}, Lemma 3.14 for details).
The functions
\begin{align}\label{c56y}
\cR(z,q,\dos,\doh)=S_0(z,q,0)\dos+\sum^{d-1}_{j=1}S_j(z,q,0)\doh
\end{align}
then have all the properties of the functions constructed in
Proposition \ref{c6k}.
\end{proof}

\begin{rem}\label{c56x}
\textup{Henceforth, we'll use the functions $\cR(z,q,\dos,\doh)$ as
in \eqref{c56y} in place of the functions constructed in Proposition
\ref{c6k}.}
\end{rem}

With $\ibfR^\pm$ as in \eqref{c53z}, for $0\leq\rho\leq \rho_0$ the
problem \eqref{c47} is equivalent to
 \begin{equation}
\label{c48} \left\{
\begin{aligned}
& u^\pm = v^\pm + \phi \ibfR^\pm
 \\
& \cL(z, q, \zeta, \partial_z) v^\pm   = f^\pm \, \quad \mathrm{
on } \;\pm z \ge 0 \,,
\\
& [v (0) ] + \phi [\ibfR(0)] = 0 \quad
 [ \partial_z v (0) ] + \phi [\partial_z \ibfR (0)]  = 0 ,
\\
 &\uw_z(0) \cdot v^+ (0)  = 0 \,.
\end{aligned}\right.
\end{equation}

As before let $\bE_\pm(q,\hzeta,\rho)$ be the set of
$V^\pm=(v^\pm_0,v^\pm_1)\in\bC^{2N}$ such that the solutions of
\begin{align}
\cL v^\pm=0\text{ on }\pm z\geq 0,\; v^\pm(0)=v^\pm_0,\;
\partial_zv^\pm(0)=v^\pm_1
\end{align}
are bounded at $\pm\infty$. Then,
$$
\widetilde \EE (q, \hzeta, \rho) = \cJ \big(\EE_- \times \EE_+
\times \CC \big)
$$
where
$$
\cJ (q, \hzeta, \rho)  : \quad (V^-, V^+, \phi) \mapsto (V^- +
\phi \rmR^-, V^+ + \phi \rmR^+, \phi)
$$
with
$$
\rmR^\pm (\hzeta, \rho)  = {}^t \big( \ibfR^\pm(0) ,
\partial_z \ibfR^\pm (0) \big) \,.
$$

Moreover,  $\ker \widetilde \Gamma  = \cJ \GG' $  with
$$
\GG'(q, \hzeta, \rho) = \big\{ (V^-, V^+,
\phi)\in\bC^{2N}\times\bC^{2N}\times\bC : V^+ - V^- = \phi ( \rmR^-
- \rmR^+) , \   \ell \cdot V^+ = 0 \big\}
$$
where $ \ell = {}^t( \uw_z(0), 0)$. Therefore, the Evans function of
the problem \eqref{c48}
\begin{equation}
\label{c49}
  D_m(q, \hzeta, \rho) :=
\det \Big(\EE_- \times \EE_+ \times \CC ,  \GG'  \Big)
\,,\;0<\rho\leq\rho_0
\end{equation}
satisfies:
\begin{equation}
 \label{c50}
 \frac{1}{C}  \vert \widetilde D(q, \hzeta, \rho) \vert \le
\vert  D_m (q, \hzeta, \rho) \vert
  \le {C} \vert  \widetilde D(q, \hzeta,\rho ) \vert\text{ for }0<\rho\leq \rho_0.
 \end{equation}

\begin{defn}\label{c51}
1.  The profile $W(z,q)$ satisfies the \emph{modified low frequency
Evans condition at $q$} when there exist positive constants $c$ and
$\rho_0$ such that
\begin{align}\label{c52}
|\tilde{D}(q,\hzeta,\rho)|\geq c \text{ for }0<\rho\leq\rho_0,
\;\hzeta\in\overline{S}^d_+.
\end{align}
Equivalently, one can replace $\tilde{D}$ by $D_m$ in \eqref{c52}.

2.   The profile $W(z,q)$ satisfies the \emph{modified uniform Evans
condition at $q$} when in addition to \eqref{c52} we have
$\tilde{D}(q,\hzeta,\rho)\neq 0$ for $\hzeta\in\oS^d_+$, $\rho>0$.

3.  Henceforth, it will be convenient to define
$D_m(q,\hzeta,\rho):=\tilde{D}(q,\hzeta,\rho)$ in $\rho>\rho_0$ and
to drop the notation $\tilde{D}$.

\end{defn}

\begin{rem}\label{c51z}
\textup{Observe that no transversality assumptions are needed to
define the standard and modified Evans functions.}
\end{rem}

\subsection{Inviscid and viscous continuity}\label{continuity}

In view  of the noncompactness of $S^d_+$ the question arises as to
whether or not the stability conditions defined in the last few
sections must necessarily hold near $\uq$ if they hold at $\uq$.  To
address this and also for later use we give the following
definition:

\begin{defn}\label{j1}
(a) We say \emph{inviscid continuity} holds at $q_0\in\cC$  when the
vector bundle $\bE_-(H_0(q,\hzeta))$ has a continuous extension from
$S^d_+$ to $\oS^d_+$ for $q$ near $q_0$.

(b)We say \emph{viscous continuity} holds at $q_0\in\cC$  when the
vector bundle $\bE_-(\check{H}(q,\hzeta,\rho))$ has a continuous
extension from $\oS^d_+\times (0,\infty)$ to $\oS^d_+\times
[0,\infty)$ for $q$ near $q_0$.

\end{defn}

\begin{rem}\label{j2}
1.   \textup{Inviscid continuity holds for both fast and slow shocks
in inviscid MHD and more generally, for Friedrichs symmetrizable
systems for which all characteristic roots of the linearized
operator are either geometrically regular or totally nonglancing
(Theorem 5.6, \cite{MZ2}).   More generally still, inviscid
continuity is a necessary condition for the existence of a smooth
$K$-family of inviscid symmetrizers.}

2.   \textup{The above structural conditions on the hyperbolic part
are not enough to guarantee viscous continuity; for example, viscous
continuity can fail when there is viscous coupling between incoming
and outgoing crossing eigenvalues (Proposition 6.5, \cite{GMWZ5}).
The existence of a smooth $K$-family of viscous symmetrizers implies
viscous continuity. Viscous continuity holds for fast (i.e.,
extreme) shocks in viscous MHD, but not for slow shocks.}

3. \textup{Since $H_0(q,\hzeta)=\check{H}(q,\hzeta,0)$, viscous
continuity implies inviscid continuity.}

4. \textup{Viscous continuity holds when the hyperbolic problem has
characteristics of constant multiplicity \cite{MZ1}.}

5. \textup{When $d=1$ variable multiplicities are impossible, so we
always have viscous continuity.}

\end{rem}

\begin{rem}\label{j31}

  \textup{1.\;\;If both the uniform Lopatinski condition and inviscid
continuity hold at $\uq$, then it is clear by compactness of
$\oS^d_+$ and continuity that the uniform Lopatinski condition holds
for $q$ near $\uq$.}
%In this case we will be able to construct
%curved inviscid shocks near a fixed $\uq$ assuming only the
%$\chi$-Lopatinski condition at $\uq$.}
   \textup{Similarly, if the low frequency standard Evans condition
and viscous continuity hold at $\uq$, then the standard low
frequency Evans condition holds for $q$ near $\uq$.}

%In this case we can prove the vanishing viscosity convergence of
%viscous shocks to a curved inviscid shock near a fixed planar $\uq$
%assuming only the low frequency standard Evans condition at $\uq$
%(and nonvanishing of $D(\uq,\hzeta,\rho)$ for $\rho>0$).}

 \textup{2.\;\;The determinants $D_{Lop}(q,\hzeta)$ and
$D_{Lop,m}(q,\zeta)$ were defined for $q$ belonging to some shock
manifold $\cC$ (or $\cC_\cB$).  The same definition makes sense for
all $q$ in a small $\bR^{2N+d}$ neighborhood of $\cC$.  Similarly,
if $q=(p_+,p_-,s,h)\in\cC_\cB$ and $(u'_+,u'_-,s',h')\in\bR^{2N+d}$
is sufficiently small, the definitions of $D_s(q,\hzeta,\rho)$ and
$D_m(q,\hzeta,\rho)$ make sense if the linearizations are computed
at
\begin{align}\label{j31z}
(W_\pm(z,q)+u'_\pm, s+s', h+h')
\end{align}
instead of at $(W(z,q),s,h)$.     This extension is important for
obtaining estimates in the linearized problems that arise in the
iteration schemes used to prove the nonlinear stability theorems.}

\end{rem}

\section{Low frequency analysis of the Evans functions}\label{low
frequency}

First,  we present an important technical lemma that is needed for
the low frequency analysis of the standard Evans function; in
particular, it is used in the proofs of Theorems \ref{K} and
\ref{g3}.

For $\cW_\pm(z,q,\zeta)$ as constructed in \ref{p4}, let us define
the variations of the profile at $\rho=0$:
\begin{align}\label{p5}
\cZ_\pm(z,q,\hzeta,0):=\partial_\rho|_{\rho=0}\cW_\pm(z,q,\zeta),\quad
\cZ_\pm=(\cZ^1_\pm,\cZ^2_\pm)\in\bC^{2N}.
\end{align}
The lemma permits us to replace $\cZ^1_\pm(z,q,\hzeta,0)$ by
$-\cbR_\pm(z,q,\hzeta,0)$ in some computations.

\begin{lem}[Variation of the extended $W_z$ at $\rho=0$]\label{p6}

(1)  For $\cZ^1_\pm$ as in \eqref{p5} and $\cbR$ as in \eqref{c48}
we have
\begin{align}\label{p7}
\cZ^1_\pm(z,q,\hzeta,0)+\cbR_\pm(z,q,\hzeta,0)=T_{12\pm}(z,q,0)e^{zP_\pm(q,0)}u_{P,z,\pm}(q,\hzeta)
\end{align}
for some $u_{P,z\pm}(q,\hzeta)\in\bE_\mp(P_\pm(q,0))$.

(2)  For $\cZ^1_\pm$ as in \eqref{p5} and $\phi_\pm$ as in Prop.
\ref{c26} we have
\begin{align}\label{p7z}
\cZ^1_\pm(z,q,\hzeta,0)+\nabla_{s,h}\phi_\pm(z,p_\pm,s,h,a_\pm)(i\htau+\hgamma,i\heta)=\nabla_{a\pm}\phi_\pm(z,p_\pm,s,h,a_\pm)
\doa_{z\pm}(q,\hzeta)
\end{align}
for some $\doa_{z\pm}(q,\hzeta)\in\bE_\mp(P_{\pm}(q,0)).$

\end{lem}

\begin{proof}
\textbf{1. }  Noting that
$\cbR_\pm(z,q,\hzeta,0)=\partial_\rho|_{\rho=0}\bfR_\pm(z,q,\zeta)$,
we apply $\partial_\rho|_{\rho=0}$ to the equation (recall
\eqref{c53})
\begin{align}\label{p9}
\cL(z,q,\zeta,\partial_z)\bfR_\pm=\cL_1(z,q,\zeta)
\end{align}
to find (for $\cL_0(z,q,\partial_z)=\cL(z,q,0,\partial_z)$ as in
\eqref{c5})
\begin{align}\label{p10}
\cL_0(z,q,\partial_z)\cbR_\pm(z,q,\hzeta,0)=A_0(W)W'(i\htau+\hgamma)+\sum^{d-1}_{j=1}A_j(W)W'i\heta_j+2\sum^{d-1}_{j=1}h_ji\heta_jW''.
\end{align}

\textbf{2. }Applying $\partial_\rho|_{\rho=0}$ to
\begin{align}\label{p10z}
\cL(z,q,\zeta,\partial_z)\cW_\pm(z,q,\zeta)=0
\end{align}
and using \eqref{p1}, we find $-\cZ^1_\pm(z,\uq,\hzeta,0)$ satisfies
the same equation:
\begin{align}\label{p11}
-\cL_0(z,q,\partial_z)\cZ^1_\pm(z,q,\hzeta,0)=A_0(W)W'(i\htau+\hgamma)+\sum^{d-1}_{j=1}A_j(W)W'i\heta_j+2\sum^{d-1}_{j=1}h_ji\heta_jW''.
\end{align}
For each $\hzeta$ the sum
$\cZ^1_\pm(z,q,\hzeta,0)+\cbR_\pm(z,q,\hzeta,0)$ is thus a solution
of $\cL_0(z,q,\partial_z)\dov=0$ which decays exponentially to zero
as $z\to\pm\infty$; hence it must be given by the right side of
\eqref{p7} for some $u_{P,z,\pm}(q,\hzeta)\in\bE_\mp(P_\pm(q,0))$
(recall \eqref{c10}).

\textbf{3. }Similarly, arguing as in the proof of Prop. \ref{c26},
we see that the image of
\begin{align}\notag
\nabla_{s,h}\phi_\pm(z,p_\pm,s,h,a_\pm)(i\htau+\hgamma,i\heta)
\end{align}
under $\cL_0(z,q,\partial_z)$ is again equal to the right side of
\eqref{p11}.  So Proposition \ref{c26} implies \eqref{p7z}.

\end{proof}

\subsection{Nonconservative Zumbrun-Serre Theorem}\label{ZS}

In the next Theorem we show that the standard low frequency Evans
condition \eqref{ce7} implies transversality and the uniform
Lopatinski condition.   Moreover, if one assumes viscous continuity
(Defn. \ref{j1}), then the converse holds.    Here we'll use the
functions $\tPsi(p,s,h,a)$ defined in \eqref{b28} as well as
Proposition \ref{b30}. In particular, we'll use the fact that
transversality of $\uw$ allows one to define $\chi$ and the manifold
$\cC_\cB$ given by $\chi(p,s,h)=0$.

%In the next theorem we establish a relation between the standard
%Evans function $\bD_s(q,\hzeta,\rho)$ defined in \eqref{ce15} and
%the Lopatinski determinant $D_{Lop,\chi}(q,\hzeta)$ \eqref{c40c}.
%We make use of the function $\tilde{\Psi}(p,s,h,a)$ \eqref{b28}
%and the manifold $\cC$ defined by $\chi(p,s,h)=0$ (see Prop.
%\ref{b30} and Remark \ref{b56y}, part 1) when $\uw(z)$ is
%transversal.

\begin{thm}\label{K}
(a).  Consider the shock profile $\uw(z)=W(z,\uq)$, where
$\uq=(\up,0,0)$,  and suppose the low frequency standard Evans
condition holds at $\uq$. That is, suppose there are positive
constants $c$ and $\rho_0$ such that
\begin{align}\label{ceq15}
|\bD_s(\uq,\hzeta,\rho)|:=|\det\begin{pmatrix}S^+_j&S^-_k&F^+_l&F^-_m\end{pmatrix}|\geq
c\rho
\end{align}
for all $\hzeta\in\oS^d_+$ and $0<\rho\leq \rho_0$.  Then $\uw(z)$
is transversal (Defn. \ref{c11})
%\begin{align}\label{cf1}
%\mathrm{rank}\nabla_{a,p,s}\tilde{\Psi}(\up,0,0,\ua)=2N+1,
%\end{align}
 and the uniform Lopatinski condition
\eqref{c40b} holds at $\uq$.  In fact, for $\hgamma>0$
\begin{align}\label{ceq16}
\bD_s(\uq,\hzeta,\rho)=\rho\beta(\uq)D_{Lop}(\uq,\hzeta)+O_{\hgamma}(\rho^2),
\end{align}
where $O_{\hgamma}(\rho^2)\leq C_{\hgamma}\rho^2$ and $\beta(\uq)$
is a constant whose nonvanishing is equivalent to a-transversality
of $\uw$.

(b).   Assume $\uw(z)$ is transversal.
% and that
%$\nabla_{a,p,s}\tilde{\Psi}(\up,0,0,\ua)$ has full rank $2N+1$.
If
$\bE^-(\check{H}(\uq,\hzeta,\rho))$ has a continuous extension from
$\oS^d_+\times (0,\rho_0)$ to $\oS^d_+\times [0,\rho_0)$ for some
$\rho_0>0$, then
\begin{align}\label{ceqa16}
\bD_s(\uq,\hzeta,\rho)=\rho\beta(\uq)D_{Lop}(\uq,\hzeta)+o(\rho),
\end{align}
with an error term that is uniform for $\hzeta\in\oS^d_+$.  Thus, if
the uniform Lopatinski condition also holds, then the low frequency
standard Evans condition holds.

\end{thm}

\begin{proof}

\textbf{1.  Part (a): Reduce to $\partial_\rho\bD_s(\uq,\hzeta,0)$.
} We can suppose $F^+_1(z,\uq,\zeta)$ and $F^-_{L_+}(z,\uq,\zeta)$
are given respectively by $\cW_+(z,\uq,\zeta)$  and
$\cW_-(z,\uq,\zeta)$ \eqref{p4}, so that
\begin{align}\label{k3}
f^+_1(z,\uq,0)=\uw_z(z)\text{ on }z\geq
0;\;f^-_{L+}(z,\uq,0)=\uw_z(z)\text{ on }z\leq 0.
\end{align}
Set $Z_\pm(z,\uq,\hzeta,0)=(z_\pm,\partial_z z_\pm)$, where
\begin{align}\label{k4}
z_+(z,\uq,\hzeta,0)=\partial_\rho|_{\rho=0}f^+_1(z,\uq,0);\;z_-(z,\uq,\hzeta,0)=\partial_\rho|_{\rho=0}f^-_{L+}(z,\uq,0).
\end{align}
Next, fix $\hgamma>0$ and subtract the $F^+_1$ column from the
$F^-_{L+}$ column to obtain
\begin{align}\label{k5}
\begin{split}
&\bD_s(\uq,\hzeta,\rho)=\rho\det\begin{pmatrix}S^+_j(0,\uq,\hzeta,0)&S^-_k&F^+_l(0,\uq,0)&F^-_m&(Z_--Z_+)(0,\uq,\hzeta,0)\end{pmatrix}+O_{\hgamma}(\rho^2),\\
&\qquad\qquad:=\rho\bD^*_s(\uq,\hzeta)+O_{\hgamma}(\rho^2).
\end{split}
\end{align}
where now $m=1,\dots,L_+-1$.  It remains to analyze
$\bD^*_s(\uq,\hzeta)$.   Observe that the assumption \eqref{ceq15}
implies
\begin{align}\label{k5a}
|\bD^*_s(\uq,\hzeta)|\geq c \text{ for all } \hzeta\in S^d_+.
\end{align}

\textbf{2.  Rewrite in terms of $\Psi$. }  Recall
\begin{align}\label{k2}
{\Psi}(p,s,h,a)=\begin{pmatrix}\phi_+-\phi_-\\\phi_{+,z}-\phi_{-,z}\end{pmatrix}(0,p,s,h,a)
\end{align}

For $\cL_0(z,q,\partial_z)$ as in \eqref{c5} we have
\begin{align}\label{k6}
\begin{split}
&\cL_0(z,\uq,\partial_z)s^+_j(z,\uq,\hzeta,0)=0\text{ and
}\lim_{z\to+\infty}s^+_j(z,\uq,\hzeta,0)=u^+_{H,j}(\uq,\hzeta,0)\in\bE_-(H_{0+}(\uq,\hzeta)),\\
&\cL_0(z,\uq,\partial_z)f^+_l(z,\uq,0)=0\text{ and
}\lim_{z\to+\infty}f^+_l(z,\uq,0)=0,
\end{split}
\end{align}
so Proposition \ref{c26} implies
\begin{align}\label{k7}
\begin{split}
&s^+_j(z,\uq,\hzeta,0)=\nabla_{p+}\phi_+(z,\up_+,0,0,\ua_+)u^+_{H,j}(\uq,\hzeta,0)+\nabla_{a+}\phi_+(z,\up_+,0,0,\ua_+) \doa_{j+}(\uq,\hzeta)\\
&f^+_l(z,\uq,0)=\nabla_{a+}\phi_+(z,\up_+,0,0,\ua_+)
\doa_{l+}(\uq)
\end{split}
\end{align}
for some $\doa_{j+}(\uq,\hzeta)$ and $\doa_{l+}(\uq)$ in
$\bE_-(P_{0+}(\uq))$. Analogous statements apply to the $s^-_k$
and $f^-_m$.

In addition, Lemma \ref{p6} implies
\begin{align}\label{k8}
z_+(z,\uq,\hzeta,0)+\nabla_{s,h}\phi_+(z,\up_+,0,0,\ua_+)(i\htau+\hgamma,i\heta)=\nabla_{a+}\phi_+(z,\up_+,0,0,\ua_+)
\doa_{z+}(\uq,\hzeta)
\end{align}
for some $\doa_{z+}(\uq,\hzeta)\in\bE_-(P_{0+}(\uq))$.   Again,
there is an analogue  for $z_-$.

We claim that, up to a sign,
\begin{align}\label{k9}
\begin{split}
&\bD^*_s(\uq,\hzeta)=\\
&\det\begin{pmatrix}\nabla_{p+}\Psi(\uq,\ua)u^+_{H,j}&\nabla_{p-}\Psi(\uq,\ua)u^-_{H,k}&\nabla_{a+}\Psi(\uq,\ua)\doa_{l+}&\nabla_{a_-}\Psi(\uq,\ua)\doa_{m-}&\nabla_{s,h}\Psi(\uq,\ua)(i\htau+\hgamma,i\heta)\end{pmatrix}.
\end{split}
\end{align}
Here we have used \eqref{k7}, \eqref{k8} and, for example, the
observation that

\begin{align}\label{k10}
\nabla_{a+}\begin{pmatrix}\phi_+\\\phi_{+,z}\end{pmatrix}(0,\up_+,0,0,\ua_+)=\nabla_{a+}\Psi(\uq,\ua).
\end{align}
Observe that  \eqref{k5a} implies the linear independence of the
$N-k$ columns of \eqref{k9} indexed by $l$ and $m$ (which is
equivalent to a-transversality of $\uw$).   In deriving \eqref{k9}
we have used those $N$ columns to eliminate the evaluations at $z=0$
of fast decaying terms like the one on the right in \eqref{k8}.

Taking $\heta=0$ in \eqref{k9} and using \eqref{k5a} again, we
conclude $\mathrm{rank}\nabla_{a,p,s}\Psi(\up,0,0,\ua)=2N$, and thus
(a,p,s)-transversality holds for $\uw$.

\textbf{3. Row and column operations. }Next we determine row and
column operations that reduce \eqref{k9} to the required block
form.  We show that after some preparation, these turn out to be
the same operations needed to compute $\chi'(\up,0,0)$ starting
from $\tilde{\Psi}'(\up,0,0,\ua)$.

Consider first the $(2N+1)\times(N+1-k)$ matrix
\begin{align}\label{k11}
\begin{pmatrix}\nabla_{a+}\tilde{\Psi}(\uq,\ua)\doa_{l+}(\uq)&\nabla_{a-}\tilde{\Psi}(\uq,\ua)\doa_{m-}(\uq)\end{pmatrix}=\begin{pmatrix}\uw_z(0)&*&*&-\uw_z(0)\\\uw_{zz}(0)&*&*&-\uw_{zz}(0)\\\uw_z(0)\cdot\uw_z(0)&*&0&0\end{pmatrix}
\end{align}
(recall \eqref{k3}).   By a linear change of coordinates in
$\bE_-(P_0(\uq))$ we arrange so that the matrix on the left in
\eqref{k11} is simply
\begin{align}\label{k12}
\begin{pmatrix}\nabla_{a+}\tilde{\Psi}(\uq,\ua)&\nabla_{a-}\tilde{\Psi}(\uq,\ua)\end{pmatrix}.
\end{align}
The $2N\times (N-k)$ submatrix of \eqref{k9} given by the columns
indexed by $l$ and $m$ then becomes
\begin{align}\label{k13}
\Psi_{\ta}(\uq,\ua):=\begin{pmatrix}\nabla_{a_+}\Psi(\uq,\ua)&\nabla_{\ta_-}\Psi(\uq,\ua)\end{pmatrix},
\text{ where }\ta_-=(a_{1-},\dots,a_{(L_+-1)-}).
\end{align}

Now, a-transversality implies that the equation
$\tilde{\Psi}(p,s,h,a)=0$  defines  a function $a(p,s,h)$ near
$(\up,0,0)$ such that for some choices $\tilde{\Psi}^1$ (resp.
$\tilde{\Psi}^2$) of $N+1-k$ (resp. $N+k$) components of
$\tilde{\Psi}$ we have
\begin{align}\label{k14}
\begin{pmatrix}\tilde{\Psi}^1(p,s,h,a(p,s,h))\\\tilde{\Psi}^2(p,s,h,a(p,s,h))\end{pmatrix}=\begin{pmatrix}0\\\chi(p,s,h)\end{pmatrix}.
\end{align}
In view of \eqref{k11} we can (and do) take the last component of
$\tilde{\Psi}$ to be the last component of $\tilde{\Psi}^1$.

Now differentiate both sides of \eqref{k14} at $(\up,0,0)$ to
obtain:
\begin{align}\label{k15}
\begin{pmatrix}0\\\chi_{p,s,h}\end{pmatrix}=\begin{pmatrix}0&0&0\\\chi_p&\chi_s&\chi_h\end{pmatrix}
=\begin{pmatrix}\tilde{\Psi}^1_{p,s,h}+\tilde{\Psi}^1_a(a_p,a_s,a_h)\\\tilde{\Psi}^2_{p,s,h}+\tilde{\Psi}^2_a(a_p,a_s,a_h)\end{pmatrix}.
\end{align}
If we start with the $(2N+1)\times (2N+d+N+1-k)$ matrix
$\tilde{\Psi}_{p,s,h,a}(\up,0,0,\ua)$, then \eqref{k15} provides
column operations  that transform this matrix (after row switches)
to
\begin{align}\label{k16}
\begin{pmatrix}0&0&0&\tilde{\Psi}^1_a\\\chi_p&\chi_s&\chi_h&\tilde{\Psi}^2_a\end{pmatrix}.
\end{align}
Inspection of \eqref{k11} shows that the transformation to
\eqref{k16} can (mostly) be performed without the help of the last
column of $\tilde{\Psi}_{p,s,h,a}(\up,0,0,\ua)$.   More precisely,
if we start with the $2N\times (2N+d+N-k)$ submatrix
$\Psi_{p,s,h,\ta}(\up,0,0,\ua)$, there are column operations that
reduce it after some row switches to the submatrix of \eqref{k16}
given by
\begin{align}\label{k18}
\begin{pmatrix}0&0&0&\Psi^1_{\ta}\\\chi_p&\chi_s&\chi_h&\Psi^2_{\ta}\end{pmatrix},
\end{align}
where $\Psi^2_{\ta}=\tPsi^2_{\ta}$ and $\Psi^1_{\ta}$ is the upper
left $(N-k)\times (N-k)$ block of $\tPsi^1_a$.

Since (after the earlier change of coordinates in $\bE_-(P_0(\uq))$)
$\Psi_{\ta}$ appears as a submatrix of \eqref{k9}, we conclude that
by adding appropriate combinations, which now involve coefficients
that depend on $\hzeta$, of the columns in \eqref{k13} to the other
columns of \eqref{k9}, we can reduce \eqref{k9} to the form
\begin{align}\label{k19}
\begin{pmatrix}0&0&\Psi^1_{\ta}&0\\\chi_{p+}(\uq)u^+_{H,j}&\chi_{p-}(\uq)u^-_{H,k}&\Psi^2_{\ta}&\chi_{s,h}(\uq)(i\htau+\hgamma,i\heta)\end{pmatrix}
\end{align}
after row switches.   With \eqref{k5} this completes the proof of
\eqref{ceq16}, and shows that $\beta(\uq)$ is a nonvanishing
multiple of the $(N-k)\times (N-k)$ determinant
$\det\Psi^1_{\ta}(\up,0,0,\ua)$.

To see that $\beta(\uq)\neq 0$, note that are obvious column
operations that we can perform on \eqref{k12} to reduce it to the
form
\begin{align}\label{k17}
\begin{pmatrix}\uw_z(0)&*&*&0\\\uw_{zz}(0)&*&*&0\\0&0&0&\uw_z(0)\cdot\uw_z(0)\end{pmatrix},
\end{align}
and these have no effect on the $2N\times (N-k)$ submatrix
\eqref{k13}. Thus, the same column operations reduce
$\tPsi^1_a(\up,0,0,\ua)$, whose determinant is nonzero by
a-transversality, to
\begin{align}\label{k17z}
\begin{pmatrix}\Psi^1_{\ta}&0\\0&\uw_z(0)\cdot\uw_z(0)\end{pmatrix}.
\end{align}

\textbf{4. Part (b).}    The hypothesis on $\uw$ allows us to define
$\chi$ and write down $D_{Lop}$.  The $u^\pm_{H,j}(\uq,\hzeta,\rho)$
in \eqref{ce13} now extend continuously to $\oS^d_+\times
[0,\rho_0)$, so \eqref{k5} holds with $O_{\hgamma}(\rho^2)$ replaced
by an error $o(\rho)$ that is uniform for $\hzeta\in\oS^d_+$.
Repetition of parts 2 and 3 of this proof yields \eqref{ceqa16}.

\end{proof}

Combining Theorem \ref{K} and Corollary \ref{c40i} we obtain the
immediate corollary:

\begin{cor}\label{k20}
When $d\geq 2$, the low frequency standard Evans condition at $\uq$
implies strong transversality  of $\uw$:
\begin{align}\label{k21}
\begin{split}
&\mathrm{rank}\nabla_{a}\tilde{\Psi}(\up,0,0,\ua)=N+1-k\\
&\mathrm{rank}\nabla_{a,p}\tilde{\Psi}(\up,0,0,\ua)=2N+1.
\end{split}
\end{align}

\end{cor}

\subsection{Low frequency analysis of the modified Evans
function}\label{LFmodified} \textbf{}

 First  we rewrite the transmission problem from
\eqref{c48} (with $f^\pm=0$) as a $2N\times 2N$ first-order system
on $\bR$.   With $V=(v,v_z)$ we obtain:
\begin{align}\label{d1}
\partial_zV=G(z,q,\zeta)V\text{ on }\pm z\geq
0,\;\Gamma(q,\hzeta,\rho)(V_+,V_-,\phi)=0\text{ on }z=0.
\end{align}
Conjugation to HP form using
\begin{align}\label{d3}
V_\pm=T_\pm(z,q,\zeta)\begin{pmatrix} u_{H\pm}\;\;\\
u_{P\pm}\end{pmatrix}\text{ recall \eqref{c6o}}
\end{align}
transforms \eqref{d1} to
\begin{align}\label{d4}
\begin{split}
&\partial_z\begin{pmatrix}u_{H\pm}\\u_{P\pm}\end{pmatrix}=\begin{pmatrix}H_\pm(q,\zeta)&0\\0&P_\pm(q,\zeta)\end{pmatrix}\begin{pmatrix}u_{H\pm}\\u_{P\pm}\end{pmatrix}\\
&\Gamma_{H,P,\cbR}(q,\hzeta,\rho)(u_H,u_P,\phi):=\Gamma_H(q,\zeta)u_H+\Gamma_P(q,\zeta)u_P+\Gamma_{\check{\bfR}}(q,\hzeta,\rho)\phi=0\text{
on }z=0,
\end{split}
\end{align}
where
\begin{align}\label{d5z}
\begin{split}
&\Gamma_{H}(q,\zeta)u_H=\begin{pmatrix}\begin{bmatrix}T_{11}(0,q,\zeta)u_H\\T_{21}(0,q,\zeta)u_H\end{bmatrix}\\(T_{11+}(0,q,\zeta)u_{H+})\cdot\uw_z(0)\end{pmatrix},
\;\Gamma_{P}(q,\zeta)u_P=\begin{pmatrix}\begin{bmatrix}T_{12}(0,q,\zeta)u_P\\T_{22}(0,q,\zeta)u_P\end{bmatrix}\\(T_{12+}(0,q,\zeta)u_{P+})\cdot\uw_z(0)\end{pmatrix},\\
&\text{ and }\Gamma_{\check{\bfR}}(q,\hzeta,\rho)\phi=\begin{pmatrix}[\check{\bfR}(0,q,\hzeta,\rho)]\\
[\check{\bfR}_z(0,q,\hzeta,\rho)]\\0\end{pmatrix}\phi.
\end{split}
\end{align}

\begin{rem}\label{d6}

\textup{For $\Gamma_\cR$ as in \eqref{c15} we have in view of
\eqref{c54}
\begin{align}\label{d8}
\Gamma_\cR(i\htau+\hgamma,\heta)=\begin{pmatrix}[\cR(0,q,i\htau+\hgamma,i\heta)]\\
[\cR_z(0,q,i\htau+\hgamma,i\heta)]\\0\end{pmatrix}=\Gamma_{\check{\bfR}}(q,\hzeta,0).
\end{align}
For $\Gamma_{0,H}(q)$, $\Gamma_{0,P}(q)$ as in \eqref{c14} we have
\begin{align}\label{d9}
\Gamma_{0,H}(q)=\Gamma_H(q,0),\;\;\Gamma_{0,P}(q)=\Gamma_{P}(q,0),
\end{align}
where the operators on the right in \eqref{d9} are as in
\eqref{d4}.}
\end{rem}

The following Lemma is used, for example, in the proofs of
Proposition \ref{d10} and Theorem \ref{d21}.

\begin{lem}[\cite{Me1}, Lemma 6.2.4] \label{d10z}
Consider a subspace $E\subset\bC^D$ with $\dim E=D_+$, and let
$\Gamma$ be a $D_+\times D$ matrix such that $\dim E +\dim \ker
\Gamma =D.$ If
\begin{align}\label{d10y}
|\det(E,\ker\Gamma)|\geq c>0,
\end{align}
then
\begin{align}\label{d10x}
|e|\leq C|\Gamma e|\text{ for all }e\in E,
\end{align}
where $C=c^{-1}|\Gamma^*(\Gamma\Gamma^*)^{-1}|$.

Conversely, if \eqref{d10x} holds, then \eqref{d10y} is satisfied
with $c=(C|\Gamma|)^{-D_+}$.

\end{lem}

\begin{prop}\label{d10}
If $W(z,q)$ satisfies the modified low frequency Evans condition
at $q\in\cC$, then $W(z,q)$ is transversal.
\end{prop}

\begin{proof}

\textbf{1. a-transversality}. For $\rho>0$ let $\bE_-(H(q,\zeta))$
and $\bE_-(P(q,\zeta))$ be as defined in \eqref{ce10}.  Note that
the boundary operator $\Gamma_{H,P,\cbR}(q,\hzeta,\rho)$ in
\eqref{d4} maintains full rank on $\bC^{4N+1}$ for $\rho$ small
since $\Gamma(q,\hzeta,\rho)$ in \eqref{d1} does. Thus, by  Lemma
\ref{d10z} the modified low frequency Evans condition holds if and
only if there exist positive constants $C$ and $\rho_0$ such that
for $\zeta\in\overline{R}^{d+1}_+$ with $0<|\zeta|\leq\rho_0$
\begin{align}\label{d11}
\begin{split}
&|u_H|+|u_P|+|\phi|\leq
C\left|\Gamma_H(q,\zeta)u_H+\Gamma_P(q,\zeta)u_P+\Gamma_{\check{\bfR}}(q,\hzeta,\rho)\phi\right|\\
&\quad\text{ for
}u_H\in\bE_-(H(q,\zeta)),\;\;u_P\in\bE_-(P(q,\zeta)),\;\;\phi\in\bC.
\end{split}
\end{align}

Take $u_H=0$, $\phi=0$ and use the smoothness of
$\bE_\mp(P_\pm(q,\zeta))$ at $\zeta=0$ to conclude
\begin{align}\label{d12}
|u_P|\leq C|\Gamma_P(q,0)u_P|\text{ on }\bE^-(P(q,0)).
\end{align}
In view of Prop. \ref{c16t} and \eqref{d9}, this implies
a-transversality.

\textbf{2. (a,p,s)-transversality.}  Fix $\hzeta=(\htau,\hgamma,0)$
and set $\zeta=\rho\hzeta$.  The estimate \eqref{d11} implies that
for all $h\in\bC^{2N+1}$ and all $0<\rho\leq\rho_0$ there is a
$(u_H(q,\zeta),u_P(q,\zeta),\phi)$ in
\[
\bE_-(H(q,\zeta))\times\bE_-(P(q,\zeta))\times\bC\subset\bC^{4N+1}
\]
such that
\[
\Gamma_{H,P,\cbR}(q,\hzeta,\rho)(u_H,u_P,\phi)=h\text{ and
}|(u_H(q,\zeta),u_P(q,\zeta),\phi)|\leq C|h|.
\]
Thus, letting $\rho\to 0$ and using compactness and continuity, we
obtain an element $(u_H^*,u_P^*,\phi^*)\in
\bC^{2N}\times\bE_-(P(q,0))\times \bC$ such that
\[
\Gamma_{H,P,\cbR}(q,\hzeta,0)(u_H^*,u_P^*,\phi^*)=h.
\]
In view of the choice of $\hzeta$ and Corollary \ref{c300}, this
implies (a,p,s)-transversality.

\end{proof}

%\begin{rem}\label{d13}
%In the boundary layer case Guy was able in Prop. 1.20 of [Feb6] to
%deduce family transversality as well.  Here I don't yet see how to
%get that. Is it true for nonconservative shocks?  Note that in
%Prop. () we show that weak transversality + uniform Lopatinski is
%equivalent to low frequency uniform Evans.

%\end{rem}

\subsubsection{Reduced  modified Evans condition}\label{reducedmod}

The reduced modified Evans function, $D_{red}(q,\hzeta,\rho)$
defined in this section is needed for the block decompositions of
the standard and modified Evans functions.  Its definition requires
only a-transversality; when viscous continuity and strong
transversality hold, $D_{red}$ provides a continuous extension of
$D_{Lop}(q,\hzeta)$ to $\rho>0$.

Suppose $W(z,q)$ is a-transversal.  Then the map
\begin{align}\label{d14}
\Gamma_P(q,0):\bE_-(P(q,0))\to\bF_P(q)
\end{align}
is an isomorphism which extends by continuity to a neighborhood of
$\zeta=0$ (recall the notation \eqref{c20}, \eqref{ce10}). Thus we
get a smooth extension of the decomposition \eqref{c19}:
\begin{align}\label{d15}
\bC^{2N+1}=\bF_{H,\cR}(q)\oplus\bF_P(q,\zeta),\;\;\bF_P(q,\zeta):=\Gamma_P(q,\zeta)\bE_-(P(q,\zeta)).
\end{align}
Denote by $\pi_{H,\cR}(q,\zeta)$ and $\pi_P(q,\zeta)$ the
associated projections, and define the \emph{reduced (transformed)
boundary operator} by
\begin{align}\label{d16}
\hGamma_{red}(q,\hzeta,\rho)(u_H,\phi):=\pi_{H,\cR}(q,\zeta)\left(\Gamma_H(q,\zeta)u_H+\Gamma_{\check{\bfR}}(q,\hzeta,\rho)\phi\right)
\end{align}
and the \emph{reduced (transformed) transmission problem}
\begin{align}\label{d17}
\partial_z u_H-H(q,\zeta)u_H=f_H\text{ on }\pm z\geq
0,\;\;\;\;\hGamma_{red}(q,\hzeta,\rho)(u_H,\phi)=g\text{ on }z=0.
\end{align}
The \emph{reduced modified Evans function} is
\begin{align}\label{d18}
D_{red}(q,\hzeta,\rho)=\det\left(\bE_-(H(q,\zeta))\times
\bC,\ker\hGamma_{red}(q,\hzeta,\rho)\right).
\end{align}

\begin{defn}\label{d18a}
The reduced modified Evans condition at $q$ is satisfied when
there exist positive constants $c$ and $\rho_0$ such that
\begin{align}\label{d19}
|D_{red}(q,\hzeta,\rho)|\geq c\text{ for }\hzeta \in
\overline{S}^d_+,\; 0<\rho\leq\rho_0.
\end{align}

\end{defn}

\begin{rem}\label{d20d}

1.   \textup{Using \eqref{d8} we see that for $\hgamma\geq0$
\begin{align}\label{d20f}
\hGamma_{0,red}(q,\hzeta) \text{ (as in  \eqref{c420})
}=\hGamma_{red}(q,\hzeta,0) \text{ (as in \eqref{d16})}.
\end{align}}

2.  \textup{In general the estimate
\begin{align}\label{d20}
\begin{split}
&\qquad |u_H|+|\phi|\leq C|\hGamma_{red}(q,\hzeta,\rho)(u_H,\phi)|\\
&\text{ for }(u_H,\phi)\in\bE_-(H(q,\zeta))\times \bC\text{ and
}\hzeta \in \overline{S}^d_+,\; 0<\rho\leq\rho_0.
\end{split}
\end{align}
implies \eqref{d19}.  By Lemma \ref{d10z} the converse holds when
the norm of
\begin{align}\label{d20z}
\left(\hGamma_{red}(q,\hzeta,\rho)\hGamma_{red}^*(q,\hzeta,\rho)\right)^{-1}
\end{align}
is uniformly bounded for $\hzeta \in \overline{S}^d_+$ and
$0<\rho\leq\rho_0$.   When $d=1$ this is the case, for example, when
$W(z,q)$ is transversal.   When $d\geq 2$, this holds when $W(z,q)$
is strongly transversal. In each case the transversality hypothesis
implies $\hGamma_{red}(q,\hzeta,\rho)$ has full rank for $\hzeta \in
\overline{S}^d_+,\; 0\leq\rho\leq\rho_0$.}

3.  \textup{For $\check{H}_\pm(q,\hzeta,\rho)$ as in \eqref{ce9},
$\hzeta\in\overline{S}^d_+$, and $\rho>0$ we clearly have
\begin{align}\label{d20e}
D_{red}(q,\hzeta,\rho)=\det\left(\bE_-(\check{H}(q,\hzeta,\rho))\times
\bC,\ker\hGamma_{red}(q,\hzeta,\rho)\right).
\end{align}}

4. \textup{The function $D_{red}(q,\cdot,\cdot)$ extends
continuously to any subset of $\oS^d_+\times [0,\rho_0]$ where
$\bE_-(\check{H}(q,\hzeta,\rho))$ is continuous and
$\hGamma_{red}(q,\hzeta,\rho)$ maintains full rank.  When such a
subset includes a point $(\hzeta,0)$ with $\hzeta\in S^d_+$, we have
\begin{align}\label{d20y}
D_{Lop,m}(q,\hzeta)=D_{red}(q,\hzeta,0).
\end{align}
(recall \eqref{ce9}(a)).}

\textup{Note that $\bE_-(\check{H}_\pm(q,\hzeta,\rho))$, which was
defined on $\overline{S}^d_+\times (0,\rho_0]$, always has a smooth
extension to $S^d_+\times [0,\rho_0]$.}
\end{rem}

\begin{thm}\label{d21}
(a)  If $W(z,q)$ satisfies the  modified low frequency Evans
condition at $q$ (Defn. \ref{c51}), then $W(z,q)$ is transversal and
the reduced modified Evans condition (Defn. \ref{d18a}) holds.

(b)  The converse holds for $d=1$.  When $d\geq 2$, the converse
holds if transversality of $W(z,q)$  is replaced by strong
transversality.

\end{thm}

\begin{proof}
\textbf{1. }By Prop. \ref{d10}  the modified low frequency Evans
condition implies transversality.

Using the splitting \eqref{d15} we see that the estimate
\eqref{d11} is equivalent to (suppress $(q,\zeta)$ dependence):
\begin{align}\label{d22}
\begin{split}
&|u_H|+|u_P|+|\phi|\leq
C\left(|\hGamma_{red}(u_H,\phi)|+|\Gamma_Pu_P+\pi_P(\Gamma_Hu_H+\Gamma_{\check{\bfR}}\phi)|\right)\\
&\quad\text{ for
}u_H\in\bE_-(H(q,\zeta)),\;\;u_P\in\bE_-(P(q,\zeta)),\;\;\phi\in\bC.
\end{split}
\end{align}
Since $\Gamma_P$ is surjective from $\bE_-(P(q,\zeta))$ onto
$\bF_P(q,\zeta)$, for all $u_H\in\bE_-(H(q,\zeta))$ and $\phi\in\bC$
there is $u_P\in\bE_-(P(q,\zeta))$ such that
\begin{align}\notag
\Gamma_Pu_P=-\pi_P(\Gamma_Hu_H+\Gamma_{\check{\bfR}}\phi),
\end{align}
so \eqref{d22} implies \eqref{d20}.

\textbf{2. }Conversely, if the profile  has the stated
transversality property, the estimate \eqref{d12} holds and remains
valid by continuity for $\zeta$ in a neighborhood of zero. Using
Remark \ref{d20d} part 2,  we deduce the estimate \eqref{d20} and
this implies \eqref{d22}. To see this note
\begin{align}\label{d23}
\begin{split}
&|u_P|\leq C|\Gamma_P u_P|\leq |\Gamma_P
u_P+\pi_P(\Gamma_Hu_H+\Gamma_{\check{\bfR}}\phi)|+|\pi_P(\Gamma_Hu_H+\Gamma_{\check{\bfR}}\phi)|\\
&\;\;\;\qquad\leq |\Gamma_P
u_P+\pi_P(\Gamma_Hu_H+\Gamma_{\check{\bfR}}\phi)|+C|\hGamma_{red}(u_H,\phi)|.
\end{split}
\end{align}

\end{proof}

%\begin{rem}\label{d24}
%1.  As in Remark 1.22 of \cite{FF} we note that the vector bundle
%$\bE^-(\check{H}_\pm(q,\hzeta,\rho))$, which was defined on
%$\overline{S}^d_+\times (0,\rho_0]$, has a smooth extension to
%$S^d_+\times [0,\rho_0]$.  Remark \ref{d20d} and \eqref{ce9}(a),
%when the uniform Lopatinski condition holds,
%$D_{red}(q,\hzeta,\rho)$ is continuous on $S^d_+\times [0,\rho_1]$
%for $\rho_1>0$ small, since then $\hGamma_{red}(q,\hzeta,\rho)$
%has full rank at $\rho=0$ and hence for small $\rho$  by
%continuity. Moreover, we have
%\begin{align}\label{d25}
%D_{Lop}(q,\hzeta)=D_{red}(q,\hzeta,0)\text{ for }\hgamma>0.
%\end{align}

%2.   Similarly, strong transversality implies
%$\hGamma_{red}(q,\hzeta,\rho)$ has full rank at $\rho=0$, so again
%$D_{red}(q,\hzeta,\rho)$ is continuous on $S^d_+\times [0,\rho_1]$
%for $\rho_1>0$ small and \eqref{d25} holds.
%\end{rem}

It is now a simple matter to prove a Zumbrun-Serre type result for
the modified Evans function associated to a nonconservative shock.

\begin{thm}\label{d26}
1.  If the modified low frequency Evans condition holds, then the
profile $W(z,q)$ is transversal and the  hyperbolic problem
\eqref{c36z} satisfies the  modified uniform Lopatinski condition.

2. Conversely, when $d=1$ suppose $W(z,q)$ is transversal and that
the hyperbolic problem \eqref{c36z} satisfies the modified uniform
Lopatinski condition.  Then the modified low frequency Evans
condition holds.

 When $d\geq 2$ assume strong
transversality of $W(z,q)$, the modified uniform Lopatinski
condition, and continuous  extendability of the vector bundle
$\bE_-(\check{H}(q,\hzeta,\rho))$ to $\overline{S}^d_+\times
[0,\rho_0]$.  Then the modified low frequency Evans condition holds.

\end{thm}

\begin{proof}
\textbf{1. }Assuming the modified low frequency Evans condition,
from Theorem \ref{d21} and its proof, we deduce transversality of
$W(z,q)$ and the estimate \eqref{d20}. If $\hgamma>0$ every term in
\eqref{d20} is continuous up to $\rho=0$ (recall Remark \ref{d20d},
part 4), so estimate \eqref{d20} implies,
\begin{align}\label{d27}
\begin{split}
&\qquad |u_H|+|\phi|\leq C|\hGamma_{red}(q,\hzeta,0)(u_H,\phi)|\\
&\text{ for }(u_H,\phi)\in\bE_-(\check{H}(q,\hzeta,0))\times
\bC\text{ and }\hzeta \in S^d_+.
\end{split}
\end{align}
From \eqref{ce9}(a) and \eqref{d20f}, we see that \eqref{d27}
implies the modified uniform Lopatinski condition.

\textbf{2. }By Remark \ref{d20d} part 4, continuity  of
$\bE_-(\check{H}_\pm(q,\hzeta,\rho))$ on
 $\overline{S}^d_+\times [0,\rho_0]$ (which is automatic when $d=1$) and the given transversality hypotheses imply that  the reduced modified Evans function $D_{red}(q,\hzeta,\rho)$ has
a continuous extension to $\overline{S}^d_+\times [0,\rho_1]$ for
some $\rho_1>0$.  The modified uniform Lopatinski condition and
\eqref{d20y} imply
\begin{align}\label{d28}
|D_{red}(q,\hzeta,\rho)|\geq c>0
\end{align}
for $\hzeta\in S^d_+$ and $\rho=0$.  By continuity this extends
first to $\hzeta\in\overline{S}^d_+$ and next to $\rho\in
[0,\rho_1]$ for some $\rho_1>0$.   Thus, the reduced modified
Evans condition holds.   An application of Theorem \ref{d21} now
yields the converse.
\end{proof}

\subsection{Block decomposition of the modified Evans
function}\label{blockmod} \textbf{}

A block decomposition was not needed to prove the Zumbrun-Serre
result for the modified Evans function, Theorem \ref{d26}, but such
decompositions are useful for understanding the relation between the
modified and standard Evans functions.  In fact, together with
Theorem \ref{K}, they yield a proof that the low frequency standard
Evans condition implies the low frequency modified Evans condition
(Proposition \ref{j3}), a proof that requires neither continuity of
decaying eigenspaces (Defn. \ref{j1}) nor constant multiplicities of
hyperbolic characteristics. This fact allows our results to be
applied, for example, to viscous MHD, even in the more difficult
case of slow shocks.
%(see Corollary ?).

% and also for
%showing that the low frequency Evans condition implies (full)
%transversality.

In this subsection and the next we obtain block decompositions
assuming strong transversality of the profile $W(z,q)$.   The
modified Evans function is $D_m(q,\hzeta,\rho)$ as in \eqref{c49}.
We'll use the notation
\begin{align}\label{f1}
\Gamma_H(q,\zeta),\; \Gamma_P(q,\zeta),\;
\Gamma_{\check{\bfR}}(q,\hzeta,\rho)
\end{align}
introduced in \eqref{d4}, as well as the notation
\begin{align}\label{f2}
\hGamma_{red}(q,\hzeta,\rho),\;D_{red}(q,\hzeta,\rho),\text{ and
}\bE_-(H(q,\zeta)),\; \bE_-(P(q,\zeta))
\end{align}
as in \eqref{d16}, \eqref{d18}, and \eqref{ce10} respectively.

For small $\rho$ we can use the conjugators $T_\pm(z,q,\zeta)$
\eqref{cee11} to rewrite $D_m$ up to a nonvanishing factor as
\begin{align}\label{f3}
D_m(q,\hzeta,\rho)=\det\left(\bE_-(H(q,\zeta))\times\bE_-(P(q,\zeta))\times\bC,\ker
\Gamma_{H,P,\check{\bfR}}(q,\hzeta,\rho)\right),
\end{align}
where
\begin{align}\label{f4}
\Gamma_{H,P,\check{\bfR}}(q,\hzeta,\rho):\bC^{2N}\times\bC^{2N}\times\bC\to\bC^{2N+1}
\end{align}
is given as in \eqref{d4} by
\begin{align}\label{f5}
\Gamma_{H,P,\check{\bfR}}(q,\hzeta,\rho)(u_H,u_P,\phi):=\Gamma_H(q,\zeta)u_H+
\Gamma_P(q,\zeta)u_P+\Gamma_{\check{\bfR}}(q,\hzeta,\rho)\phi.
\end{align}

Clearly,  the map in \eqref{f4} is surjective, so
\begin{align}\label{f5a}
\dim\ker\Gamma_{H,P,\check{\bfR}}(q,\hzeta,\rho)=2N.
\end{align}

\begin{prop}\label{f6}
(a)  Assume the profile $W(z,q)$ is strongly transversal.  Then
for $\rho_0$ small enough, $0<\rho\leq\rho_0$, and
$\hzeta\in\overline{S}^d_+$, the modified Evans determinant
satisfies
\begin{align}\label{f7}
D_m(q,\hzeta,\rho)=\beta(q,\zeta)D_{red}(q,\hzeta,\rho),
\end{align}
where $\beta(q,\cdot)$ is a nonvanishing smooth function on a
neighborhood of the origin.

(b)  When $d=1$ the same result holds if we just assume that
$W(z,q)$ is transversal.

\end{prop}

\begin{proof}
As in the proof of \cite{Me1}, Theorem 6.4.1 (which gives an
analogous result for Dirichlet boundary layers), the key is to
decompose the kernel of $\Gamma_{H,P,\check{\bfR}}$ as in \eqref{f4}
in a suitable way.

\textbf{1. }Define the map $\Pi_{H,\phi}$ by
\begin{align}\label{f8}
\Pi_{H,\phi}(u_H,u_P,\phi)=(u_H,\phi),
\end{align}
and for $\rho$ small set
\begin{align}\label{f9}
\hat{\cC}(q,\hzeta,\rho):=\Pi_{H,\phi}\left(\ker\Gamma_{H,P,\check{\bfR}}(q,\hzeta,\rho)\cap(\bC^{2N}\times\bE_-(P(q,\zeta))\times\bC)\right).
\end{align}
Observe that
\begin{align}\label{f10}
\hat{\cC}(q,\hzeta,\rho)=\ker\hGamma_{red}(q,\hzeta,\rho),
\end{align}
where
\begin{align}\label{f11}
\hGamma_{red}(q,\hzeta,\rho):\bC^{2N}\times\bC\to\bF_{H,\cR}(q)
\text{ (recall \eqref{d15})}.
\end{align}
Moreover, we have
\begin{align}\label{f12}
\begin{split}
&(u_H,u_P,\phi)\in\ker\Gamma_{H,P,\check{\bfR}}(q,\hzeta,\rho)\cap(\bC^{2N}\times\bE_-(P(q,\zeta))\times\bC)\Leftrightarrow\\
&\qquad(u_H,\phi)\in\hat{\cC}(q,\hzeta,\rho)\text{ and
}u_P=K(q,\hzeta,\rho)(u_H,\phi),
\end{split}
\end{align}
where
\begin{align}\label{f13}
K(q,\hzeta,\rho)(u_H,\phi):=-\Gamma_P(q,\zeta)^{-1}\pi_P(q,\zeta)\left(\Gamma_H(q,\zeta)u_H+
+\Gamma_{\check{\bfR}}(q,\hzeta,\rho)\phi\right).
\end{align}
This gives the parametrization
\begin{align}\label{f14}
\begin{split}
&\bK(q,\hzeta,\rho):=\ker\Gamma_{H,P,\check{\bfR}}(q,\hzeta,\rho)\cap(\bC^{2N}\times\bE_-(P(q,\zeta))\times\bC)=\\
&\qquad\{(u_H,K(q,\hzeta,\rho)(u_H,\phi),\phi):(u_H,\phi)\in\hat{\cC}(q,\hzeta,\rho)\}.
\end{split}
\end{align}

\textbf{2. }(a,p)-transversality implies that the intersection in
\eqref{f12} is transversal for all $\hzeta$ and hence of dimension
$N+1-k$ when $\rho=0$ (use Prop. \ref{c16t}).  Thus,
$\hat{\cC}(q,\hzeta,0)$ has dimension $N+1-k$ and by continuity,
these properties persist for $\rho>0$ small.

When $d=1$, Corollary \ref{c300} and \eqref{c54}  imply the same can
be said about $\hat{\cC}$ when (a,p)-transversality is replaced by
the weaker assumption of (a,p,s)-transversality.

\textbf{3. }Note that
\begin{align}\label{f15}
\bC^{2N}=\bE_-(P(q,\zeta))\oplus\bE_+(P(q,\zeta)),
\end{align}
and consider the map
\begin{align}\label{f16}
\omega(q,\hzeta,\rho):(u_H,u_{P-}+u_{P+},\phi)\to u_{P+}
\end{align}
from
$\ker\Gamma_{H,P,\check{\bfR}}(q,\hzeta,\rho)\subset\bC^{2N}\times\bC^{2N}\times\bC$
to $\bE_+(P(q,\zeta))$.   The kernel of $\omega$ is the subspace of
dimension $N+1-k$ given by \eqref{f14}, so $\omega$ is surjective.
Thus, there is a map
\begin{align}\label{f17}
K'(q,\hzeta,\rho):\bE_+(P(q,\zeta))\to\ker\Gamma_{H,P,\check{\bfR}}(q,\hzeta,\rho)
\end{align}
such that $\omega K'=I$.  Setting
$\bK'(q,\hzeta,\rho):=K'(q,\hzeta,\rho)\bE_+(P(q,\zeta))$ we have
\begin{align}\label{f18}
\ker\Gamma_{H,P,\check{\bfR}}(q,\hzeta,\rho)=\bK(q,\hzeta,\rho)\oplus\bK'(q,\hzeta,\rho).
\end{align}

\textbf{4. }We are now ready to obtain the block decomposition of
$D_m(q,\hzeta,\rho)$.  Suppressing evaluation at $(q,\hzeta,
\rho)$, we choose bases
\begin{align}\label{f19}
\{u_{H,j}\}_{j=1,\dots,N-1+k},\;\{u_{P,k}\}_{k=1,\dots,N+1-k},\;\{(v_{H,l},\phi_l)\}_{l=1,\dots,N+1-k},\;\{w_{P,m}\}_{m=1,\dots,N-1+k}
\end{align}
of $\bE_-(H)$, $\bE_-(P)$, $\hat{\cC}$, and $\bE_+(P)$,
respectively. Using \eqref{f18} and writing
$K'=(K'_H,K'_P,K'_\phi)$, we can compute $D_m$ as the
$(4N+1)\times(4N+1)$ determinant
\begin{align}\label{f20}
D_m(q,\hzeta,\rho)=\det\begin{pmatrix}u_{H,j}&0&0&v_{H,l}&K'_H(w_{P,m})\\0&u_{P,k}&0&K(v_{H,l},\phi_l)&K'_P(w_{P,m})\\0&0&1&\phi_l&K'_\phi(w_{P,m})\end{pmatrix}
\end{align}

The terms $K(v_{H,l},\phi_l)$ lie in the span of the $u_{P,k}$, so
we can eliminate them in the determinant.  Switching rows and
columns, this shows that, up to a sign,
\begin{align}\label{f21}
\begin{split}
&D_m(q,\hzeta,\rho)=\det\begin{pmatrix}u_{H,j}&0&v_{H,l}\\0&1&\phi_l\end{pmatrix}\det\begin{pmatrix}u_{P,k}&K'_P(w_{P,m})\end{pmatrix}\\
&\qquad=D_{red}(q,\hzeta,\rho)\det\begin{pmatrix}u_{P,k}&w_{P,m}\end{pmatrix},
\end{split}
\end{align}
where we've used $\omega K'=I$ for the last equality.  This gives
\eqref{f7} with
\begin{align}\label{f22}
\beta(q,\zeta)=\det\begin{pmatrix}u_{P,k}&w_{P,m}\end{pmatrix}
\end{align}
up to a sign.

\end{proof}

\begin{rem}\label{f22a}
\textup{Note that for  bases $\{u_{H,j}\}_{j=1,\dots,N-1+k}$,
$\{(v_{H,l},\phi_l)\}_{l=1,\dots,N+1-k}$ of $\bE_-(H(q,\zeta))$ and
$\hat{\cal C}(q,\hzeta,\rho)$ as above, we have
\begin{align}\notag
D_{red}(q,\hzeta,\rho)=\det\begin{pmatrix}u_{H,j}&v_{H,l}\end{pmatrix},
\end{align}
so under the assumptions of Proposition \ref{f6}:
\begin{align}\label{f22b}
D_m(q,\hzeta,\rho)=\det\begin{pmatrix}u_{H,j}&v_{H,l}\end{pmatrix}\beta(q,\zeta)
\end{align}
for $\hzeta\in\oS^d_+$ and $0<\rho\leq\rho_0$.}
\end{rem}

Using Remark \ref{d20d} part 4, we obtain the following immediate
corollary generalizing Prop. 3.13 of \cite{GMWZ3} to the
nonconservative case:
\begin{cor}\label{f23}
Assume $W(z,q)$ is strongly transversal and that the vector bundle
$\bE_-(\check{H}(q,\hzeta,\rho))$ has a continuous extension to
$\oS^d_+\times[0,\rho_0]$.  Then, up to a sign,
\begin{align}\label{f24}
D_m(q,\hzeta,0)=\beta(q,0)D_{Lop,m}(q,\hzeta),
\end{align}
where $\beta(q,0)\neq 0$ is given by \eqref{f22}.

When $d=1$, strong tranversality can be replaced by transversality
in the above statement.
\end{cor}

\subsection{Block decomposition of the standard Evans
function}\label{blockstandard} \textbf{}

The standard Evans function $D_s(q,\hzeta,\rho)$ and the standard
low frequency Evans condition were defined in Definition \ref{cc1}.
The main result of this subsection is the following theorem, which
we need in order to show that the low frequency standard Evans
assumption implies the low frequency modified Evans assumption.
%As
%noted before, the absence of any continuity assumptions is important
%for later applications to variable multiplicity problems like
%viscous MHD.
We shall work now with the alternative form $\tD_s(q,\hzeta,\rho)$
defined in Remark \ref{ce5}, and we recall the notations $\Gamma_s$
\eqref{ce2}, $\tilde{\Gamma}_{H,P}$ \eqref{ce6}.

\begin{thm}\label{g3}
  Assume the profile $W(z,q)$ is strongly transversal.  Then, up to a sign,
\begin{align}\label{ga3}
\tD_s(q,\hzeta,\rho)=\rho\beta(q,\zeta)D_{red}(q,\hzeta,\rho)+O(\rho^2),
\end{align}
where $\beta$ is given by \eqref{f22} and, for some $\rho_0>0$,
the error is uniform for $\hzeta\in\oS^d_+$, $0<\rho\leq\rho_0$.
So in particular we have
\begin{align}\label{ga5}
\tD_s(q,\hzeta,\rho)=\rho D_m(q,\hzeta,\rho)+O(\rho^2)
\end{align}
with the same kind of error.

\end{thm}

Observe that \eqref{ga5} follows immediately from \eqref{ga3} and
Proposition \ref{f6}.  Before giving the proof of \eqref{ga3}, we
need some preparation. The kernel of
\begin{align}\label{g6}
\tGamma_{H,P}(q,\zeta):\bC^{2N}\times\bC^{2N}\to\bC^{2N}
\end{align}
has dimension $2N$, but in decomposing the kernel we now have to
contend with the fact that $\tGamma_P(q,\zeta)$ degenerates on a
one-dimensional subspace of $\bE_-(P(q,\zeta))$ as $\zeta\to 0$.

First, recall from \eqref{p4} the extensions of
$(W_z(z,q),W_{zz}(z,q))$ to nonzero frequencies:
\begin{align}\label{g10}
\begin{split}
&(a)\;\cW_\pm(z,q,\zeta)=T_\pm(z,q,\zeta)\begin{pmatrix}0\\e^{zP_\pm(q,\zeta)}\pi_\pm(q,\zeta)c_\pm(q)\end{pmatrix}\\
&(b)\;\Gamma_s(\cW_+,\cW_-)(0,q,0)=0.
\end{split}
\end{align}
For $c_\pm(q)$ as in \eqref{g10} let
\begin{align}\label{g11}
\bE^-_s(P(q,0)):=\mathrm{span}\begin{pmatrix}c_+(q)\\c_-(q)\end{pmatrix}\subset\bC^{2N}\text{
($s$ for singular)}
\end{align}
and let $\bE^-_n(P(q,0))$ ($n$ for nonsingular) be any complementary
subspace satisfying
\begin{align}\label{g12}
\bE_-(P(q,0))=\bE^-_n(P(q,0))\oplus\bE^-_s(P(q,0))\subset\bC^{2N}
\end{align}
with basis
\begin{align}\label{g13}
\begin{pmatrix}d_{+,i}(q)\\d_{-,i}(q)\end{pmatrix},\;i=1,\dots,N-k.
\end{align}
Setting
\begin{align}\label{g14}
\nu_{P,s}(q,\zeta):=\begin{pmatrix}\pi_+(q,\zeta)c_+(q)\\\pi_-(q,\zeta)c_-(q)\end{pmatrix}\text{
and
}\nu_{P,n,i}(q,\zeta):=\begin{pmatrix}\pi_+(q,\zeta)d_{+,i}(q)\\\pi_-(q,\zeta)d_{-,i}(q)\end{pmatrix},i=1,\dots,N-k,
\end{align}
we smoothly extend the decomposition \eqref{g12} to small
$|\zeta|$:
\begin{align}\label{g15}
\bE_-(P(q,\zeta))=\bE^-_n(P(q,\zeta))\oplus\bE^-_s(P(q,\zeta))\subset\bC^{2N},
\end{align}
where
\begin{align}\label{g16}
\bE^-_n(P(q,\zeta)):=\mathrm{span}\left(\nu_{P,n,i}(q,\zeta),\;i=1,\dots,N-k\right)\text{
and
}\bE^-_s(P(q,\zeta)):=\mathrm{span}\left(\nu_{P,s}(q,\zeta)\right).
\end{align}

Assume now that $W(z,q)$ is a-transversal. $\tGamma_P(q,0)$ vanishes
by \eqref{g10}(b) on $\bE^-_s(P(q,0))$, but is nonsingular (by
a-transversality) on $\bE^-_n(P(q,0))$.   So we can write
\begin{align}\label{g17}
\bC^{2N}=\bF_n(q)\oplus\bF_{H,s}(q),
\end{align}
where
\begin{align}\label{g18}
\bF_n:=\tGamma_P(q,0)\bE^-_n(P(q,0))
\end{align}
and $F_{H,s}(q)$ is any complementary ($N+k$ dimensional) subspace.
The decomposition extends smoothly to $|\zeta|$ small:
\begin{align}\label{g19}
\bC^{2N}=\bF_n(q,\zeta)\oplus\bF_{H,s}(q),\qquad\bF_n(q,\zeta):=\tGamma_P(q,\zeta)\bE^-_n(P(q,\zeta)),
\end{align}
and we let $\pi_n(q,\zeta)$, $\pi_{H,s}(q,\zeta)$ be the
associated projections.

We can now define a new reduced boundary operator,
$\Gamma_*(q,\zeta)$, that will serve as a replacement for
$\hGamma_{red}(q,\hzeta,\rho)$ in the proof of Theorem \ref{g3}.
Using the decomposition \eqref{g15} to write (with obvious
notation)
\begin{align}\label{g20a}
u_P=u_{P,n}+u_{P,s},\qquad u_P\in\bE_-(P(q,\zeta)),
\end{align}
we define
\begin{align}\label{g20}
\Gamma_*(q,\zeta):\bC^{2N}\times\bE^-_s(P(q,\zeta))\to\bF_{H,s}(q)\;\;\text{
and }\cC^*(q,\zeta)\subset\bC^{2N}\times\bE^-_s(P(q,\zeta))
\end{align}
by
\begin{align}\notag
\Gamma_*(q,\zeta)(u_H,u_{P,s}):=\pi_{H,s}(q,\zeta)\left(\tGamma_H(q,\zeta)u_H+\tGamma_P(q,\zeta)u_{P,s}\right),\qquad\cC^*(q,\zeta)=\ker\Gamma_*(q,\zeta).
\end{align}
Then, parallel to \eqref{f14}, we have the parametrization
\begin{align}\label{g21}
\begin{split}
&\bK_*(q,\zeta):=\ker\tGamma_{H,P}(q,\zeta)\cap\left(\bC^{2N}\times\bE_-(P(q,\zeta))\right)\\
&\qquad\qquad=\{\left(u_H,K_*(q,\zeta)(u_H,u_{P,s})+u_{P,s}\right):(u_H,u_{P,s})\in\cC^*(q,\zeta)\},
\end{split}
\end{align}
where
\begin{align}\label{g22}
K_*(q,\zeta)(u_H,u_{P,s}):=-\tGamma_P(q,\zeta)^{-1}\pi_n(q,\zeta)\left(\tGamma_H(q,\zeta)u_H+\tGamma_P(q,\zeta)u_{P,s}\right)\in\bE^-_n(P(q,\zeta)).
\end{align}

A key step in the proof of Theorem \ref{g3} will be to set up a
correspondence between $\cC^*(q,\zeta)$ and
$\hat{\cC}(q,\hzeta,\rho)$. With this preparation we can now give
the proof of Theorem \ref{g3}.

\begin{proof}[Proof of Theorem \ref{g3}.]
\textbf{}

First we decompose $\ker\tGamma_{H,P}(q,\zeta)$ in a suitable way.

\textbf{1. } Strong transversality implies that for $\zeta=0$ the
intersection $\bK_*(q,\zeta)$ \eqref{g21} is transversal and thus of
dimension $N+1-k$; hence, $\cC^*(q,0)$ has dimension $N+1-k$. By
continuity these properties persist for $\rho>0$ small.

\textbf{2. }Again using \eqref{f15}, we consider the map
\begin{align}\label{g30}
\omega_*(q,\zeta):(u_H,u_{P-}+u_{P+})\to u_{P+}
\end{align}
from $\ker\tGamma_{H,P}(q,\zeta)\subset\bC^{2N}\times\bC^{2N}$ to
$\bE_+(P(q,\zeta))$.   The kernel of $\omega_*(q,\zeta)$ is the
subspace of $\ker\tGamma_{H,P}(q,\zeta)$ of dimension $N+1-k$ given
by $\bK_*(q,\zeta)$ \eqref{g21}, so $\omega_*$ is surjective. Thus,
there is a map
\begin{align}\label{g31}
K'_*(q,\zeta):\bE_+(P(q,\zeta))\to\ker\tGamma_{H,P}(q,\zeta)
\end{align}
such that $\omega_* K'_*=I$.  Setting
$\bK'_*(q,\zeta):=K'_*(q,\zeta)\bE_+(P(q,\zeta))$ we have
\begin{align}\label{g32}
\ker\tGamma_{H,P}(q,\zeta)=\bK_*(q,\zeta)\oplus\bK'_*(q,\zeta).
\end{align}

\textbf{3. } Suppressing evaluation at $(q,\hzeta, \rho)$ and
recalling \eqref{cee6}, we choose bases
\begin{align}\label{g33}
\{u_{H,j}\}_{j=1,\dots,N-1+k},\;\{u_{P,k}\}_{k=1,\dots,N+1-k},\;\{(\nu_{H,l},\nu_{P,s,l})\}_{l=1,\dots,N+1-k},\;\{w_{P,m}\}_{m=1,\dots,N-1+k}
\end{align}
of $\bE_-(\check{H})$, $\bE_-(P)$, $\cC^*$, and $\bE_+(P)$,
respectively. Using \eqref{g32} and writing
$K'_*=(K'_{*,H},K'_{*,P})$, we can compute $\tD_s$ as the
$4N\times4N$ determinant
\begin{align}\label{g34}
\tD_s(q,\hzeta,\rho)=\det\begin{pmatrix}u_{H,j}&0&\nu_{H,l}&K'_{*,H}(w_{P,m})\\0&u_{P,k}&\nu_{P,s,l}+K_*(\nu_{H,l},\nu_{P,s,l})&K'_{*,P}(w_{P,m})\end{pmatrix}.
\end{align}

The terms $\nu_{P,s,l}+K_*(\nu_{H,l},\nu_{P,s,l})$ lie in the span
of the $u_{P,k}$, so we can eliminate them in the determinant.
Switching rows and columns, this shows that, up to a sign,
\begin{align}\label{g35a}
\tD_s(q,\hzeta,\rho)=\det\begin{pmatrix}u_{H,j}&\nu_{H,l}&0&K'_{*,H}(w_{P,m})\\0&0&u_{P,k}&K'_{*,P}(w_{P,m})\end{pmatrix}=\det\begin{pmatrix}u_{H,j}&\nu_{H,l}\end{pmatrix}\beta(q,\zeta),
\end{align}
where we've used $\omega_* K'_*=I$ for the last equality and
$\beta$ is as in \eqref{f22}.  To finish we will show
\begin{align}\label{g35}
\det\begin{pmatrix}u_{H,j}&\nu_{H,l}\end{pmatrix}(q,\hzeta,\rho)=\rho
D_{red}(q,\hzeta,\rho)+O(\rho^2).
\end{align}

\textbf{4. } To set up a correspondence between $\cC^*$ and
$\hat{\cC}$ it is helpful to write out the boundary conditions
more explicitly.  We have
\begin{align}\label{g37}
\begin{split}
&(u_H,u_{P,s})\in\cC^*(q,\zeta)\Leftrightarrow \exists \;u_{P,n}\in\bE^-_n(P(q,\zeta)) \text{ such that } \\
&\qquad\begin{bmatrix}T_{11}u_H\\T_{21}u_H\end{bmatrix}+\begin{bmatrix}T_{12}u_{P,s}\\T_{22}u_{P,s}\end{bmatrix}+\begin{bmatrix}T_{12}u_{P,n}\\T_{22}u_{P,n}\end{bmatrix}=0.
\end{split}
\end{align}
On the other hand
\begin{align}\label{g38}
\begin{split}
&(u_H,\phi)\in\hat{\cC}(q,\hzeta,\rho)\Leftrightarrow\exists
\;u_P\in\bE_-(P(q,\zeta))\text{ such that }\\
&(a)\;\qquad\begin{bmatrix}T_{11}u_H\\T_{21}u_H\end{bmatrix}+\begin{bmatrix}T_{21}u_P\\T_{22}u_P\end{bmatrix}
+\begin{bmatrix}\phi\cbR(0,q,\hzeta,\rho)\\\phi\cbR_z(0,q,\hzeta,\rho)\end{bmatrix}=0\text{
and }\\
&(b)\;\qquad
(T_{11+}u_{H+})\cdot\uw_z(0)+(T_{12+}u_{P+})\cdot\uw_z(0)=0.
\end{split}
\end{align}

To proceed further we make a more explicit choice of basis of
$\cC^*(q,\zeta)$.   Strong transversality allows us to choose near
$\zeta=0$ a smooth basis of $\cC^*(q,\zeta)$ of the form
\begin{align}\label{g39}
\left\{\begin{pmatrix}\nu_{H,l}(q,\zeta)\\0\end{pmatrix},l=1,\dots,N-k\right\}\cup\left\{\begin{pmatrix}\nu_H(q,\zeta)\\\nu_{P,s}(q,\zeta)\end{pmatrix}\right\},
\end{align}
where $\nu_{P,s}(q,0)=(c_+(q),c_-(q))$.   Moreover, with
\begin{align}\label{g39a}
\nu_{P,n}(q,\zeta):=K_*(q,\zeta)(\nu_H,\nu_{P,s})\text{ and }
c(q):=(c_+(q),c_-(q)),
\end{align}
we have
\begin{align}\label{g40}
\uV(q,\zeta):=\begin{pmatrix}\nu_H(q,\zeta)\\\nu_{P,s}(q,\zeta)+\nu_{P,n}(q,\zeta)\end{pmatrix}\in\bK_*(q,\zeta),\;\begin{pmatrix}\nu_H(q,0)\\\nu_{P,s}(q,0)\end{pmatrix}=\begin{pmatrix}0\\c(q)\end{pmatrix},\;\nu_{P,n}(q,0)=0.
\end{align}

\textbf{5. } Starting with the basis of $\cC^*(q,\zeta)$ given by
\eqref{g39}, we next derive from it a basis for
$\hat{\cC}(q,\hzeta,\rho)$.

Given $(\nu_{H,l},0)$ as in \eqref{g39}, there exists
$\nu_{P,n,l}\in\bE^-_n(P(q,\zeta))$ such that part (a) of
\eqref{g38} holds with $u_H=\nu_{H,l}$, $u_P=\nu_{P,n,l}$, and
$\phi=0$. Using \eqref{g10}, \eqref{g40}  and setting
\begin{align}\label{g41}
\begin{pmatrix}v_{H,l}(q,\zeta)\\v_{P,l}(q,\zeta)\end{pmatrix}:=\begin{pmatrix}\nu_{H,l}(q,\zeta)\\\nu_{P,n,l}(q,\zeta)\end{pmatrix}+\alpha(q,\zeta)\uV(q,\zeta),\;l=1,\dots,N-k,
\end{align}
for an appropriate smooth scalar $\alpha(q,\zeta)$, we see that
both (a) and (b) of \eqref{g38} are satisfied with $u_H=v_{H,l}$,
$u_P=v_{P,l}$, $\phi=0$. Thus,
\begin{align}\label{g42}
(v_{H,l}(q,\zeta),0),\; l=1,\dots,N-k
\end{align}
are linearly independent elements of $\hat{\cC}(q,\hzeta,\rho)$
for $\rho\geq 0$ small.

To obtain the last basis element of $\hat{\cC}(q,\hzeta,0)$ we
apply $\partial_\rho|_{\rho=0}$ to the equation
\begin{align}\label{g43}
\begin{bmatrix}T_{11}\nu_H\\T_{21}\nu_H\end{bmatrix}+\begin{bmatrix}T_{12}\nu_{P,s}\\T_{22}\nu_{P,s}\end{bmatrix}+\begin{bmatrix}T_{12}\nu_{P,n}\\T_{22}\nu_{P,n}\end{bmatrix}=0.
\end{align}
With
\begin{align}\label{g44}
\nu_H^{\sharp}(q,\hzeta):=\partial_\rho|_{\rho=0}\nu_H(q,\zeta),\;\;\nu_{P,n}^{\sharp}(q,\hzeta):=\partial_\rho|_{\rho=0}\nu_{P,n}(q,\zeta)
\end{align}
we obtain using \eqref{g40} and Lemma \ref{p6}
\begin{align}\label{g45}
\begin{bmatrix}T_{11}\nu_H^{\sharp}\\T_{21}\nu_H^{\sharp}\end{bmatrix}-\begin{bmatrix}\cbR(0,q,\hzeta,0)\\\cbR_z(0,q,\hzeta,0)\end{bmatrix}+\begin{bmatrix}T_{12}u_{P,z}\\T_{22}u_{P,z}\end{bmatrix}+\begin{bmatrix}T_{12}\nu_{P,n}^\sharp\\T_{22}\nu_{P,n}^\sharp\end{bmatrix}=0,
\end{align}
where $u_{P,z,\pm}$ was defined in \eqref{p7}.

Observe that since we can use the basis of $\bE^-_n(P(q,\zeta))$
given in \eqref{g16} to write
\begin{align}\label{g46}
\nu_{P,n}(q,\zeta)=\sum^{N-k}_{i=1}c_i(q,\zeta)\nu_{P,n,i}(q,\zeta)
\end{align}
with $c_i(q,0)=0$ for all $i$ (by \eqref{g40}), we have
$\nu_{P,n}^{\sharp}(q,\hzeta)\in\bE^-_n(P(q,0))$.

In view of \eqref{g45}, by setting
\begin{align}\label{g47}
\begin{pmatrix}v_{H,N+1-k}(q,\hzeta)\\v_{P,N+1-k}(q,\hzeta)\end{pmatrix}:=\begin{pmatrix}\nu_{H}^\sharp(q,\hzeta)\\u_{P,z}(q,\hzeta)+\nu_{P,n}^\sharp(q,\hzeta)\end{pmatrix}+\beta(q,\hzeta)\uV(q,0)
\end{align}
for an appropriate scalar $\beta(q,\hzeta)$, we can arrange so that
both (a) and (b) of \eqref{g38} are satisfied at $\rho=0$ with
$u_H=v_{H,N+1-k}$, $u_P=v_{P,N+1-k}$, $\phi=-1$.   Thus,
$(\nu^\#_H(q,\hzeta),-1)^t\in\hat{\cC}(q,\hzeta,0)$.   Choosing a
smooth extension of this element,
\begin{align}\notag
\begin{pmatrix}\nu^\#_{H,e}(q,\hzeta,\rho)\\-1+O(\rho)\end{pmatrix}\in\hat{\cC}(q,\hzeta,\rho),
\end{align}
we see that a basis of $\hat{\cC}(q,\hzeta,\rho)$ for $\rho$ small
is given by
\begin{align}\label{g48}
\left\{\begin{pmatrix}v_{H,l}(q,\zeta)\\0\end{pmatrix},l=1,\dots,N-k\right\}\cup\left\{\begin{pmatrix}\nu_{H,e}^\sharp(q,\hzeta,\rho)\\-1+O(\rho)\end{pmatrix}\right\},
\end{align}
where $v_{H,l}$, $l=1,\dots,N-k$ are as in \eqref{g42}.

\textbf{6. }We can now finish the proof by showing \eqref{g35}. We
work with the  basis \eqref{g39} of $\cC^*(q,\zeta)$ and the
associated basis \eqref{g48} of $\hat{\cC}(q,\hzeta,\rho)$.  From
\eqref{g41} we obtain
\begin{align}\label{g49}
\nu_{H,l'}(q,\zeta)=v_{H,l'}(q,\zeta)+O(\rho),\; l'=1,\dots,N-k,
\end{align}
and from \eqref{g40}, \eqref{g44} we have
\begin{align}\label{g50}
\nu_{H,N+1-k}(q,\zeta):=\nu_H(q,\zeta)=\rho\nu^\sharp_H(q,\hzeta)+O(\rho^2).
\end{align}
 Using \eqref{g49}, \eqref{g50},
and Remark \ref{f22a}, we find
\begin{align}\label{g51}
\begin{split}
&\det\begin{pmatrix}u_{H,j}&\nu_{H,l}\end{pmatrix}(q,\hzeta,\rho)=\det\begin{pmatrix}u_{H,j}&v_{H,l'}(q,\zeta)+O(\rho)&\rho\nu^\sharp_{H,e}(q,\hzeta,\rho)+O(\rho^2)\end{pmatrix}\\
&\qquad\qquad=\rho D_{red}(q,\hzeta,\rho)+O(\rho^2).
\end{split}
\end{align}

\end{proof}

\subsection{Summary of low frequency  results}\label{summary}

%For applications to variable multiplicity problems like viscous
%MHD and its nonconservative analogues,   we need to understand the
%relationship between the low frequency standard and modified Evans
%conditions when no assumptions (or only weak assumptions) are made
%about continuity of the decaying eigenspaces
%$\bE^-(H_0(q,\hzeta))$ \eqref{c39} and $\bE^-(\check{H}(q,\hzeta,
%\rho))$ \eqref{d20e}.    Here we summarize our results on this
%point.

%First note that there are two relevant notions of continuity.

%TO:Does inviscid continuity play a role?

The following theorem, which ties together the results of sections 5
and 6,  is our main low frequency stability result.

\begin{thm}\label{j3}

(a)  The low frequency standard Evans condition \eqref{ce7}
implies the low frequency modified Evans condition \eqref{c52}.
The converse holds when $d=1$.

(b) Assume the profile $W(z,q)$ is strongly transversal.   Then
the low frequency standard Evans condition is equivalent to the
low frequency modified Evans condition.

(c)   Assume the profile $W(z,q)$ is strongly transversal and that
viscous continuity holds. Then for $\rho_0$ small enough,
$0<\rho\leq\rho_0$, and $\hzeta\in\overline{S}^d_+$, we have
\begin{align}\label{g53}
\begin{split}
&(a)\;D_m(q,\hzeta,\rho)=\beta(q,0)D_{Lop}(q,\hzeta)+o(1)\\
&(b)\;\tD_s(q,\hzeta,\rho)=\rho\beta(q,0)D_{Lop}(q,\hzeta)+o(\rho),
\end{split}
\end{align}
where $\beta(q,0)$ is given by \eqref{f22} and the errors are
uniform for $(\hzeta,\rho)\in\oS^d_+\times[0,\rho_0]$. Consequently,
under these assumptions the modified and standard low frequency
Evans conditions are both equivalent to the uniform Lopatinski
condition for the inviscid hyperbolic problem.

(d)  The low frequency standard Evans condition implies
transversality and the uniform Lopatinski condition.  When $d\geq
2$, the  low frequency standard Evans condition implies strong
transversality.

(e)  The low frequency modified Evans condition implies
transversality and the modified uniform Lopatinski condition.

\end{thm}

\begin{proof}
\textbf{1. Parts (b) and (c).}  Part (b) follows immediately from
Theorem \ref{g3}. Part (c) follows from the block decompositions of
the modified and standard Evans functions, Remark \ref{d20d}, part
4, and Proposition \ref{c408}.

\textbf{2. Part (a).}  When $d\geq 2$, it follows from Theorem
\ref{K} and Corollary \ref{c40i} that the low frequency standard
Evans condition implies strong transversality.  So we can apply (b)
to deduce the low frequency modified Evans condition.   When $d=1$,
the low frequency standard Evans condition   implies transversality
and the uniform Lopatinski condition (Theorem \ref{K}). Hence, the
low frequency modified Evans condition follows from part 2 of
Theorem \ref{d26}.

The low frequency modified Evans condition implies the modified
uniform Lopatinski condition and transversality (Thm. \ref{d26}).
When $d=1$ we can apply Prop. \ref{cb3} to deduce that the uniform
Lopatinski condition holds,  so the standard low frequency Evans
condition follows from part (b) of Theorem \ref{K}.

\textbf{3. Part (d)}.  This is contained in the statement of Theorem
\ref{K} and Corollary \ref{k20}.

\textbf{4. Part (e)}.  This is contained in Theorem \ref{d26}.

\end{proof}

\begin{rem}\label{sim}
\textup{1. The standard Evans condition is the one that is easier to
verify numerically \cite{B,HZ} or analytically \cite{PZ,FS}, while
the modified Evans condition is the one needed for the rigorous
study of small viscosity limits.   Thus the implication (a) in
Theorem \ref{j3} is especially important.   The implication (d)
allows us to construct  curved inviscid $\cC_\cB$ shocks near a
given planar shock when the low frequency standard Evans condition
is satisfied (as long as there exists a $K$-family of smooth
inviscid symmetrizers). The next section shows that implication (d)
also permits the construction of high order approximate solutions to
the nonlinear parabolic transmission problem \ref{h1}.}

\textup{2. When viscous continuity holds, the implication in part
(a) of Theorem \ref{j3} can be proved without the use of the block
decompositions given by Proposition \ref{f6} and Theorem \ref{g3}.
In fact, it follows immediately from Theorem \ref{K}, Corollary
\ref{k20}, and Theorem \ref{d26}.   This type of argument was used
in \cite{GMWZ3} in the conservative setting. }

\end{rem}

\section{Approximate viscous shocks}\label{approximate}

 In this section we construct high
order approximate solutions to the nonlinear small viscosity
transmission problem:
\begin{align}\label{h1}
\begin{split}
&(a)\;\cE(u,d\psi):=\sum^{d-1}_{j=0}A_j(u)\partial_ju+\cA_d(u,d\psi)\partial_du-\epsilon\sum^d_{j=1}(\partial_j-\partial_j\psi\partial_d)^2u=0\text{ on }\pm x_d\geq 0\\
&(b)\;[u]=0,\;[\partial_du]=0,\\
&(c)\;\partial_t\psi-\epsilon\triangle_y\psi+\ell(t,y)\cdot
u=\partial_t\psi^0-\epsilon\triangle_y\psi^0+\ell(t,y)\cdot
\cU^0(t,y,0,0),
\end{split}
\end{align}
where $t=x_0$, $y=(x_1,\dots,x_{d-1})$.  Here we suppose that we are
given an inviscid $\cC_\cB$-shock $(U^0,\psi^0)$ on
$[-T_0,T_0]\times \bR^d_{y,x_d}$ satisfying
%with $(U^0_+(t,y,0),U^0_-(t,y,0),d\psi^0(t,y))$ near
%$\uq=(\up_+,\up_-,0,0)$ for all $(t,y)$
\begin{align}\label{h2}
\begin{split}
&\sum^{d-1}_{j=0}A_j(U^0)\partial_jU^0+\cA_d(U^0,d\psi^0)\partial_dU^0=0\text{ on }\pm x_d\geq 0\\
&(U^0_+(t,y,0),U^0_-(t,y,0),d\psi^0(t,y))\in\cC_\cB
\end{split}
\end{align}
%and the $\chi$-Lopatinski condition.
In \eqref{h1} we have set
\begin{align}\label{h3}
\ell(t,y):=W_z(0,q(t,y)),\text{ where
}q(t,y):=(U^0_+(t,y,0),U^0_-(t,y,0),d\psi^0(t,y))\text{ and }
\end{align}
$W(z,q)$ is the profile associated to $q$.     Moreover, we suppose
that the inviscid shock satisfies the hypotheses of Theorem
\ref{i27}, with the modifications that the standard uniform Evans
condition is replaced by the standard low frequency Evans condition
(Definition \ref{cc1}), and the existence of a $K$-family of smooth
viscous symmetrizers is replaced by the existence of a $K$-family of
smooth \emph{inviscid} symmetrizers.

The function $\cU^0$ appearing in \eqref{h1}(c) can be written as
$\cU^0(t,y,x_d,\frac{x_d}{\epsilon})$ where
\begin{align}\label{h4}
\cU^0(t,y,x_d,z):=U^0(t,y,x_d)+\left(W(z,q(t,y))-U^0(t,y,0)\right),
\end{align}

% Since $q(t,y)$ is near $\uq$ for all $(t,y)$, there is a
%constant $c$ such that
%\begin{align}\label{h5}
%\uw_z(0)\cdot\partial_z\cU^0(t,y,0,0)\geq c>0 \text{ for all
%}(t,y).
%\end{align}

We seek an approximate solution $(u^a,\psi^a)$ of the form
\begin{align}\label{h6}
\psi^a=\psi^0(t,y)+\epsilon\psi^1(t,y)
+\cdots+\epsilon^{M}\psi^{M}(t,y),
\end{align}
\begin{align}\label{h7}
u^a=\left(\mathcal{U}^0(t,y,x_d,z)+\epsilon\mathcal{U}^1(t,y,x_d,z)+\cdots+\epsilon^{M}\mathcal{U}^{M}(t,y,x_d,z)\right)|_{z=\frac{x_d}{\epsilon}},
\end{align}
where
\begin{align}\label{h7a}
\mathcal{U}^j(t,y,x_d,z)=U^j(t,y,x_d)+V^j(t,y,z).
\end{align}
Here $V^0$ is already determined and is given by
\begin{align}\label{h8}
V^0(t,y,z)=W(z,q(t,y))-U^0(t,y,0).
\end{align}
The $V^j_\pm(t,y,z)$ are boundary layer profiles constructed to be
exponentially decreasing to $0$ as $z\to \pm\infty$.    For the
moment we just assume enough regularity so that all the operations
involved in the construction make sense.   A precise statement is
given in Prop. \ref{esoln}.

\subsection{Profile equations}\label{profile}

We substitute \eqref{h6}, \eqref{h7} into \eqref{h1} and write the
result as
\begin{align}\label{h9}
\sum^M_{-1}\epsilon^j\mathcal{F}^j(t,y,x_d,z)|_{z=\frac{x_d}{\epsilon}}+\epsilon^M
R^{\epsilon,M}(t,y,x_d),
\end{align}
where we separate $\cF^j$ into slow and fast parts
\begin{align}\label{h10}
\mathcal{F}^j(t,y,x_d,z)=F^j(t,y,x_d)+G^j(t,y,z),
\end{align}
and the $G^j$ decrease exponentially to $0$ as $z\to\pm\infty$.

The interior profile equations are obtained by setting the $F^j,G^j$
equal to zero.  In the following expressions for $G^j(t,y,z)$, the
functions $U^j(t,y,x_d)$ and their derivatives are evaluated at
$(t,y,0)$. Let $\cL_0(z,q,\partial_z)$ and $\cL_{0,1}(z,q)$ be the
operators defined in \eqref{c5} and set
\begin{align}\label{h11}
L_0v:=\sum^{d-1}_{j=0}A_j(U^0)\partial_jv+\cA_d(U^0,d\psi^0)\partial_dv.
\end{align}
We have
\begin{align}\label{h12}
\begin{split}
&F^{-1}(t,y,x_d)=0\\
&G^{-1}(t,y,z)=-(1+|d\psi^0|^2)\partial_z^2\cU^0+\cA_d(\cU^0,d\psi^0)\partial_z\cU^0,
\end{split}
\end{align}
\begin{align}\label{h13}
\begin{split}
&F^0(t,y,x_d)=L_0U^0\\
&G^0(t,y,z)=\cL_0(z,q(t,y),\partial_z)\cU^1-\cL_{0,1}(z,q(t,y))d\psi^1-Q^0(t,y,z),
\end{split}
\end{align}
where $Q^0$ decays exponentially as $z\to\pm\infty$ and depends
only on $(U^0,V^0,d\psi^0)$. For $j\geq 1$ we have
\begin{align}\label{h14}
\begin{split}
&F^j(t,y,x_d)=L_0U^j-P^{j-1}(t,y,x_d)\\
&G^j(t,y,z)=\cL_0(z,q(t,y),\partial_z)\cU^{j+1}-\cL_{0,1}(z,q(t,y))d\psi^{j+1}-Q^j(t,y,z),
\end{split}
\end{align}
where $Q^j$ decays exponentially as $z\to\pm\infty$ and $P^j$,
$Q^j$ depend only on $(U^k,V^k,d\psi^k)$ for $k\leq j$.

Similarly, we obtain the boundary profile equations in which
$(t,y,x_d,z)$ is evaluated at $(t,y,0,0)$:
\begin{align}\label{h15}
\begin{split}
&(a)\;[\cU^0]=0\\
&(b)\;[\cU^0_z]=0\\
&(c)\;\partial_t\psi^0-\ell(t,y)\cdot\cU^0=\partial_t\psi^0-\ell(t,y)\cdot\cU^0,
\end{split}
\end{align}
\begin{align}\label{h16}
\begin{split}
&(a)\;[\cU^1]=0\\
&(b)\;[\cU^1_z]=-[\partial_{x_d}U^0]\\
&(c)\;\partial_t\psi^1-\triangle_y\psi^0+\ell(t,y)\cdot\cU^1=-\triangle_y\psi^0,
\end{split}
\end{align}
and for $j\geq 2$
\begin{align}\label{h17}
\begin{split}
&(a)\;[\cU^j]=0\\
&(b)\;[\cU^j_z]=-[\partial_{x_d}U^{j-1}]\\
&(c)\;\partial_t\psi^j-\triangle_y\psi^{j-1}+\ell(t,y)\cdot\cU^j=0.
\end{split}
\end{align}

\subsection{Solution of the profile
equations}\label{solutionprofile}

The solution of the profile equations given below assumes strong
transversality and the uniform Lopatinski condition, as well as the
existence of a $K$-family of smooth inviscid symmetrizers. Recall
from Theorem \ref{j3} and Corollary \ref{k20} that when $d\geq 2$,
the first two conditions both follow from the low frequency standard
Evans condition.  Strong tranversality can be replaced by
transversality as explained in Remark \ref{h25}.

\textbf{1. } The interior equations $G^{-1}=0$ and $F^0=0$ and the
boundary equations \eqref{h15} are satisfied because of our
assumptions about $U^0$, $\psi^0$ and $W(z,q)$.

\textbf{2. Construction of $(\cU^1,U^1,\psi^1)$.}     We construct
the function $\cU^1(t,y,x_d,z)$ as in \eqref{h7a} from the
equations $G^0=0$, $F^1=0$, and the boundary equations
\eqref{h16}.  $\cU^1$ will be a sum of three parts
\begin{align}\label{h18}
\begin{split}
&\cU^1(t,y,x_d,z)=\cU^1_a+\cU^1_b+\cU^1_c,\text{ where }\\
&\cU^1_k(t,y,x_d,z)=U^1_k(t,y,x_d)+V^1_k(t,y,z),\;k=a,b,c,
\end{split}
\end{align}
where we suppress $\pm$ subscripts.

First use the exponential decay of $Q^0$ to find an exponentially
decaying solution $V^1_a(t,y,z)$ to
\begin{align}\label{h19}
\begin{split}
&\cL_0(z,q(t,y),\partial_z)V^1_a=Q^0(t,y,z)\text{ on }\pm z\geq 0\\
&V^1_a\to 0\text{ as }z\to\pm\infty,
\end{split}
\end{align}
and define $U^1_a(t,y,x_d)\equiv 0$.   This is the same type of
problem as \eqref{c6n}, which we solved by conjugating the
corresponding first order system with $T_\pm(z,q,0)$.

%Such solutions are constructed, for example, in Prop. 5.5.6 of
%\cite{Me1}.

Next, for $\cU^1_a$ fixed as above,  use (a,p)-transversality
(recall Prop. \ref{c10a}) to solve for $\cU^1_b(t,y,0,z)$ in
\begin{align}\label{h20}
\begin{split}
&\cL_0(z,q(t,y),\partial_z)\cU^1_b=0\text{ on }\pm z\geq 0\\
&[\cU^1_a+\cU^1_b]=0\\
&[\partial_z\cU^1_a+\partial_z\cU^1_b]=-[\partial_{x_d}U^0]\\
&\ell(t,y)\cdot(\cU^1_a+\cU^1_b)=0.
\end{split}
\end{align}
Using \eqref{c10} we see that $\cU^1_b$ has limits as
$z\to\pm\infty$. Define
\begin{align}\label{h21}
\begin{split}
&U^1_b(t,y,0):=\lim_{z\to\pm\infty}\cU^1_b(t,y,0,z), \\
&V^1_b(t,y,z):=\cU^1_b(t,y,0,z)-U^1_b(t,y,0),
\end{split}
\end{align}
and let $U^1_b(t,y,x_d)$ be any smooth extension of $U^1_b(t,y,0)$
with compact support in $x_d$.

Finally, for an appropriate choice of $(U^1_c(t,y,0),\psi^1)$ we
need $\cU^1_c(t,y,0,x_d)$ to satisfy
\begin{align}\label{h22}
\begin{split}
&\cL_0(z,q(t,y),\partial_z)\cU^1_c=\cL_{0,1}(z,q(t,y))d\psi^1\\
&[\cU^1_
c]=0,\;[\partial_z\cU^1_c]=0,\;\partial_t\psi^1+\ell(t,y)\cdot\cU^1_c=0\\
&\lim_{z\to\pm\infty}\cU^1_c(t,y,0,z)=U^1_{c\pm}(t,y,0).
\end{split}
\end{align}
According to the characterization of $T_q\cC_\cB$ given in
Proposition \ref{c30}, this is possible if and only if
$(U^1_c(t,y,0),d\psi^1)\in T_{q(t,y)}\cC_\cB$.   Thus, we first
solve for $(U^1_c(t,y,x_d),\psi^1)$ satisfying the linearized
inviscid problem
\begin{align}\label{h23}
\begin{split}
&L_0U^1_c=P^0-L_0U^1_b\text{ on }\pm x_d\geq 0\\
&(U^1_c(t,y,0),d\psi^1(t,y))\in T_{q(t,y)}\cC_\cB
%\Leftrightarrow\chi'(q(t,y))(U^1_c(t,y,0),d\psi^1(t,y))=0.
\end{split}
\end{align}

This problem requires an initial condition in order to be
well-posed.  The right side in the interior equation of \eqref{h23}
is initially defined just for $t\in [-T_0,T_0]$. With a $C^\infty$
cutoff that is identically one in $t\geq -T_0/2$, we can modify the
right side to be zero in $t\leq -T_0+\delta$, say. Requiring
$(U^1_c,\psi^1)$ to be identically zero in $t\leq -T_0+\delta$, we
thereby obtain a problem for $(U^1_c,\psi^1)$ that is forward
well-posed since $(U^0,\psi^0)$ satisfies the uniform Lopatinski
condition  and $\chi_p$ has full rank $N+k$ (see the proof of
Theorem \ref{i8p} in section \ref{existence}). Thus, there exists a
solution to \eqref{h23} on $[-\frac{T_0}{2},T_0]$. This allows us to
obtain $\cU^1_c(t,y,0,z)$ satisfying \eqref{h22} and to define
\begin{align}\label{h24}
V^1_c(t,y,z):=\cU^1_c(t,y,0,z)-U^1_c(t,y,0).
\end{align}

By construction the functions  $(\cU^1,U^1,\psi^1)$ satisfy the
equations $G^0=0$, $F^1=0$, and the boundary conditions
\eqref{h16}.

\textbf{3. Contruction of $(\cU^j,U^j,\Psi^j)$, $j\geq 2$.}  In the
same way, for $j\geq 2$ we use the equations $G^{j-1}=0$, $F^j=0$,
and the boundary conditions \eqref{h17} to determine the functions
$(\cU^j,U^j,\psi^j)$.

\begin{rem}\label{h25}

%  \textup{In order to carry out the profile construction, first we need
%transversality to define $\chi$ and the manifold $\cC$.  We also
%need the $\chi$-uniform Lopatinski condition in order to prove
%estimates for \eqref{h23}.   Together these conditions are roughly
%equivalent to the low frequency standard Evans condition (see
%Theorem \ref{K} for a precise statement).   (Also need inviscid
%symmetrizers.)}

%\textup{The additional condition of (a,p)-transversality was used to
%solve for $\cU^1_b(t,y,0,z)$ and, indirectly, to solve for
%$U^1_c(t,y,x_d)$ through the cited linear inviscid theory of Mokrane
%\cite{Mo} and Coulombel \cite{Cou}, where it was used to prove
%existence: specifically, to define a suitable adjoint problem.  In
%the case $d\ge 2$  we have shown (a,p)-transversality to be a
%consequence of transversality and the $\chi$-uniform Lopatinski
%condition (Prop. \ref{c40i}) and not an additional hypothesis, thus
%resolving a point left open in the treatment of \cite{Mo, Cou}.}

   \textup{It turns out that the  profile construction can be
carried out with only slight changes (e.g., a nonhomogeneous
boundary condition in \eqref{h23}) assuming just transversality, the
uniform Lopatinski condition, and the existence of a $K$-family of
smooth inviscid symmetrizers. In solving for $\cU^1_b(t,y,0,z)$ one
uses the formulation of (a,p,s)-transversality  given in Corollary
\ref{c300}. The inviscid problem \eqref{h23} can be solved in $d\geq
1$ without explicit reliance on (a,p)-transversality using the
approach described in \cite{GMWZ6}.}
% In particular, we can construct
%profiles by this method when $d=1$.}

\end{rem}

In the next Proposition we formulate a precise statement summarizing
the construction of this section.  The regularity assertions in the
Proposition are justified as in \cite{GMWZ4}, Prop. 5.7, except that
now we use the estimates of \cite{Mo,Cou} in place of those of
\cite{Ma1}. Regularity is expressed in terms of the following
spaces:
\begin{defn}\label{espaces}
1.  Let $H^s$ be the set of functions $U(t,y,x_d)$ on
$[-T_0,T_0]\times\bR^d$ such that the restrictions $U_\pm$ belong to
$H^s([-T_0,T_0]\times\overline{\bR}^d_\pm)$.

2. Let $\tH^s$ be the set of functions $V(t,y,z)$ on
$[-T_0,T_0]\times\bR^{d-1}\times\bR$ such that the restrictions
$V_\pm$ belong to $C^\infty(\overline{\bR}_\pm,H^s(t,y))$ and
satisfy
\begin{align}\label{es1}
|\partial_z^kV(t,y,z)|_{H^s(t,y)}\leq C_{k,s}e^{-\delta|z|}\text{
for all }k
\end{align}
for some $\delta>0$.

\end{defn}

\begin{prop}[Approximate solutions]\label{esoln}
For given integers $m\geq 0$ and $M\geq 1$ let
\begin{align}\label{es3}
s_0>m+\frac{7}{2}+2M+\frac{d+1}{2}.
\end{align}
Suppose the given inviscid shock $(U^0,\psi^0)$  has the properties
assumed in Theorem \ref{i27}, but now require only the standard low
frequency Evans condition and the existence of a $K$-family of
smooth inviscid symmetrizers.  Assume $U^0\in H^{s_0}$,
$U^0_\pm(t,y,0)\in H^{s_0}(t,y)$, and $\psi^0(t,y)\in
H^{s_0+1}(t,y)$.  Then one can construct $(u^a,\psi^a)$ as in
\eqref{h6}, \eqref{h7}
\begin{align}\label{es4}
\psi^a=\psi^0(t,y)+\epsilon\psi^1(t,y)
+\cdots+\epsilon^{M}\psi^{M}(t,y),
\end{align}
\begin{align}\label{es5}
u^a=\left(\mathcal{U}^0(t,y,x_d,z)+\epsilon\mathcal{U}^1(t,y,x_d,z)+\cdots+\epsilon^{M}\mathcal{U}^{M}(t,y,x_d,z)\right)|_{z=\frac{x_d}{\epsilon}},
\end{align}
Let $\cE$ denote the operator on the left side of \eqref{h1}(a).
The approximate solution $(u^a,\psi^a)$ satisfies
\begin{align}\label{es6}
\begin{split}
&\mathcal{E}(u^a,\psi^a)=\epsilon^M R^M(t,y,x_d)\text{ on }[-\frac{T_0}{2},T_0]\times\overline{\bR}^d_\pm\\
&[u^a]=0; \;\;[\partial_du^a]=0\text{ on }x_d=0\\
&\begin{aligned}\partial_t\psi^a- &
\epsilon\triangle_{y}\psi^a+\ell(t,y)\cdot
u^a\\
&=\partial_t\psi^0-\epsilon\triangle_{y}\psi^0+\ell(t,y)\cdot\mathcal{U}^{0}(t,y,0,0)\text{
on }x_d=0.\end{aligned}
\end{split}
\end{align}
We have
\begin{align}\label{es7}
\begin{split}
&U^j(t,y,x_d)\in H^{s_0-2j},\;\psi^j(t,y)\in H^{s_0-2j+1}(t,y)\\
&V^j(t,y,z)\in\tH^{s_0-2j},
%&r^M(t,y)\in H^{s_0-2M-\frac{3}{2}}(t,y),
\end{split}
\end{align}
and $R^M(t,y,x_d)$ satisfies
\begin{align}\label{es8}
\begin{split}
&(a)\;|(\partial_t,\partial_y,\epsilon\partial_{x_d})^\alpha
R^M|_{L^2(t,y,x_d)}\leq
C_\alpha\text{ for }|\alpha|\leq m+\frac{d+1}{2}\\
&(b)\;|(\partial_t,\partial_y,\epsilon\partial_{x_d})^\alpha
R^M|_{L^\infty(t,y,x_d)}\leq C_\alpha\text{ for }|\alpha|\leq m.
\end{split}
\end{align}
\end{prop}

\begin{rem}\label{es100}
\textup{Observe that Proposition \ref{esoln} can be applied to give
high order approximate solutions even for slow MHD shocks.}

\end{rem}

\section{Existence of nonclassical inviscid and viscous
shocks}\label{existence}

\textbf{\qquad}In this section we present the proofs of Theorems
\ref{i8p} and \ref{i27}.   The discussion will be brief since the
proofs mainly involve combining results proved in previous sections
with results proved elsewhere.

\begin{proof}[Proof of Theorem \ref{i8p}]\;\;  We will show how the proofs given
in \cite{Mo,Cou,Ma2}  can be adapted to our case. We discuss the
case $d\geq 2$.  The case $d=1$ can be treated using results of
\cite{LY}, for example.

\textbf{Part 1. } Let $\chi$ be a local defining function for the
shock manifold $\cC$, and suppose $q\in\cC$ is a point where the
uniform Lopatinski condition holds.  The existence of a smooth
$K$-family of inviscid symmetrizers implies inviscid continuity, so
from the uniform Lopatinski condition at $q$ we deduce that there
exists a $c>0$ such that
\begin{align}\label{m1}
\begin{split}
&|\chi'_p(q)u+\chi'_s(q)(i\htau+\hgamma)\psi+\chi'_h(q)i\heta\psi|\geq
c|u,\psi|\\
&\qquad \text{ for all }u=(u_+,u_-)\in\bE_-(H_0(q,\hzeta)),\;
\psi\in\bC,\; \hzeta\in\oS^d_+.
\end{split}
\end{align}
This implies
\begin{align}\label{m2}
b(q,\hzeta):=\chi'_s(q)(i\htau+\hgamma)+\chi'_h(q)i\heta\neq 0\text{
for all }\hzeta\in\oS^d_+,
\end{align}
so we may define $\Pi(q,\hzeta):\bC^{N+k}\to\bC^{N+k}$, the smooth
orthogonal projector onto $b(q,\hzeta)^\perp$.   Setting
\begin{align}\label{m3}
M(q)u:=\chi'_p(q)u
\end{align}
and applying $\Pi(q,\hzeta)$ to the boundary condition in
\eqref{c37z} we obtain the projected boundary condition (whose
analogue in \cite{Cou} is equation (17a)):
\begin{align}\label{m4}
\Pi(q,\hzeta)M(q)u=0.
\end{align}
Moreover, we deduce from \eqref{m1}  the corresponding uniform
Lopatinski condition for the problem \eqref{c37z}(a) with the new
boundary condition \eqref{m4}:
\begin{align}\label{m5}
|\Pi(q,\hzeta)M(q)u|\geq c|u|\text{ for all
}u\in\bE_-(H_0(q,\hzeta)).
\end{align}
Recall from Proposition \ref{i41} that
\begin{align}\label{m6}
M(q)=\chi'_p(q) \text{ has full rank }N+k
\end{align}
as a consequence of the uniform Lopatinski condition.  The analogue
of \eqref{m6} in the conservative case was stated as ``Assumption 4"
in \cite{Cou}, and was used there to construct a suitable adjoint
boundary condition as a step in proving existence of solutions to
the linearized shock problem by a duality argument.   We are now in
a position to complete the proof of the first part of Theorem
\ref{i8p} by repeating arguments used in \cite{Mo,Cou} to establish
Theorem 5.2 of that paper.

\textbf{Part 2. }First we must define and construct \emph{shock
front initial data compatible to order $s-1$}.  Such data is
constructed in Proposition 2.2 of \cite{Ma2} in the conservative Lax
case. To carry out a similar construction here we need to use local
defining functions, and in order to patch together locally defined
initial data, we must check that compatibility is independent of the
choice of local defining function.

One obtains local corner compatibility conditions as in \cite{Ma2}
by supposing $(u_+,u_-,\psi)$ is a smooth solution to the nonlinear
problem
\begin{align}\label{m7}
\begin{split}
&(a)\;\sum^{d-1}_{j=0}A_j(u_\pm)\partial_ju_\pm+\cA_d(u_\pm,d\psi)\partial_du_\pm=0\text{
in
}\pm x\geq 0\\
&(b)\;\chi(u_+,u_-,\partial_t\psi,\partial_y\psi)=0\text{ on }x=0,\\
&(c)\;u_\pm(0,y,x_d)=u^0_\pm(y,x_d),\;\partial_t\psi(0,y)=\sigma(y),\;\psi(0,y)=\psi^0(y),
\end{split}
\end{align}
computing relations between $\partial_t^{k+1}\psi(0,y)$ and
$\partial_t^ku_\pm(0,y,0)$ by differentiating \eqref{m7}(b), and
using \eqref{m7}(a) to express time derivatives of $u$ in terms of
space derivatives of $u$ (and derivatives of $\psi$ involving
$\partial_t^j$, $j\leq k$).   The relations at $t=0$, $x_d=0$ have
the form
\begin{align}\label{m8}
\partial_t^{k+1}\psi\;\chi'_s+\chi'_{p_+}\;\left(-A_0(u^0_+)^{-1}\cA_d\right)^k\partial_{x_d}^ku^0_++\chi'_{p_-}\;\left(-A_0(u^0_-)^{-1}\cA_d\right)^k\partial_{x_d}^ku^0_-=I_k,
\end{align}
where derivatives of $\chi$ are evaluated at
$(u^0(y,0),\sigma(y),\partial_y\psi^0(y))$, and $I_k$ is an
explicitly computable function involving no $t$ (resp. $x_d$)
derivatives of $\psi$ (resp. $u_\pm$) of order greater than $k$
(resp. $k-1$).

\begin{defn}\label{sfd}
A function $(u^0_+(y,x_d),u^0_-(y,x_d),\sigma(y),\psi^0(y))$
satisfying \eqref{i8i} is said to determine compatible shock front
initial data to order $s-1$ when there exists a function $\psi(t,y)$
such that
\begin{align}\label{m10}
\psi(0,y)=\psi^0(y),\;\partial_t\psi(0,y)=\sigma(y)
\end{align}
and the relations \eqref{m8} hold for $1\leq k\leq s-1$.
\end{defn}

It is not immediately clear that the notion of compatibility is
well-defined.  To see that it is, note that if $\chi_1$ and $\chi_2$
are two local defining functions for $\cC$ near $q$, we must have
\begin{align}\label{m11}
\chi_2(q)=B(q)\chi_1(q)
\end{align}
for some invertible $(N+k)\times (N+k)$ matrix $B(q)$.  Using
\eqref{m11} and a straightforward induction on $k$, one shows that
if $\psi(t,y)$ and $u^0_\pm(y,x_d)$ satisfy the relations \eqref{m8}
for $1\leq k\leq s-1$ with $\chi=\chi_1$,  then these same functions
satisfy the relations with $\chi=\chi_2$.

Majda's construction of compatible data relies on the observation
that the uniform Lopatinski condition implies \emph{one dimensional
stability} (Lemma 2.1 of \cite{Ma2}).   In our context this is the
statement that for $\hzeta'=(\htau',\hgamma',\heta'):=(0,1,0)$, the
map
\begin{align}\label{m12}
\hGamma_\chi(q,\hzeta'):\bE_-(H_0(q,\hzeta'))\times\bC\to \bR^{N+k}
\end{align}
is invertible, which follows immediately from the uniform Lopatinski
condition \eqref{m1}.   This allows us to carry out Majda's
construction of compatible initial data for nonclassical shocks. The
rest of the proof of part (2) of Theorem \ref{i8p} can now be
completed by following \cite{Mo,Cou}.

\end{proof}

\begin{proof}[Proof of Theorem \ref{i27}]
The theorem is proved by constructing exact solutions of the
nonlinear transmission problem \eqref{h1} which are close to an
approximate solution $(u^a,\psi^a)$ as constructed in Proposition
\ref{esoln}.  The key step is to obtain good estimates for solutions
$(v,\phi)$ to the transmission problem obtained by linearizing
\eqref{h1} with respect to both $u$ and $\psi$ at $(u^a,\psi^a)$:
\begin{align}\label{m13}
\begin{split}
&\cE_u'(u^a,\psi^a)v+\cE_\psi'(u^a,\psi^a)\phi=f\text{ on
}[-T_0,T_0]\times\obR^d_\pm\\
&[v]=0,\;[\partial_{d}v]=0,\;\partial_t\phi-\eps\triangle_y\phi+\ell(t,y)\cdot
v=0\text{ on }x_d=0,\\
&v=0,\;\phi=0,\;f=0\text{ in }t<-T_0/3.
\end{split}
\end{align}

The desired estimate is stated (in the conservative, Lax, constant
multiplicity case) in Theorem 7.2 of \cite{GMWZ3}. We can use the
same argument to prove the identical estimate in our context
provided the linearized problem \eqref{m13} satisfies the modified
uniform Evans condition.  The modified uniform Evans condition then
allows us to construct a Kreiss-type symmetrizer for \eqref{m13}
from the $K$-family of smooth viscous symmetrizers by taking $K$
large enough.

By Theorem \ref{j3}(a) of this paper, the low frequency standard
Evans condition implies the low frequency modified Evans condition.
The fact that nonvanishing of $D_s(q,\hzeta,\rho)$ for $\rho>0$
implies nonvanishing of $D_m(q,\hzeta,\rho)$ for $\rho>0$ can easily
be proved just as in Proposition 2.16 of \cite{GMWZ3}.  Observe that
 since $(u^a,d\psi^a)$ differs from $(W(\frac{x_d}{\eps},q(t,y)),d\psi^0)$ by an error which is small in $L^\infty$ for $\eps$ and $|x_d|$ small (recall \eqref{h4}-\eqref{h7}), we need to deduce nonvanishing of
$D_m(q,\hzeta,\rho)$ for $q$ near $\cC_\cB$ but not in $\cC_\cB$;
but this follows from viscous continuity.

With the linear estimate of Theorem 7.2 of \cite{GMWZ3} in hand, the
proof is completed by the same fixed point argument used to prove
Theorem 7.7 of \cite{GMWZ3}.

\end{proof}

\section{Appendix A: Extension to real viscosity}\label{extension}

In this appendix we describe the changes needed to treat real
viscosities.   Consider the $N\times N$ viscous system on
$\bR^{d+1}$ given by \eqref{i9}.

\textbf{1. Structural assumptions. }   The assumptions made here in
order to define and treat the case of real, or partially parabolic,
viscosities are all satisfied, for example, by the compressible
Navier-Stokes and viscous MHD equations.    A detailed discussion of
the assumptions is given in \cite{GMWZ4}.

We again make Assumptions \ref{i1a} and \ref{i3}, but now  add the
following block form requirement:
\begin{ass}\label{w1}
Possibly after a change of variables $u$ and multiplying the system
on the left by an invertible constant coefficient matrix, there is
$N'\in\{1,\dots,N\}$ and there are coordinates
$u=(u^1,u^2)\in\bR^{N-N'}\times\bR^{N'}$ such that
\begin{align}\label{w2}
A_0(u)=\begin{pmatrix}A^{11}_0&0\\A^{21}_0&A^{22}_0\end{pmatrix},\,\;B_{j,k}=\begin{pmatrix}0&0\\0&B^{22}_{j,k}\end{pmatrix}.
\end{align}
\end{ass}

Assumption \ref{i10} should now be replaced by
\begin{ass}\label{w3}
(H2')(Partial parabolicity.)\;The $B^{22}_{j,k}$ are $C^\infty$
functions on $\cU^*$ valued in $\bR^{N'\times N'}$.  There is $c>0$
such that for all $u\in\cU^*$ and $\xi\in\bR^d$ the eigenvalues of
$\oB^{22}(u,\xi)=\sum^d_{j,k=1}\oB^{22}_{j,k}(u)$ satisfy
$\Re\mu\geq c|\xi|^2$.

(H3')(Strict dissipativity.)\;There is $c>0$ such that for all
$u\in\cU$ and $\xi\in\bR^d$ the eigenvalues $\mu$ of
$i\oA(u,\xi)+\oB(u,\xi)$ satisfy
\begin{align}\label{w4}
\Re\mu\geq \frac{c|\xi|^2}{1+|\xi|^2}.
\end{align}

\end{ass}

\begin{ass} For the low
and medium frequency analysis it is enough to assume:

$(H4)_{\ell,m}$\label{w5} \;For all $u\in\cU^*$ and all
$\xi\in\bR^d\setminus 0$,
$\oA^{11}(u,\xi):=\sum^d_{j=1}\oA^{11}_j(u)\xi_j$ has only real
eigenvalues.

For the high frequency analysis we must strengthen this to:

$(H4)_h$\; For all $u\in\cU^*$ and $\xi\in\bR^d\setminus 0$ the
eigenvalues of $\oA^{11}(u,\xi)$ are real and semisimple with
constant multiplicities.
\end{ass}

The profile equation in the case of a real viscosity $\cB(u)$ can
still be written as \eqref{i14}, and profiles $W(z,q)$ are defined
and associated to points $q$ of a shock  manifold $\cC_\cB$ as
before. We discuss the local construction of  $\cC_\cB$ for real
viscosities $\cB(u)$ below.  One difference is that now we restrict
the undercompressive index $k$ to satisfy
\begin{align}\label{w6}
0\leq k\leq N'-1.
\end{align}
In addition we now add the following assumption:
\begin{ass}\label{w7}
Suppose we are given a shock manifold $\cC_\cB$.   For the low and
medium frequency analysis it is enough to suppose:

$(H5)_{\ell,m}$\;\;For any planar shock $q=(p_+,p_-,s,h)\in\cC_\cB$
with normal direction $\nu=\nu(s,h)=(-s,-h,1)$,
\begin{align}\label{w8}
\det\left(\sum^d_{j=0}A_j^{11}(W(z,q))\nu_j\right)\neq 0 \text{ for
all } z\in\bR\cup\{\pm\infty\}.
\end{align}
For the high frequency analysis we must strengthen this to

$(H5)_h$\;\;For any $q\in\cC_\cB$ and $z\in\bR\cup\{\pm\infty\}$ the
polynomial in $\xi$
\begin{align}\label{w9}
\det\left(\sum^d_{j=0}A_j^{11}(W(z,q))\xi_j\right)
\end{align}
is hyperbolic in the direction $\nu(s,h)$.

\end{ass}

\begin{rem}\label{w10}

\textup{1.  With $\cA_d(u,s,h)=\sum^d_{j=0}A_j(u)\nu_j$ as before,
let $\overline{\cA}_d=(A_0)^{-1}\cA_d$.  By Assumptions \ref{w5} and
\ref{w7} the eigenvalues of $\overline{\cA}^{11}_d(W(z,q),s,h)$ are
real and nonzero.  Let $N^1_+$ be the number of positive eigenvalues
of $\overline{\cA}^{11}_d(W(z,q),s,h)$.  By connectedness this
number is independent of $q\in\cC_\cB$, $z\in\bR\cup\{\pm\infty\}$.}

\textup{2.  If instead of a shock manifold we are given just a
single transversal profile $\uw(z)$ associated to a planar shock
$\uq$, then we make Assumption \ref{w7} with $W(z,q)$ replaced by
$\uw(z)$. As explained below this will allow us to construct a shock
manifold near $\uq$.}
\end{rem}

\textbf{2. Construction of $\phi_\pm(z,p_\pm,s,h,a_\pm)$. }With
$w=(w^1,w^2)$ and $w^3:=\partial_zw^2$ the profile transmission
problem equivalent to \eqref{i14} can be written
\begin{align}\label{w11}
\begin{split}
&\partial_z w^1=-(\cA_d^{11})^{-1}\cA_d^{12}w^3\\
&\partial_zw^2=w^3\\
&\partial_z(\cB^{22}_{d,d}
w^2)=\left(\cA^{22}_d-\cA^{21}_d(\cA^{11}_d)^{-1}\cA^{12}_d\right)w^3\\
&[w]=0,\;[w^3]=0\text{ on }z=0,
\end{split}
\end{align}
where the matrices are evaluated at $(w_\pm(z),s,h)$.   For
$q=(p_+,p_-,s,h)$ one again looks for profiles $w=W(z,q)$, with
endstates $p_\pm$ and satisfying \eqref{w11}, near a given profile
$\uw=W(z,\uq)$ where $\uq=(\up_+,\up_-,0,0)$.   In place of
\eqref{b14} we consider
\begin{align}\label{w12}
\begin{pmatrix}w^1\\w^2\\w^3\end{pmatrix}'=\begin{pmatrix}-(\cA_d^{11})^{-1}\cA_d^{12}w^3\\w^3\\G_dw^3\end{pmatrix},
\end{align}
where
\begin{align}\label{w13}
G_d:=(\cB^{22}_{d,d})^{-1}\left(\cA^{22}_d-\cA^{21}_d(\cA^{11}_d)^{-1}\cA^{12}_d\right)
\end{align}
and the matrices are now evaluated at $(p_\pm,s,h)$.  The interior
problems \eqref{w11} in $\pm z\geq 0$ are solved by considering them
as perturbations, quadratic in $(w_\pm-p_\pm,w^3_\pm)$, of
\eqref{w12}.
 $G_d$ is clearly nonsingular and, in fact, strict dissipativity
implies $G_d(p_\pm,s,h)$ has no purely imaginary eigenvalues
(\cite{GMWZ4}, Lemma 3.39).

\begin{defn}\label{w14}
1.  Let $r_-$ (resp. $\ell_+$) denote the number of eigenvalues
$\mu$ of $G_d(p_+,s,h)$ (resp. $G_d(p_-,s,h)$) with $\Re\mu<0$
(resp. $\Re\mu>0$).

2. Parallel to \eqref{b14g}, \eqref{b15} we define invariant
subspaces $\bE_\mp(G_d(p_\pm,s,h))\subset\bR^{N'}$ of dimensions
$r_-$ (resp. $\ell_+$) with  corresponding projections
$\Pi_\mp(p_\pm,s,h):\bR^{N'}\to\bE_\mp(p_\pm,s,h)$, and fix
isomorphisms linear in $a_\pm\in\bE_\mp(G_d(\up_\pm,0,0))$:
\begin{align}\label{w15}
\alpha_\pm(p_\pm,s,h;a_\pm):\bE_\mp(G_d(\up_\pm,0,0))\to\bE_\mp(G_d(p_\pm,s,h))
\end{align}
\end{defn}

Functions $\phi_\pm(z,p_\pm,s,h,a_\pm)$ as in \eqref{bn2} satisfying
the obvious analogue of Proposition \ref{b17a} are now constructed
just as before.  The function $\tPsi(p,s,h,a)$ corresponding to
\eqref{b28} is given by
\begin{align}\label{w16}
\tPsi(p,s,h,a):=\begin{pmatrix}\phi_+(0,\cdot)-\phi_-(0,\cdot)\\\phi^2_{+,z}(0,\cdot)-\phi^2_{-,z}(0,\cdot)\\s+\phi^2_+(0,\cdot)\cdot\uw^2_z(0)-\uw^2(0)\cdot\uw^2_z(0)\end{pmatrix}(p,s,h,a)\in\bR^{N+N'+1}.
\end{align}

\textbf{3. Linearization and HP form. }The rescaled transmission
problem corresponding to \eqref{i42} is now
\begin{align}\label{w17}
\begin{split}
&\sum^{d-1}_{j=0}A_j(u)\partial_ju+\cA_d(u,d\psi)\partial_z
u-\sum^d_{j,k=1}D_j(B_{j,k}(u)D_ku)=0\text{ on }\pm z\geq 0\\
&[u]=0,\;[\partial_z u^2]=0\text{ on }z=0,
\end{split}
\end{align}
where $D_j:=\partial_j-(\partial_j\psi)\partial_z$ for
$j=1,\dots,d-1$ and $D_d=\partial_z$.   Assuming we have an exact
solution of \eqref{w17} given by a profile $W(z,q)$ and front
$\psi=st+hy$, consider the partially and fully linearized
(Fourier-Laplace transformed) transmission problems
\begin{align}\label{w18}
\begin{split}
&\cL(z,q,\zeta,\partial_z)u=f\text{ on }\pm z\geq 0\\
&[u]=0, [u^2_z]=0\text{ on }z=0,
\end{split}
\end{align}
and
\begin{align}\label{w19}
\begin{split}
&\cL(z,q,\zeta,\partial_z)u-\psi\cL_1(z,q,\zeta)=f\\
&[u]=0,\;[u^2_z]=0,\; c_0(\zeta)\psi+\uw^2_z(0)\cdot u^2_+=0\text{
on }z=0,
\end{split}
\end{align}
where $\cL(z,q,\zeta,\partial_z)$ and $\cL_1(z,q,\zeta)$ are written
out explicitly in \cite{GMWZ4}, equations (3.14) and (2.65)  (in the
latter case replace $\partial_zf_j(W)$ by $A_j(W)\partial_zW$,
$j=0,\dots,d-1$).

With $U=(u,u^2_z)$ we rewrite \eqref{w18} as a first-order
$(N+N')\times(N+N')$ transmission problem
\begin{align}\label{w20}
\begin{split}
&\partial_z U-G(z,q,\zeta)U=F\\
&[U]=0\text{ on }z=0,
\end{split}
\end{align}
where the components of $G$ are written out in \cite{GMWZ4},
equation (3.36).    As in \eqref{c6i}, \eqref{ce4} for $|\zeta|$
small we conjugate \eqref{w20} to HP form
\begin{align}\label{w21}
\begin{split}
&\partial_z\begin{pmatrix}u_{H\pm}\\u_{P\pm}\end{pmatrix}=\begin{pmatrix}H_\pm(q,\zeta)&0\\0&P_\pm(q,\zeta)\end{pmatrix}\begin{pmatrix}u_{H\pm}\\u_{P\pm}\end{pmatrix}+\tilde{F}\\
&\tGamma_H(q,\zeta)u_H+\tGamma_P(q,\zeta)u_P=0\text{ on }z=0,
\end{split}
\end{align}
where, as before,
\begin{align}\label{w22}
H_\pm(q,\zeta)=-\cA_d(p_\pm,s,h)^{-1}\left(A_0(p_\pm)(i\tau+\gamma)+\sum^{d-1}_{j=1}A_j(p_\pm)i\eta_j\right)+O(\rho^2),
\end{align}
but now the lower right $N'\times N'$ block is
\begin{align}\label{w23}
P_\pm(q,\zeta)=G_d(p_\pm,s,h)+O(\rho)
\end{align}
for $G_d$ as in \eqref{w13}.

\begin{prop}[Relations between indices]\label{w17a}
Consider the indices $R_-$, $L_+$ (as in Assumption \ref{i3}),
$r_-$, $\ell_+$ (Defn. \ref{w14}),  $N^1_+$ (Remark \ref{w10}), and
the undercompressive index $k$. We have
\begin{align}\label{w24}
\begin{split}
&(a)\;N'+N^1_+=r_-+(N-R_-)\\
&(b)\;N-N^1_+=\ell_++(N-L_+)\\
&(c)\;r_-+\ell_+-N'=R_-+L_+-N=1-k.
\end{split}
\end{align}
\end{prop}

\begin{proof}
Part (c) follows by adding (a) and (b).  To prove (a), let
\begin{align}\label{w25}
G_+(q,\zeta)=\lim_{z\to+\infty}G(z,q,\zeta).
\end{align}
Strict dissipativity (H3') implies that for $\zeta\neq 0$,
$G_+(q,\zeta)$ has no eigenvalues on the imaginary axis.  For $\rho$
small, we can count the eigenvalues $\mu$ with $\Re\mu<0$ using
\eqref{w21} (recall \eqref{c6f}, \eqref{c6g}), and the number is
clearly $r_-+(N-R_-)$.  For $\rho$ large one can show as in
\cite{GMWZ4}, Lemma 3.38 that this number is $N'+N^1_+$.  This
implies (a), and (b) is proved similarly.
\end{proof}

\textbf{4. Transversality and the manifold $\cC_\cB$. }With $\tPsi$
as defined in \eqref{w16}, we have now in place of \eqref{b31}
\begin{align}\label{w26}
\begin{split}
&\mathrm{rank} \nabla_a\tPsi(\up,0,0,\ua)=r_-+\ell_+=N'+1-k\\
&\mathrm{rank} \nabla_{a,p}\tPsi(\up,0,0,\ua)=N+N'+1,
\end{split}
\end{align}
as the conditions that define a-transverality and
(a,p)-transversality respectively.    To define
(a,p,s)-transversality we just replace $\nabla_{a,p}$ by
$\nabla_{a,p,s}$ in \eqref{w26}.  The notions of
\emph{transversality} and \emph{strong transversality} may now be
defined just as before (Defn. \ref{c11}).  Moreover, the proof of
Prop. \ref{b30} can now be repeated to yield an $N+d-k$ dimensional
shock manifold $\cC_\cB$ near $\uq$ when transversality holds. In
place of Prop. \ref{bz1} we now have
\begin{align}\label{w27}
\begin{split}
&\dim\cS_+=N+r_-,\;\;\dim\cS_-=N+\ell_+,\;\;\dim\cS=2N+(N'+1-k)\\
&\dim\cS^0_+=r_-,\;\;\dim\cS^0_-=\ell_+.
\end{split}
\end{align}
\\

Sections 3-6, and Appendix B now extend in a straightforward way to
the case of real viscosity.  In many cases only a change of index is
needed (e.g., $N'$ in place of $N$), or $u_z$ should be replaced by
$u^2_z$ in a boundary condition.   We'll briefly indicate some of
these changes below.\\

\textbf{5. Reduced transmission conditions. }In place of \eqref{c19}
we have
\begin{align}\label{w28}
\bC^{N+N'+1}=\bF_{H,\cR}(q)\oplus\bF_P(q),
\end{align}
where
\begin{align}\label{w29}
\dim\bF_P(q)=N'+1-k,\;\;\dim\bF_{H,\cR}(q)=N+k.
\end{align}
We still have
\begin{align}\label{w30}
\Gamma_{0,red}(q)(u_H,\dos,\doh):=\pi_{H,\cR}(q)\left(\Gamma_{0,H}(q)u_H+\Gamma_\cR(\dos,\doh)\right),
\end{align}
but, of course, $\cR_z$ is replaced by $\cR^2_z$ in the definition
\eqref{c15} of $\Gamma_\cR$.\\

\textbf{6. Stability determinants. }Using Proposition \ref{w17a} we
see that the spaces in \eqref{ce10} now have dimensions
\begin{align}\label{w31}
\begin{split}
&\dim\bE_-(H(q,\zeta))=(N-R_-)+(N-L_+)=N+k-1\text{ as before},\\
&\dim\bE_-(P(q,\zeta))=r_-+\ell_+=N'+1-k,
\end{split}
\end{align}
so $\bD_s(q,\hzeta,\rho)$ \eqref{ce15} is an $(N+N')\times(N+N')$
determinant. Similarly, $D_s$ \eqref{cc2} and $\tD_s$ \eqref{cee6}
are of size $2(N+N')\times 2(N+N')$.    The modified Evans functions
$\tD$ \eqref{c47a} and $D_m$ \eqref{c49} are of size
$(2(N+N')+1)\times (2(N+N')+1)$.

With $\tPsi$  as in \eqref{w16}, the proof of the nonconservative
Zumbrun-Serre works without any substantial change.  Now, for
example, the block $\Psi^1_{\ta}$ in \eqref{k18} is of size
$(N'-k)\times(N'-k)$.

The spaces $\hat{\cC}$ and $\bE_+(P)$ that enter into the block
decomposition of $D_m$ (see \eqref{f19},\eqref{f20}) are now
respectively of dimensions $N+1-k$ (as before) and $N'-1+k$.  Thus,
$\beta(q,\zeta)$ \eqref{f22} is given by a $2N'\times 2N'$
determinant.

The space $\cC^*$ that enters into the definition of $\tD_s$ (see
\eqref{g33}, \eqref{g34}) has dimension $N+1-k$ as before.

Theorem \ref{j3} summarizing the low frequency results can be
repeated verbatim in the case of real viscosity.\\

\textbf{7. Viscous $\cC_\cB$-shocks. }The construction of
approximate viscous shocks in section \ref{approximate} can be
repeated with no significant changes.  Theorem \ref{i27}, giving the
existence of a family of viscous $\cC_\cB$-shocks converging to a
given inviscid $\cC_\cB$-shock as $\eps\to 0$, remains true as
stated when we substitute the new structural Assumptions \ref{w1},
\ref{w3}, \ref{w5},  and \ref{w7}.  The estimates and iteration
scheme of \cite{GMWZ4} (in particular, the high frequency analysis,
which is the main difference between the fully and partially
parabolic cases) can be repeated without change to complete the
proof of Theorem \ref{i27} for real viscosities.   We recall that
the high frequency analysis of \cite{GMWZ4} requires both hypotheses
$(H4)_h$ and $(H5)_h$.

\section{Appendix B: Uniqueness of $\cC_\cB$}\label{uniqueness}

%\subsection{Uniqueness of $\cC_\cB$}

$\;\;\;$Some choices were made in the definition of $\cC_\cB$, and
it is not immediately clear that $\cC_\cB$ is independent of these
choices. We show this now using properties of the functions
$\Phi(z,p,s,h,a)$ arising from translation invariance of the profile
equation.

For example, recall the choice of $\uz$, which determines $\ua$, in
\eqref{b24}.  By the argument of Prop. \ref{b17a}, part (b),
$-\uz\geq 0$ should be large enough so that for all $z_0\geq-\uz$
\begin{align}\label{e1}
\begin{split}
&|\partial_z\uw|_{L^1(|z|\geq z_0)}\leq
R,\;|\partial_z\uw|_{L^\infty(|z|\geq z_0)}\leq R\\
&|\Pi_\mp(\up_\pm,0,0)\partial_z\uw(\pm z_0)|\leq r,
\end{split}
\end{align}
where $R$ and $r$ are the constants in \eqref{b17b}.

\begin{prop}\label{e2} The manifold $\cC_\cB$ defined in Proposition
\ref{b30} is, in a sufficiently small neighborhood of $(\up,0,0)$,
independent of the choice of $\uz$ as above.
\end{prop}

There is also some freedom in the choice of the third component of
$\tilde{\Psi}(p,s,h,a)$ in \eqref{b28}; for example, the $s$ in
the third component could be replaced by any smooth function
$g(s,h)$ such that $g(0,0)=0$.

\begin{prop}\label{e3} The manifold $\cC_\cB$ defined in Proposition
\ref{b30} is, in a sufficiently small neighborhood of $(\up,0,0)$,
independent of the choice of $g(s,h)$ as above replacing $s$ in the
third component of $\tilde{\Psi}(p,s,h,a)$.  More generally, the
same manifold $\cC_\cB$ near $(\up,0,0)$ is obtained for any choice
of the third component of $\tilde{\Psi}$, so long as the rank
conditions assumed in Prop. \ref{b30}  hold and
$\tilde{\Psi}(\up,0,0,\ua)=0$.
\end{prop}

We present the proof of Proposition \ref{e3}.  The proof of
Proposition \ref{e2} is similar but easier.

\begin{proof}[Proof of Proposition \ref{e3}]
\textbf{1. }At first we take $g(s,h)=s$ as in the definition
\eqref{b28} of $\tilde{\Psi}$.

The proof of Prop. \ref{b30} provides us with functions
\begin{align}\label{e5}
p_+(p_\beta,s,h), p_-(p_\beta,s,h), a(p_\beta,s,h)
\end{align}
such that $\cC_\cB$ consists of points of the form
\begin{align}\label{e6}
(p_+(p_\beta,s,h), p_-(p_\beta,s,h),s,h),
\end{align}
where $(p_\beta,s,h)$ varies in a neighborhood $\cN_1$ determined
by the implicit function theorem.   Define
\begin{align}\label{e7}
\cC_{\tilde{\Psi}}=\{(p_+(p_\beta,s,h), p_-(p_\beta,s,h),s,h,
a(p_\beta,s,h):(p_\beta,s,h)\in\cN_1\},
\end{align}
which is precisely the zero set near $(\up,0,0,\ua)$ for the
equation
\begin{align}\label{e8}
\tilde{\Psi}(p,s,h,a)=0.
\end{align}

 Similarly, for $\Psi$ as in \eqref{b25} and using
part 3 of Remark \ref{b56y}, we can use the implicit function
theorem to define the zero set near $(\up,0,0,\ua)$ for the equation
\begin{align}\label{e9}
\Psi(p,s,h,a)=0,
\end{align}
which we write
\begin{align}\label{e10}
\cC_{\Psi}=\{(p_+(p_\gamma,s,h,a_1), p_-(p_\gamma,s,h,a_1),s,h,
a_1,a'(p_\gamma,s,h,a_1):(p_\gamma,s,h,a_1)\in\cN_2\}.
\end{align}
Here we set $a=(a_1,a')$ and take $a_i$ in Remark \ref{b56y} to be
$a_1$.

Let $\pi$ be the projection
\begin{align}\label{e11}
\pi(p_+,p_-,s,h,a)=(p_+,p_-,s,h).
\end{align}
It suffices to show
\begin{align}\label{e12}
\pi\cC_{\tilde{\Psi}}=\pi\cC_{\Psi}\text{ near }(\up,0,0).
\end{align}

\textbf{2. }We clearly have
\begin{align}\label{e13}
\cC_{\tilde{\Psi}}\subset\cC_\Psi\text{ near }(\up,0,0,\ua),\text{
and thus }\pi\cC_{\tilde{\Psi}}\subset\pi\cC_{\Psi}\text{ near
}(\up,0,0),
\end{align}
simply because $\Psi$ gives the first two components of
$\tilde{\Psi}$.

\textbf{3. }A particular point
\begin{align}\label{e14a}
P=(p_+(p_\gamma,s,h,a_1), p_-(p_\gamma,s,h,a_1),s,h,
a_1,a'(p_\gamma,s,h,a_1)
\end{align}
of $\cC_\Psi$ has the property that the function of $z$ given by
(recall \eqref{bn2})
\begin{align}\label{e14}
f_P(z):=\phi(z,p_+(p_\gamma,s,h,a_1), p_-(p_\gamma,s,h,a_1),s,h,
a_1,a'(p_\gamma,s,h,a_1))
\end{align}
satisfies at $z=0$ the first two parts of the boundary condition
\eqref{b27}, but not necessarily the third.  In particular it
defines a smooth solution of \eqref{i14} on $\bR$.

        For $|g(s,h)|=|s|$ small enough, the
function
\begin{align}\label{e16}
f_P(z+r)=\phi(z+r,p_+(p_\gamma,s,h,a_1),
p_-(p_\gamma,s,h,a_1),s,h, a_1,a'(p_\gamma,s,h,a_1))
\end{align}
will satisfy the third boundary condition at $z=0$ for some
translate $r$ (here we suppress the dependence of $r$ on the other
parameters).  Thus,  by translation invariance of \eqref{i14},
$f_P(z+r)$ is a smooth solution of \eqref{i14} on $\bR$ and
satisfies the third boundary condition at $z=0$.  In other words, it
satisfies \eqref{t1}(a), \eqref{b13}, and all three boundary
conditions at $z=0$.

Now, for $-\uz$ initially chosen large enough, the argument of
\eqref{b17a} yields an $a^*$ near $\ua$ such that
\begin{align}\label{e14b}
f_P(z+r)=\phi(z,p_+(p_\gamma,s,h,a_1), p_-(p_\gamma,s,h,a_1),s,h,
a^*).
\end{align}
But then
\begin{align}\label{e15a}
P^*=(p_+(p_\gamma,s,h,a_1), p_-(p_\gamma,s,h,a_1),s,h, a^*)
\end{align}
is a point near $(\up,0,0,\ua)$ satisfying $\tilde{\Psi}=0$.

Hence $P^*\in \cC_{\tilde{\Psi}}$, and this implies that for $P$
as in \eqref{e14a}
\begin{align}\label{e16a}
\pi P\in\pi\cC_{\tilde{\Psi}}.
\end{align}
Thus, $\pi\cC_\Psi\subset\pi\cC_{\tilde{\Psi}}$ near $(\up,0,0)$,
so we obtain \eqref{e12}.

\textbf{4. }The same argument applies if the third boundary
condition (and hence $\tilde{\Psi}$) is redefined by replacing $s$
with any smooth function $g(s,h)$ satisfying $g(0,0)=0$.   Indeed,
the argument works for any choice of the third component of
$\tilde{\Psi}$, so long as the rank conditions \eqref{b31} hold
and $\tilde{\Psi}(\up,0,0,\ua)=0$.

\end{proof}

The following description of $\cC_\cB$ in terms of $\Psi$ provided
by the above proof is worth emphasizing:

\begin{cor}\label{e20}
Assume $\uw$ is transversal (Defn. \ref{c11}) and let $\cC_\cB$ be
the manifold defined in Prop. \ref{b30}.  For $\cC_\Psi$ as in
\eqref{e10} and the projection $\pi$ as in \eqref{c11} we have
\begin{align}
\cC_\cB=\pi\cC_\Psi\text{ near }(\up,0,0).
\end{align}
\end{cor}

\end{document}